\documentclass[svgnames, oneside, 11pt]{amsart}

\usepackage[english]{babel}
\usepackage[T1]{fontenc}

\usepackage{color}
\usepackage{hyperref}
\hypersetup{colorlinks,citecolor=Black,linkcolor=Black,urlcolor=Blue}

\usepackage{tabularx, ifthen}
\usepackage{amssymb} \usepackage{amsfonts} \usepackage{amsmath}
\usepackage{amsthm} \usepackage{epsfig, subfig} \usepackage{mathabx}
\usepackage{ amscd, amsxtra, latexsym,pb-diagram,pb-xy}
\usepackage[all]{xy}
\usepackage{caption}
\usepackage{enumitem}
\usepackage{comment}

\usepackage{geometry}

\usepackage{epigraph}

\usepackage{tikz}\usetikzlibrary{cd}

\newtheorem{lemma}{Lemma}[section]
\newtheorem{thm}[lemma]{Theorem}
\newtheorem{prop}[lemma]{Proposition}
\newtheorem{cor}[lemma]{Corollary}

\newtheorem{thmintro}{Theorem}

\newtheorem{claim}[lemma]{Claim}

\theoremstyle{definition}
\newtheorem{defn}[lemma]{Definition}

\newtheorem{quest}[lemma]{Question}

\newtheorem{property}[lemma]{Property}
\newtheorem{rem}[lemma]{Remark}
\newtheorem{remintro}{Remark}
\newtheorem{ass}[lemma]{Assumption}
\newtheorem{notation}[lemma]{Notation}

\newenvironment{claimproof}{\begin{proof}}{\end{proof}}

\definecolor{darkgreen}{cmyk}{1,0,1,.2}

\newcommand{\showcomments}{yes}

\newsavebox{\commentbox}
{\ifthenelse{\equal{\showcomments}{yes}}%
{\footnotemark
        \begin{lrbox}{\commentbox}
        \begin{minipage}[t]{1.25in}\raggedright\sffamily\tiny
        \footnotemark[\arabic{footnote}]}
{\begin{lrbox}{\commentbox}}}%
{\ifthenelse{\equal{\showcomments}{yes}}%
{\end{minipage}\end{lrbox}\marginpar{\usebox{\commentbox}}}
{\end{lrbox}}}

\newcommand{\stabilizer}{\mathrm{Stab}}

\newcommand{\gate}{\mathfrak g}

\newcommand{\nest}{\sqsubseteq}
\newcommand{\propnest}{\sqsubsetneq}
\newcommand{\orth}{\bot}
\newcommand{\transverse}{\pitchfork}

\newcommand{\duaug}[2]{{#1}^{+{#2}}}
\newcommand{\link}{\mathrm{Lk}}
\newcommand{\Sat}{\mathrm{Sat}}
\newcommand{\dist}{\mathrm{d}}
\newcommand{\diam}{\mathrm{diam}}
\newcommand{\neb}{\mathcal N}

\newcommand{\cuco}[1]{{\mathcal #1}}

\newcommand{\tsh}[1]{\left\{\kern-.7ex\left\{#1\right\}\kern-.7ex\right\}}
\newcommand{\Tsh}[2]{\tsh{#2}_{#1}}
\newcommand{\ignore}[2]{\Tsh{#2}{#1}}

\newcommand{\frakS}{\mathfrak S}
\newcommand{\ww}{\wedge_{\text{min}}}

\makeatletter
\@tfor\next:=abcdefghijklmnopqrstuvwxyzABCDEFGHIJKLMNOPQRSTUVWXYZ\do{%
  \def\command@factory#1{%
    \expandafter\def\csname cal#1\endcsname{\mathcal{#1}}
    \expandafter\def\csname frak#1\endcsname{\mathfrak{#1}}
    \expandafter\def\csname scr#1\endcsname{\mathscr{#1}}
    \expandafter\def\csname bb#1\endcsname{\mathbb{#1}}
    \expandafter\def\csname rm#1\endcsname{\mathrm{#1}}
      \expandafter\def\csname bf#1\endcsname{\mathbf{#1}}
  }
 \expandafter\command@factory\next
}
\makeatother

\newcommand{\fontact}{{\mathcal C}}

\title[Many HHS are Combinatorial HHS]{A Combinatorial Structure for Many Hierarchically Hyperbolic Spaces}
\author[M. Hagen]{Mark Hagen}
    \address{(Mark Hagen) School of Mathematics, University of Bristol, Bristol, UK}
    \email{markfhagen@posteo.net}

\author[G. Mangioni]{Giorgio Mangioni}
    \address{(Giorgio Mangioni) Maxwell Institute and Department of Mathematics, Heriot-Watt University, Edinburgh, UK}
    \email{gm2070@hw.ac.uk}

\author[A. Sisto]{Alessandro Sisto}
    \address{(Alessandro Sisto) Maxwell Institute and Department of Mathematics, Heriot-Watt University, Edinburgh, UK}
    \email{a.sisto@hw.ac.uk}
    
\begin{document}

\setlength\parindent{0pt}

\begin{abstract}
The combinatorial hierarchical hyperbolicity criterion is a very useful way of constructing new hierarchically hyperbolic spaces (HHSs). We show that, conversely, HHSs satisfying natural assumptions (satisfied, for example, by mapping class groups) admit a combinatorial HHS structure. This can be useful in constructions of new HHSs, and also our construction clarifies how to apply the combinatorial HHS criterion to suspected examples. We also uncover connections between HHS notions and lattice theory notions.
\end{abstract}

\maketitle

\epigraph{This has nothing to do with links.}{A. S.}

\setcounter{tocdepth}{1}
\tableofcontents

\section*{Introduction} Showing that a given space or group is hierarchically hyperbolic yields a lot of information about it (\cite{HHS:quasiflats,HHP:coarse,ANS:UEG,HHL:proximal,DMS:stable} is a highly non-exhaustive list), but it is quite challenging to check the definition directly. To remedy this, a criterion was devised in \cite{BHMS} to show that a space is hierarchically hyperbolic, and roughly speaking the criterion involves checking that a certain simplicial complex has links which are hyperbolic, along with some combinatorial conditions. Such a complex is called a \emph{combinatorial HHS} structure (see Section~\ref{sec:what_is_CHHS} for background on combinatorial HHSs). This criterion has proven useful to show that various spaces and groups are HHS, see \cite{BHMS,ELTAG_HHS,HRSS,Russell:multicurve,DDLS}. It is natural to wonder whether there is a converse to this, namely whether every HHS admits a combinatorial HHS structure. We show that is true under mild assumptions on the HHS, see Theorem \ref{thmintro:main} below.

Besides the intrinsic interest of a converse statement, our theorem has proven useful in constructing new examples of HHSs, by introducing a combinatorial structure which is easier to manipulate. For instance, the second author used the flexibility of combinatorial structures coming from the main theorem to construct uncountably many coarse median structures on certain HHGs including the mapping class group of the 5-holed sphere \cite{short_I}. The same flexibility has in turn been used in \cite{short_II} to construct Dehn-filling-like quotients of certain HHGs, leading to a proof of the Hopf property for generic Artin groups. Additionally, our construction provides a blueprint for identifying candidate combinatorial HHS structures, and for instance this is one of the ingredients in \cite{Ragosta:artin}.

We defer the precise statements of the additional conditions to Section~\ref{sec:converse}; the reader can find a list of the conditions and where to find them below, together with a short discussion. For now, we just mention that the properties hold in many natural examples, including mapping class groups, see Section~\ref{sec:mcg}, and we now state (a stripped down version of) our main theorem:

\begin{thmintro}[see Theorem~\ref{thm:main}]\label{thmintro:main}
Let $(\cuco Z,\frakS)$ be a hierarchically hyperbolic space with weak wedges, clean containers, the orthogonals for non-split domains property, and dense product regions. Then there exists a combinatorial HHS $(X,\mathcal{W})$ such that $\cuco Z$ is quasi-isometric to $\mathcal{W}$.
\end{thmintro}

Theorem~\ref{thmintro:main} only contains the non-equivariant part of Theorem~\ref{thm:main}, but our constructions are equivariant in a suitable sense and compatible with the notion of hierarchically hyperbolic group, rather than just space, as would be needed for the application to quotients mentioned above; see Theorem~\ref{thm:comb_hhg} for the exact statement.

Our construction also clarifies how to construct, starting with a space or group that one suspects to be hierarchically hyperbolic, a candidate combinatorial HHS structure for it, which can then be used to show that the given space or group is indeed HHS. We explain this in Section~\ref{sec:construction}. We do not know whether the additional conditions we have to impose are necessary, but we provide an example where our construction fails to yield a combinatorial HHS structure in the absence of the additional conditions, in Section~\ref{sec:CHHS_for_orth_sets}. In fact, we uncovered intriguing connections with lattice theory that arise from these considerations, see in particular Remark~\ref{rem:lattice}. Roughly, modifying an HHS structure to ensure the additional condition reduces to a problem in lattice theory, see Question~\ref{question:lattice}. To highlight the connections with lattice theory, we note that the notion of an ortholattice is very closely related to those of wedges and clean containers (which have appeared in the HHS literature repeatedly, see e.g. \cite{BerlaiRobbio,ABD,Rcubes, Russell:relative, AbbottBehrstock, Hagen:non-colourable}), see Definition~\ref{defn:orth_set}.

To conclude this subsection, we suggest that a possible application of the aforementioned construction could be showing that mapping class groups of finite type non-orientable surfaces are hierarchically hyperbolic.

\subsection*{A true converse}
A second aim of this paper is to clarify how various conditions on an HHS structure relate to each other and to properties of combinatorial HHS structure. As this is fairly technical we do not make precise statements in the introduction, but we refer the reader to Section~\ref{sec:additional_hyp} and in particular Lemma~\ref{lem:comparison_between_prop}; we believe this can also be useful for applications. An especially striking output of this study is the following “true converse” theorem, which provides an actual equivalence between combinatorial HHS structures and HHS structures, each satisfying natural conditions (see list below for where to find each condition).

\begin{thmintro}\label{thmintro:strong-equivalence}
    Let $({\cuco Z},\frakS)$ be a hierarchically hyperbolic space. Then ${\cuco Z}$ has wedges, clean containers and the strong orthogonal property if and only if there exists a CHHS $(X,\mathcal{W})$ with simplicial wedges and simplicial containers such that $\mathcal{W}$ is quasi-isometric to ${\cuco Z}$.
\end{thmintro}

We note that the standard HHS structure on mapping class groups does not satisfy the stronger conditions above, but we construct a different structure which does in Section~\ref{sec:mcg}. We now state this result, together with the aforementioned result on the standard HHS structure, which in fact confirms the speculations from \cite[Subsection 1.6]{BHMS}.

The following combines Theorem~\ref{thm:mcg_cHHS} and Theorem~\ref{thm:strong_mcg} (which in fact give combinatorial HHG structures):
\begin{thmintro}
    Let $S$ be obtained from a closed connected oriented surface of finite genus by deleting finitely many points and open discs.
    
    There exists a combinatorial HHS structure $(X,\mathcal{W})$ for $\mathcal{MCG}(S)$, where $X$ is the blow-up of the curve graph of $S$, obtained by replacing every curve with the cone over its annular curve graph.

     Moreover, there exists a (different) combinatorial HHS structure for $\mathcal{MCG}(S)$ with simplicial wedges and simplicial containers.
\end{thmintro}
In the new HHS structures, the elements of the index set are still isotopy classes of certain essential subsurfaces, as in \cite{HHS_II}, and in particular the unbounded hyperbolic spaces involved in the HHS structure all still come from curve graphs of connected subsurfaces.  The difference between the HHS structures in Section~\ref{sec:mcg}, which are used to construct combinatorial HHS structures, and the “standard” one is in which subsurfaces with bounded curve graph are present; care in this choice is what allows us to arrange extra combinatorial properties of the index set needed to make a CHHS structure.

As a final note on mapping class groups, we summarise in the following remark the connection between our construction and clean markings:
\begin{remintro}[Relation with clean markings]\label{rem:clean_markings_for_mcg}
    The reader familiar with Masur and Minsky's graph of complete clean markings from \cite{MasurMinsky2} will notice that the combinatorial HHS structure for $\mathcal{MCG}(S)$ that our main Theorem provides has the same flavour of the marking graph. Indeed, a maximal simplex of the graph $X$ from Definition~\ref{defn:blow-up} will correspond to a choice of a maximal collection of disjoint annuli $A_1,\ldots, A_k$ (that is, a pants decomposition), plus a choice of a point $x_i$ inside the annular curve graph associated to $A_i$ for every $i=1,\ldots, k$ (that is, a transversal for every curve in the pants decomposition). Hence, a maximal simplex corresponds to a complete marking. Moreover, some of the $\mathcal{W}$-edges we define in Definition~\ref{defn:W-edge} correspond to elementary moves. Indeed, let $\Sigma,\Delta\subset X$ be two maximal simplices, and suppose that their supports differ by a single curve (say, the support of $\Sigma$ is $\alpha\cup P$ and the support of $\Delta$ is $\beta\cup P$, for some almost-maximal collection of pairwise disjoint curves $P$). Then these simplices are $\mathcal{W}$-adjacent if and only if:
    \begin{itemize}
        \item $\alpha$ and $\beta$ are close in the curve graph of the subsurface of complexity $1$ cut out by $P$;
        \item $\alpha$ projects close to the coordinate prescribed by $\Delta$ in the annular curve graph of $\beta$;
        \item the same holds with $\alpha$ and $\beta$ swapped.
    \end{itemize}
    In other words, our $\mathcal{W}$--edges detect when one replaces a curve with one of its transversals, and this corresponds to one of the elementary moves from \cite{MasurMinsky2}.
\end{remintro}

\subsection*{Additional conditions}
The additional conditions on an HHS structure in Theorems~\ref{thmintro:main},\ref{thmintro:strong-equivalence} are stated precisely in the following places:

\begin{itemize}
    \item Wedges are defined in Property~\ref{defn:wedge_property}.
    \item Clean Containers are defined in Property~\ref{property:clean_containers}.
    \item Dense product regions are defined in Property~\ref{property:dpr}.
    \item Orthogonals for non-split domains are defined in Property~\eqref{property:orthogonal_for_non_split}.
    \item Strong orthogonal property is defined in Property~\ref{property:strong_orth} (and see also Remark~\ref{rem:lattice}).
\end{itemize}

Wedges and clean containers are standard assumptions on an HHS structure introduced in \cite{BerlaiRobbio} and \cite{ABD} respectively --- wedges make the nesting poset into a lattice, while clean containers make it into a complemented poset.  Dense product regions is a coboundedness assumption automatically satisfied by HHGs.  Strong orthogonality corresponds to \emph{orthomodularity} of the nesting lattice, and is the source of the lattice-theoretic question mentioned earlier.  It is one way of verifying orthogonals for non-split domains, which is really the main enabling assumption in Theorem~\ref{thm:main} and is modelled on the role of boundary annuli in the HHS structure on mapping class groups.

The additional conditions on a combinatorial HHS are:

\begin{itemize}
        \item Simplicial wedges are defined in Definition~\ref{defn:simp_wedge}.
    \item Simplicial containers are defined in Definition~\ref{defn:simp_cont}.
\end{itemize}

These conditions are very natural properties one might hope for from a simplicial complex, along the lines that containment of links of simplices corresponds to reverse containment of the simplices.

\subsection*{Outline of the paper}
Sections~\ref{sec:background} and~\ref{sec:what_is_CHHS} contain all relevant definitions and facts about (combinatorial) HHS. Section~\ref{sec:converse} gathers the hypotheses of the main result of this paper, which is Theorem~\ref{thm:main} and states that a HHS $(\cuco Z, \frakS)$ with some additional assumption is quasi-isometric to a combinatorial HHS. The actual construction of the candidate combinatorial HHS $(X,\mathcal{W})$ is in Section~\ref{sec:construction} (see in particular Definitions~\ref{defn:blow-up} and~\ref{defn:W-edge}). The quasi-isometry $f:\,\mathcal{W} \to \cuco Z$ is constructed in Definition~\ref{defn:realisation_map}, and maps any maximal simplex of $X$ to one of its realisation points (in the sense of the partial realisation axiom~\eqref{item:dfs_partial_realisation}).

In Section~\ref{sec:verify_axioms} we verify that, under the assumptions from Section~\ref{sec:converse}, the pair $(X,\mathcal{W})$ is a combinatorial HHS and $f$ is a quasi-isometry (see Assumption~\ref{ass:standing_assumption} and Theorem~\ref{thm:XWf}). In Section~\ref{sec:converse_group} we show that our construction is equivariant, meaning that whenever a group $G$ acts on $(\cuco Z,\frakS)$ by hieromorphisms then it also acts on $(X,\mathcal{W})$ (see Theorem~\ref{thm:G-equivariance}). Then we use this fact to prove that, whenever a group acts metrically properly and coboundedly on $(\cuco Z,\frakS)$ and some other mild assumptions hold, then $G$ has a structure of hierarchically hyperbolic group coming from the action on a combinatorial HHS (see Theorem~\ref{thm:comb_hhg}).

In Section~\ref{sec:additional_hyp} we present some more “natural” hypotheses that one could require on $(\cuco Z, \frakS)$, and we show how they relate to each other and to the ones from Section~\ref{sec:converse}. In Section~\ref{sec:HHS_iff_CHHS} we establish an equivalence between strong orthogonality properties on the HHS structure of $(\cuco Z, \frakS)$ and some strong intersection properties on the links of the associated combinatorial HHS $(X,\mathcal{W})$ (see Theorem~\ref{thm:strong_iff}). 

In Section~\ref{sec:mcg} we apply our results to the mapping class group of a compact orientable surface, with the usual HHS structure (the one from, e.g., \cite[Section 11]{HHS_II}), showing that it admits a combinatorial HHS structure whose underlying graph is a certain blow-up of the curve graph. This confirms the speculations from \cite[Subsection 1.6]{BHMS}. Moreover, we show that Theorem~\ref{thm:strong_iff} applies if one adds to the index set some non-essential subsurfaces, including pairs of pants (see Theorem~\ref{thm:strong_mcg}).

Finally, in Section~\ref{sec:CHHS_for_orth_sets} we illustrate the necessity of the hypotheses of Theorem~\ref{thm:main}, by providing a counterexample of an unbounded space $\cuco Z$ for which the construction from Section~\ref{sec:construction} can only yield a bounded CHHS. Remarkably, $\cuco Z$ can be chosen to be a CAT(0) cube complex with a factor system, with the usual HHS structure (i.e. the one from \cite[Remark 13.2]{HHS_I}). Then we speculate on which conditions on the factor system could allow one to modify the HHS structure in order to satisfy our hypotheses.

\subsection*{Acknowledgements}  We thank Carolyn Abbott and Alexandre Martin for helpful discussions.  MH thanks Montserrat Casals-Ruiz for being a strong proponent of the “orthogonal set” viewpoint during work on \cite{Rcubes}, which influenced the ideas in Section~\ref{sec:CHHS_for_orth_sets} here.  We thank Jason Behrstock for some useful comments on an earlier version, and the referee for numerous very helpful comments.

\section{Background on hierarchical hyperbolicity}\label{sec:background}
\subsection{Axioms}\label{subsec:axioms}
We recall from~\cite{HHS_II} the definition of a hierarchically hyperbolic space.
\begin{defn}[HHS]\label{defn:HHS}
The $q$--quasigeodesic space  $(\cuco Z,\dist_{\cuco Z})$ is a \emph{hierarchically hyperbolic space} if there exists $E\geq0$, called the \emph{HHS constant}, an index set $\frakS$, whose elements will be referred to as \emph{domains}, and a set $\{\fontact  U\mid U\in\frakS\}$ of $E$--hyperbolic spaces $(\fontact  U,\dist_U)$, called \emph{coordinate spaces},  such that the following conditions are satisfied:
\begin{enumerate}
\item\textbf{(Projections.)}\label{item:dfs_curve_complexes}
There is a set $\{\pi_U: \cuco Z\rightarrow 2^{\fontact  U}\mid U\in\frakS\}$ of \emph{projections} sending points in $\cuco Z$ to sets of diameter bounded by $E$ in the various $\fontact  U\in\frakS$. Moreover, for all $U\in\frakS$, the coarse map $\pi_U$ is $(E,E)$--coarsely Lipschitz and $\pi_U(\cuco Z)$ is $E$--quasiconvex in $\fontact  U$.

\item \textbf{(Nesting.)} \label{item:dfs_nesting}
$\frakS$ is equipped with a partial order $\nest$, and either $\frakS=\emptyset$ or $\frakS$ contains a unique $\nest$--maximal element, denoted by $S$. When $V\nest U$, we say $V$ is \emph{nested} in $U$. For each $U\in\frakS$, we denote by $\frakS_U$ the set of $V\in\frakS$ such that $V\nest U$. Moreover, for all $U,V\in\frakS$ with $V\propnest U$ there is a specified subset $\rho^V_U\subset\fontact  U$ with $\diam_{\fontact  U}(\rho^V_U)\leq E$. There is also a \emph{projection} $\rho^U_V: \fontact U\rightarrow 2^{\fontact V}$. (The similarity in notation is justified by viewing $\rho^V_U$ as a coarsely constant map $\fontact V\rightarrow 2^{\fontact  U}$.)
 
\item \textbf{(Orthogonality.)} \label{item:dfs_orthogonal}
$\frakS$ has a symmetric and anti-reflexive relation called \emph{orthogonality}: we write $U\orth V$ when $U,V$ are orthogonal. Also, whenever $V\nest U$ and $U\orth W$, we require that $V\orth W$. We require that for each $T\in\frakS$ and each $U\in\frakS_T$ such that $\{V\in\frakS_T\mid V\orth U\}\neq\emptyset$, there exists $W\in\frakS_T-\{T\}$, which we call a \emph{container} for $U$ inside $T$, so that whenever $V\orth U$ and $V\nest T$, we have $V\nest W$. Finally, if $U \orth V$, then $U,V$ are not $\nest$--comparable.

\item \textbf{(Transversality and consistency.)}\label{item:dfs_transversal}
If $U,V\in\frakS$ are not orthogonal and neither is nested in the other, then we say $U,V$ are \emph{transverse}, denoted $U\transverse V$. In this case there are sets $\rho^V_U\subseteq\fontact U$ and $\rho^U_V\subseteq\fontact  V$, each of diameter at most $E$ and satisfying:
$$\min\left\{\dist_{U}(\pi_U(z),\rho^V_U),\dist_{V}(\pi_V(z),\rho^U_V)\right\}\leq E$$
for all $z\in \cuco Z$. Furthermore, for $U,V\in\frakS$ satisfying $V\nest U$ and for all $z\in\cuco Z$, we have: 
$$\min\left\{\dist_{U}(\pi_U(z),\rho^V_U),\diam_{\fontact V}(\pi_V(z)\cup\rho^U_V(\pi_U(z)))\right\}\leq E.$$ 
The preceding two inequalities are the \emph{consistency inequalities} for points in $\cuco Z$.
 
Finally, if $U\nest V$, then $\dist_W(\rho^U_W,\rho^V_W)\leq E$ whenever $W\in\frakS$ satisfies either $V\propnest W$ or $V\transverse W$ and $W\not\bot U$.
 
\item \textbf{(Finite complexity.)} \label{item:dfs_complexity}
There exists $n\geq0$, the \emph{complexity} of $\cuco Z$ (with respect to $\frakS$), so that any set of pairwise--$\nest$--comparable elements has cardinality at most $n$.
  
\item \textbf{(Large links.)} \label{item:dfs_large_link_lemma}
Let $U\in\frakS$, let $z,z'\in\cuco Z$ and let $N=\dist_{_U}(\pi_U(z),\pi_U(z'))$. Then there exists $\{T_i\}_{i=1,\dots,\lfloor N\rfloor}\subseteq\frakS_U- \{U\}$ such that, for any domain $T\in\mathfrak S_U-\{U\}$, either $T\in\frakS_{T_i}$ for some $i$, or $\dist_{T}(\pi_T(z),\pi_T(z'))<E$.  Also, $\dist_{U}(\pi_U(z),\rho^{T_i}_U)\leq N$ for each $i$.

\item \textbf{(Bounded geodesic image.)}\label{item:dfs:bounded_geodesic_image}
For all $U\in\frakS$, all $V\in\frakS_U- \{U\}$, and all geodesics $\gamma$ of $\fontact  U$, either $\diam_{\fontact  V}(\rho^U_V(\gamma))\leq E$ or $\gamma\cap\neb_E(\rho^V_U)\neq\emptyset$.
 
\item \textbf{(Partial realisation.)} \label{item:dfs_partial_realisation}
Let $\{V_j\}$ be a family of pairwise orthogonal elements of $\frakS$, and let $p_j\in \pi_{V_j}(\cuco Z)\subseteq \fontact  V_j$. Then there exists $z\in \cuco Z$, which we call a \emph{partial realisation point} for the family, so that:
\begin{itemize}
\item $\dist_{V_j}(z,p_j)\leq E$ for all $j$,
\item for each $j$ and 
each $V\in\frakS$ with $V_j\nest V$, we have 
$\dist_{V}(z,\rho^{V_j}_V)\leq E$, and
\item for each $j$ and 
each $V\in\frakS$ with $V_j\transverse V$, we have $\dist_V(z,\rho^{V_j}_V)\leq E$.
\end{itemize}

\item\textbf{(Uniqueness.)} For each $\kappa\geq 0$, there exists
$\theta_u=\theta_u(\kappa)$ such that if $x,y\in\cuco Z$ and
$\dist_{\cuco Z}(x,y)\geq\theta_u$, then there exists $V\in\frakS$ such
that $\dist_V(x,y)\geq \kappa$.\label{item:dfs_uniqueness}
\end{enumerate}
 
We often refer to $\frakS$, together with the nesting and orthogonality relations, and the projections as a \emph{hierarchically hyperbolic structure} for the space $\cuco Z$.  Observe that $\cuco Z$ is hierarchically hyperbolic with respect to $\frakS=\emptyset$, i.e., hierarchically hyperbolic of complexity $0$, if and only if $\cuco Z$ is bounded. Similarly, $\cuco Z$ is hierarchically hyperbolic of complexity $1$ with respect to the index set $\frakS=\{\cuco Z\}$, if and only if $\cuco Z$ is hyperbolic.
\end{defn}

\begin{notation}\label{notation:suppress_pi}
Where it will not cause confusion, given $U\in\frakS$, we will often suppress the projection map $\pi_U$ when writing distances in $\fontact  U$, i.e., given $x,y\in\cuco Z$ and $p\in\fontact  U$  we write $\dist_U(x,y)$ for $\dist_U(\pi_U(x),\pi_U(y))$ and $\dist_U(x,p)$ for $\dist_U(\pi_U(x),p)$. Note that when we measure distance between a pair of sets (typically both of bounded diameter) we are taking the minimum distance between the two sets. Given $A\subset \cuco Z$ and $U\in\frakS$ we set $$\pi_{U}(A)=\bigcup_{a\in A}\pi_{U}(a).$$
\end{notation}

\subsection{Useful facts about HHS}\label{subsec:facts_about_HHS}
We now recall results from \cite{HHS_II} that will be useful later on. 

\begin{lemma}[{\cite[Lemma 1.5]{DHS}}]\label{lem:close_proj_of_orthogonals}
Let $U, V, W\in\frakS$ satisfying $U\orth V$, and $U,V\not\bot W$, and $W\not \nest U,V$. Then $\dist_V(\rho^U_W,\rho^V_W)\le 2E$.
\end{lemma}

\begin{rem}[Normalisation]\label{rem:normalisation}
    As argued in \cite[Remark 1.3]{HHS_II}, it is always possible to assume that the HHS structure is \emph{normalised}, that is, for every $U\in \frakS$ the projection $\pi_U:\,\cuco Z\to \fontact U$ is uniformly coarsely surjective. In order to do so, one roughly replaces every $\fontact U$ with $\pi_U(\cuco Z)$, which is itself hyperbolic since it is quasiconvex in $\fontact U$, and then replaces every projection $\rho^V_U$ with the composition $p_{U}\circ\rho^V_U $, where $p_U:\,\fontact U \to \pi_U(\cuco Z)$ is the coarse closest point retraction. The resulting space, which is again hierarchically hyperbolic, has the same set of domains $\frakS$ with the same relations of nesting and orthogonality.
\end{rem}

\begin{ass}\label{assump:normalised}
    In view of the remark above, we will always assume that the HHS structures we consider are normalised.
\end{ass}

\begin{defn}[Consistent tuple]\label{defn:consistent_tuple}
Let $\kappa\geq1$ and let $(b_U)_{U\in\frakS}\in\prod_{U\in\frakS}2^{\fontact  U}$ be a tuple such that for each $U\in\frakS$, the $U$--coordinate  $b_U$ has diameter $\leq\kappa$.  Then $(b_U)_{U\in\frakS}$ is \emph{$\kappa$--consistent} if for all $V,W\in\frakS$, we have $$\min\{\dist_V(b_V,\rho^W_V),\dist_W(b_W,\rho^V_W)\}\leq\kappa$$
whenever $V\transverse W$ and 
$$\min\{\dist_W(b_W,\rho^V_W),\diam_V(b_V\cup\rho^W_V(b_W))\}\leq\kappa$$
whenever $V\propnest W$.
\end{defn}

\begin{thm}[Realisation {\cite[Theorem~3.1]{HHS_II}}]\label{thm:realisation}
Let $(\cuco Z,\frakS)$ be a hierarchically hyperbolic space. Then for each $\kappa\geq1$, there exists $\theta=\theta(\kappa)$ so that,  for any $\kappa$--consistent tuple $(b_U)_{U\in\frakS}$, there exists $x\in\cuco Z$ such that $\dist_V(x,b_V)\leq\theta$ for all $V\in\frakS$.
\end{thm}

\noindent Observe that the uniqueness axiom (Definition~\eqref{item:dfs_uniqueness}) implies that the \emph{realisation point} $x$ for $(b_U)_{U\in\frakS}$ provided by Theorem~\ref{thm:realisation} is coarsely unique.  

\begin{defn}[Product regions and factors]\label{defn:factor}
Fix a constant $\kappa\ge1$. For any domain $U$, let $F_U$ be the set of $\kappa$-consistent tuples \emph{for $U$}, that is, all tuples $(b_V)_{V\in\frakS_U}$ that satisfy the consistency inequalities involving only domains nested in $U$. Similarly, one can define $E_U$ as the set of $\kappa$-consistent tuples of the form $(b_V)_{V\orth U}$.

Now let $P_U=F_U\times E_U$, which we call the \emph{product region} associated to $U$. By the realisation Theorem~\ref{thm:realisation} there is a coarsely well-defined map $\phi:P_U \to \cuco Z$. If we fix $e\in E_U$, the image of the \emph{factor} $F_U\times\{e\}$, which we will still denote by $F_U$ when the dependence on $e$ is irrelevant, can be endowed with the subspace metric, which makes it a sub-HHS of $\cuco Z$ with domain set $\frakS_U=\{V\in\frakS \mid V\nest U\}$. Two parallel copies $F_U\times\{e\}$ and $F_U\times\{e'\}$ are quasi-isometric (see e.g. \cite[Section 2.2]{DHScorrection} for more details), thus the metric structure on $F_U$ is well-defined up to quasi-isometry.

A similar argument holds for $E_U$. For more details on product regions, see \cite[Section 15]{Rcubes}.
\end{defn}

It will often be convenient to think of $F_U$ as an abstract space, instead of as a subspace of $\cuco Z$. This way, whenever $V\nest U$, we have a (non-unique) embedding $F_V\to F_U$, given as follows. Let $E^U_{V}$ be the set of $\kappa$-consistent tuples of the form 
$$(b_W)_{W\propnest U,\,W\orth V_i\,\forall i=1,\ldots, k},$$
and choose $e\in E^U_V$. Then define a map $F_V\to F_U$ by sending a tuple $(y_W)_{W\nest V}$ to the tuple $(x_W)_{W\nest U}$, defined as follows:
$$x_W=\begin{cases}
    y_W \mbox{ if }W\nest V;\\
    \rho^V_W \mbox{ if }V\propnest W \mbox{ or }V\transverse W;\\
    e_W \mbox{ if }W\orth V.
\end{cases}$$
In other words, we extend the tuple $y$ “naturally” whenever we have a well-defined projection from $V$ to $W$, and then we choose consistent coordinates whenever $W\orth V$. This kind of arguments will be recurrent throughout the paper.

\begin{defn}[Relative product regions]\label{defn:rel_prod_region}
Fix a constant $\kappa\ge 0$, and let $U,V_1,\ldots, V_k\in \frakS$ be such that $V_i\nest U$ and $V_i\orth V_j$ for every $i,j\le k$. The \emph{relative} product region associated to $V_1,\ldots, V_k$ inside $U$ is defined as
$$P^U_{\{V_i\}}=F_{V_1}\times\ldots\times F_{V_k}\times E^U_{\{V_i\}}\subset F_U,$$
where $E^U_{\{V_i\}}$ is the set of $\kappa$-consistent tuples of the form 
$$(b_W)_{W\propnest U,\,W\orth V_i\,\forall i=1,\ldots, k}.$$
With a slight abuse of notation, whenever the ambient domain $U$ is clear we will drop the superscript and refer to the relative product region simply as $P_{\{V_i\}}$.
\end{defn}
\begin{rem}
    For the rest of the paper, unless otherwise stated, all factors and relative product regions are with respect to $\kappa=20E$, where $E$ is the HHS constant. This specific choice of $\kappa$ will be motivated in Assumption~\ref{ass:standing_assumption}.
\end{rem}

\begin{thm}[{Distance formula \cite[Theorem~4.5]{HHS_II}}]\label{thm:distance_formula}
Let $(\cuco Z,\frakS)$ be a hierarchically hyperbolic space.  Then
there exists $s_0$ such that for all $s\geq s_0$, there exist $C,K$ so
that for all $x,y\in\cuco Z$,
$$\dist(x,y)\asymp_{_{K,C}}\sum_{U\in\mathfrak
S}\ignore{\dist_U(x,y)}{s}.$$
\end{thm}

\noindent (The notation $\ignore{A}{B}$ denotes the quantity which is $A$ if $A\geq B$ and $0$ otherwise. The notation $A\asymp_{K,C} B$ means $A\leq K B+C$ and $B\leq K A+C$.)

\subsection{Groups acting on HHS}\label{subsec:HHG}

First, we need to discuss which group actions we allow on a hierarchically hyperbolic space. The following are some definitions from \cite{HHS_II} and \cite{BHMS}:

\begin{defn}[Automorphism]\label{defn:auto}
Let $({\cuco Z},\frakS)$ be a HHS. An \emph{automorphism} consists of a map $g:\,{\cuco Z}\to {\cuco Z}$, a bijection $g^\sharp:\, \frakS\to \frakS$ preserving nesting and orthogonality, and, for each $U\in\frakS$, an isometry $g^\diamond(U):\,\fontact U\to \fontact(g^\sharp(U))$ for which the following two diagrams commute for all $U,V\in\frakS$ such that $U\propnest V$ or $U\transverse V$:
$$\begin{tikzcd}
{\cuco Z}\ar{r}{g}\ar{d}{\pi_U}&{\cuco Z}\ar{d}{\pi_{g^\sharp (U)}}\\
\fontact U\ar{r}{g^\diamond (U)}&\fontact (g^\sharp (U))\\
\end{tikzcd}$$
and
$$\begin{tikzcd}
\fontact U\ar{r}{g^\diamond (U)}\ar{d}{\rho^U_V}&\fontact (g^\sharp (U))\ar{d}{\rho^{g^\sharp (U)}_{g^\sharp (V)}}\\
\fontact V\ar{r}{g^\diamond (V)}&\fontact (g^\sharp (V))\\
\end{tikzcd}$$
Notice that $g$ must be a uniform quasi-isometry by the distance formula, Theorem~\ref{thm:distance_formula}. Whenever it will not cause ambiguity, we will abuse notation by dropping the superscripts and just calling all maps $g$.
\end{defn}

We say that two automorphisms $g,g'$ are \emph{equivalent}, and we write $g\sim g'$, if $g^\sharp=(g')^\sharp$ and $g^\diamond(U)=(g')^\diamond(U)$ for each $U\in\frakS$.  Given an automorphism $g$, a quasi-inverse $\bar g$ for $g$ is an automorphism with $\bar g^\sharp=(g^\sharp)^{-1}$ and such that, for every $U\in \frakS$, $\bar g^\diamond(U)=g^\diamond(U)^{-1}$. Since the composition of two automorphisms is an automorphism, the set of equivalence classes of automorphisms forms a group, denoted $\text{Aut}(\frakS)$.

\begin{defn}\label{defn:action_on_HHS}
    A finitely generated group $G$ \emph{acts} on a HHS $({\cuco Z},\frakS)$ by automorphisms if there is a homomorphism $G\to \text{Aut}(\frakS)$.
\end{defn}

\begin{rem}\label{rem:commutativity_of_diag_in_hhg}
The original definition of an automorphism, which is \cite[Definition 1.20]{HHS_II}, only requires the diagrams from Definition~\ref{defn:auto} to coarsely commute (with uniform constants). However,
as shown in \cite[Section 2.1]{DHScorrection}, if $G$ acts on $({\cuco Z},\frakS)$ in the sense of \cite{HHS_II} then one can ensure that the diagrams genuinely commute by perturbing every $\pi_U:\,{\cuco Z}\to \fontact U$ and every $\rho^U_V$, whenever the quantity is defined, by a uniformly bounded amount. This way, up to a single initial change in the constant $E$, the HHS structure is unaffected, meaning that the new structure has the same domain set $\frakS$ with the same relations and the same coordinate spaces. 
\end{rem}

\begin{defn}[HHG]\label{defn:HHG}
A finitely generated group $G$ is \emph{hierarchically hyperbolic} if there exists a hierarchically hyperbolic space $({\cuco Z},\frakS)$ and an action $G\rightarrow\text{Aut}(\frakS)$ so that the uniform quasi-action of $G$ on ${\cuco Z}$ is metrically proper and cobounded and $\frakS$ contains finitely many $G$--orbits. Then we can equip $G$ with a HHS structure, whose domains and coordinate spaces are the same as the ones for ${\cuco Z}$ and whose projections are obtained by precomposing the projections for $({\cuco Z},\frakS)$ with the $G$-equivariant quasi-isometry $G\to {\cuco Z}$ given by the Milnor-\v{S}varc lemma. 
\end{defn}

\section{Combinatorial HHSs}\label{sec:what_is_CHHS} 
In this section we recall the definition of a combinatorial HHS and its hierarchically hyperbolic structure, as first introduced in \cite{BHMS}. 

\subsection{Basic definitions}
Let $X$ be a simplicial graph. 

\begin{defn}[Induced subgraph]
    Given a subset $S\subseteq X^{(0)}$ of the set of vertices of $X$, the subgraph \emph{spanned} by $S$ is the complete subgraph of $X$ with vertex set $S$.
\end{defn}

\begin{defn}[Join, link, star]\label{defn:join_link_star}
Given disjoint simplices $\Delta,\Delta'$ of $X$, we let $\Delta\star\Delta'$ denote the simplex spanned by $\Delta^{(0)}\cup\Delta'^{(0)}$, if it exists. 

For each simplex $\Delta$, the \emph{link} $\link(\Delta)$ is the union of 
all simplices $\Sigma$ of $X$ such that $\Sigma\cap\Delta=\emptyset$ and $\Sigma\star\Delta$ is a simplex of $X$.  Observe that $\link(\Delta)=\emptyset$ if and only if $\Delta$ is a maximal simplex. The link of a subgraph of $X$ is the intersection of the links of its vertices.

The \emph{star} of $\Delta$ is $\text{Star}(\Delta):=\link(\Delta)\star\Delta$, i.e. the union of all simplices of $X$ that contain $\Delta$.
\end{defn}

\begin{defn}[$X$--graph, $\mathcal{W}$--augmented graph]\label{defn:X_graph}
An \emph{$X$--graph} is a graph $\mathcal{W}$ whose vertex set is the set of all maximal simplices of $X$.

For a simplicial graph $X$ and an $X$--graph $\mathcal{W}$, the \emph{$\mathcal{W}$--augmented graph} $\duaug{X}{\mathcal{W}}$ is the graph defined as follows:
\begin{itemize}
     \item the $0$--skeleton of $\duaug{X}{\mathcal{W}}$ is $X^{(0)}$;
     \item if $v,w\in X^{(0)}$ are adjacent in $X$, then they are adjacent in $\duaug{X}{\mathcal W}$; 	
     \item if two vertices in $\mathcal W$ are adjacent, then we consider $\sigma,\rho$, the associated maximal simplices of $X$, and in $\duaug{X}{\mathcal W}$ we connect each vertex of $\sigma$ to each vertex of $\rho$.
\end{itemize}
We equip $\mathcal W$ with the usual path-metric, in which each edge has unit length, and do the same for $\duaug{X}{\mathcal W}$.
\end{defn}

\subsection{HHS structure}\label{subsec:W_structure}
\cite[Theorem 1.18]{BHMS} states that, under some assumptions on the pair $(X,\mathcal W)$, $\mathcal W$ has the hierarchically hyperbolic structure described below. First, we define what will be the index set. 

\begin{defn}[Equivalence between simplices, saturation]\label{defn:simplex_equivalence}
For $\Delta,\Delta'$ simplices of $X$, we write $\Delta\sim\Delta'$ to mean $\link(\Delta)=\link(\Delta')$. We denote the $\sim$--equivalence class of $\Delta$ by $[\Delta]$.
Let $\Sat(\Delta)$ denote the set of vertices $v\in X$ for which there exists a simplex $\Delta'$ of $X$ such that $v\in\Delta'$ and $\Delta'\sim\Delta$, i.e.
$$\Sat(\Delta)=\left(\bigcup_{\Delta'\in[\Delta]}\Delta'\right)^{(0)}.$$
We denote by $\frakS$ the set of $\sim$--classes of \textbf{non-maximal} simplices in $X$.
\end{defn}

Next we introduce the candidate coordinate spaces:
\begin{defn}[Complement, link subgraph]\label{defn:complement}
Let $\mathcal W$ be an $X$--graph.  For each simplex $\Delta$ of $X$, let $Y_\Delta$ be the subgraph of $\duaug{X}{\mathcal W}$ induced by the set $(\duaug{X}{\mathcal W})^{(0)}-\Sat(\Delta)$ of vertices.

Let $\fontact (\Delta)$ be the induced subgraph of $Y_\Delta$ spanned by $\link(\Delta)^{(0)}$.  Note that $\fontact (\Delta)=\fontact (\Delta')$ whenever $\Delta\sim\Delta'$. (We emphasise that we are taking links in $X$, not in $\duaug{X}{\mathcal W}$, and then considering the subgraphs of $Y_\Delta$ induced by those links.)
\end{defn}

The following is the equivalent of the finite complexity Axiom~\eqref{item:dfs_complexity} in the combinatorial framework:
\begin{defn}[Finite complexity]\label{defn:finite_complexity_cHHS}
The simplicial complex $X$ has \emph{finite complexity} if there exists $n\in\mathbb{N}$ so that any chain $\link(\Delta_1)\subsetneq\dots\subsetneq\link(\Delta_i)$, where each $\Delta_j$ is a simplex of $X$, has length at most $n$; the minimal such $n$ is the \emph{complexity} of $X$.
\end{defn}

The following is the main definition from \cite{BHMS}:

\begin{defn}[Combinatorial HHS]\label{defn:combinatorial_HHS}
A \emph{combinatorial HHS} $(X,\mathcal W)$ consists of a simplicial graph $X$ and an $X$--graph $\mathcal W$ satisfying the following conditions:
\begin{enumerate}
    \item \label{item:cHHS_flag} $X$ has complexity $n<+\infty$, as in Definition~\ref{defn:finite_complexity_cHHS};
    \item \label{item:cHHS_delta} There is a constant $\delta$ so that for each non-maximal simplex $\Delta$, the subgraph $\fontact (\Delta)$ is $\delta$--hyperbolic and $(\delta,\delta)$--quasi-isometrically embedded in $Y_\Delta$, where $Y_\Delta$ is as in Definition~\ref{defn:complement};
    \item \label{item:cHHS_join} Whenever $\Delta$ and $\Sigma$ are non-maximal simplices for which there exists a non-maximal simplex $\Gamma$ such that $\link(\Gamma)\subseteq\link(\Delta)\cap \link(\Sigma)$, and $\diam(\fontact  (\Gamma))\geq \delta$, then there exists a simplex $\Pi$ which extends $\Sigma$ such that $\link(\Pi)\subseteq \link(\Delta)$, and all $\Gamma$ as above satisfy $\link(\Gamma)\subseteq\link(\Pi)$;
    \item \label{item:C_0=C} If $v,w$ are distinct non-adjacent vertices of $\link(\Delta)$, for some simplex $\Delta$ of $X$, contained in $\mathcal W$-adjacent maximal simplices, then they are contained in $\mathcal W$-adjacent simplices of the form $\Delta\star\Sigma$.
\end{enumerate}
\end{defn}

In order to complete the HHS structure on $\mathcal{W}$ we are left to define nesting and orthogonality relations on $\frakS$, and projections between coordinate spaces.

\begin{defn}[Nesting, orthogonality, transversality, complexity]\label{defn:nest_orth}
Let $X$ be a simplicial graph.  Let $\Delta,\Delta'$ be non-maximal simplices of $X$.  Then:
\begin{itemize}
     \item $[\Delta]\nest[\Delta']$ if $\link(\Delta)\subseteq\link(\Delta')$;
     \item $[\Delta]\orth[\Delta']$ if $\link(\Delta')\subseteq \link(\link(\Delta))$.
\end{itemize}
If $[\Delta]$ and $[\Delta']$ are neither $\orth$--related nor $\nest$--related, we write 
$[\Delta]\transverse[\Delta']$.

Note that $[\emptyset]$ is the unique $\nest$--maximal $\sim$--class of simplices in $X$ and that $\nest$ is a partial ordering on the set of $\sim$--classes of simplices in $X$.  Notice that the simplicial graph $X$ has finite complexity, in the sense of Definition~\ref{defn:finite_complexity_cHHS}, if there exists $n\in\mathbb{N}$ so that any $\nest$--chain has length at most $n$; the minimal such $n$ is the complexity of $X$.  
\end{defn}

\begin{rem}
 The definition of $\orth$ says that any vertex in the link of $\Delta'$ is connected to any vertex in the link of $\Delta$.
\end{rem}

One might be tempted to think of nesting as being equivalent to inclusion of simplices, but this only works in one direction, namely:

\begin{rem}\label{rem:set_theory}
Let $\Delta,\Delta'$ be simplices of $X$. If $\Delta\subseteq \Delta'$, then $[\Delta']\nest [\Delta]$.
\end{rem}

Notice that Definition~\ref{defn:combinatorial_HHS}.\eqref{item:cHHS_join} can be rephrased as follows:

\begin{enumerate}[start=3]
 \item Whenever $\Delta$ and $\Sigma$ are non-maximal simplices for which there exists 
a non-maximal simplex $\Gamma$ such that $[\Gamma]\nest[\Delta]$, $[\Gamma]\nest[\Sigma]$, and $\diam(\fontact  
(\Gamma))\geq \delta$, then there exists a simplex $\Pi$ such that:
\begin{itemize}
    \item $\Pi$ extends $\Sigma$ (in particular $[\Pi]\nest[\Sigma]$);
    \item $[\Pi]\nest [\Delta]$;
    \item all $[\Gamma]$ as above satisfy $[\Gamma]\nest[\Pi]$.
\end{itemize}
\end{enumerate}

Our next goal is to define projections from $\mathcal W$ to $\fontact ([\Delta])$ for $[\Delta]\in\frakS$.

\begin{defn}[Projections]\label{defn:projections}
Let $(X,W,\delta,n)$ be a combinatorial HHS.  

Fix $[\Delta]\in\frakS$ and define a map $\pi_{[\Delta]}:W\to 2^{\fontact ([\Delta])}$ as follows. Let
$$p:Y_\Delta\to2^{\fontact ([\Delta])}$$
be the coarse closest point projection, i.e. 
$$p(x)=\{y\in\fontact ([\Delta]):\dist_{Y_\Delta}(x,y)\le\dist_{Y_\Delta}(x,\fontact ([\Delta]))+1\}.$$

As explained in \cite[Remark 1.17]{BHMS}, $p$ is a coarsely Lipschitz, coarse map. Roughly, this is because either $\diam(\fontact ([\Delta]))\le \delta$, and the conclusion is immediate, or $Y_\Delta$ is uniformly hyperbolic, and therefore $p$ is the coarse retraction onto the quasiconvex subset $\fontact ([\Delta])$.

Now suppose that $w$ is a vertex of $\mathcal W$, so $w$ corresponds to a unique simplex $\Sigma_w$ of $X$. By \cite[Lemma 1.15]{BHMS}, the intersection $\Sigma_w\cap Y_\Delta$ is non-empty and has diameter at most $1$. Define $$\pi_{[\Delta]}(w)=p(\Sigma_w\cap Y_\Delta).$$

We have thus defined  $\pi_{[\Delta]}:W^{(0)}\to 2^{\fontact ([\Delta])}$. If $v,w\in W$ are joined by an edge $e$ of $\mathcal W$, then $\Sigma_v,\Sigma_w$ are joined by edges in $\duaug{X}{\mathcal W}$, and we let
$$\pi_{[\Delta]}(e)=\pi_{[\Delta]}(v)\cup\pi_{[\Delta]}(w).$$

Now let $[\Delta],[\Delta']\in\frakS$ satisfy $[\Delta]\transverse[\Delta']$ or $[\Delta']\propnest [\Delta]$. Let $$\rho^{[\Delta']}_{[\Delta]}=p(\Sat(\Delta')\cap Y_\Delta),$$ where $p:Y_\Delta\to\fontact ([\Delta])$ is coarse closest-point projection. 

Let $[\Delta]\propnest [\Delta']$. Let $\rho^{[\Delta']}_{[\Delta]}:\fontact ([\Delta'])\to \fontact ([\Delta])$ be defined as follows.  On $\fontact ([\Delta'])\cap Y_\Delta$, it is the restriction of $p$ to $\fontact ([\Delta'])\cap Y_\Delta$. Otherwise, it takes the value $\emptyset$.
\end{defn}

We are finally ready to state the main theorem of \cite{BHMS}:
\begin{thm}[HHS structures for $X$--graphs]\label{thm:HHS_links}
Let $(X,W)$ be a combinatorial HHS. Let $\frakS$ be as in Definition~\ref{defn:simplex_equivalence}, define nesting and orthogonality relations on $\frakS$ as in Definition~\ref{defn:nest_orth}, let the associated hyperbolic spaces be as in Definition~\ref{defn:combinatorial_HHS}, and define projections as in Definition~\ref{defn:projections}. 

Then $(W,\frakS)$ is a hierarchically hyperbolic space, and the HHS constants only depend on $\delta,n$ as in Definition 
\ref{defn:combinatorial_HHS}. 
\end{thm}

The aim of the present paper is, morally, to establish a “converse” of the previous result, by showing that any HHS satisfying reasonable hypotheses has a hierarchically hyperbolic structure that comes from a combinatorial HHS.

\section{Combinatorial hyperbolicity from hierarchical hyperbolicity}\label{sec:converse}

Fix a hierarchically hyperbolic space $(\cuco Z,\frakS)$. The goal of this section is to construct a combinatorial HHS structure $(X,\cuco Z)$ for the space $\cuco Z$.  The exact statement is Theorem~\ref{thm:main}, which will require the additional mild assumptions on $(\cuco Z,\frakS)$ that we now present.

\subsection{(Weak) wedges}
The following property was first articulated in~\cite{BerlaiRobbio}. It is a fairly natural requirement, satisfied by all reasonable naturally occurring examples.

\begin{property}[Wedges]\label{defn:wedge_property}
 The HHS $({\cuco Z},\frakS)$ has \emph{wedges} if for all $U,V\in\frakS$, one of the following holds:
 \begin{itemize}
  \item there exists a unique $\nest$--maximal $T\in\frakS$ such that $T\nest U$ and $T\nest V$, and we write $T=U\wedge V$;
  \item there does not exist $T\in\frakS$ with $T\nest U$ and $T\nest V$, and we formally write $U\wedge V=\emptyset$.
 \end{itemize}
\end{property}

What we will actually need is the following weak version of the wedge property:
\begin{property}[Weak wedges]\label{defn:weak_wedge}
 The HHS $(\cuco Z,\frakS)$ has \emph{weak wedges} if for all $U,V\in\frakS$, one of the following holds:
 \begin{enumerate}
  \item\label{item:ww1} there exists a $T\in\frakS$ such that $T\nest U$, $T\nest V$ and whenever $W\in \frakS$ is a $\nest$-minimal domain that is nested in both $U$ and $V$ then $W\nest T$;
  \item there does not exist $T\in\frakS$ with $T\nest U$ and $T\nest V$.
 \end{enumerate}
\end{property}

\begin{rem}\label{rem:find_THE_weak_wedge}
If $(\cuco Z,\frakS)$ has weak wedges and $U,V\in\frakS$ share a common nested domain then we can find a unique $\bar T$ satisfying the properties of Item~\eqref{item:ww1} which is $\nest$-minimal among all domains with the same properties. To do so, let $\mathfrak T= \{T_i\}_{i\in I}$ be the family of domains which are nested in both $U$ and $V$ and that contain every $\nest$-minimal domain $W$ which is nested in both $U$ and $V$. For any two $T_i,T_j\in \mathfrak T$ there exists some $T\nest T_i$, $T\nest T_j$ which satisfies the properties in Definition~\ref{defn:weak_wedge} for $T_i$ and $T_j$. Hence $T$ is again an element of $\mathfrak T$, since it must contain all $\nest$-minimal domains which are nested in both $T_i$ and $T_j$ (and therefore in both $U$ and $V$). Now, there must be an element $\bar T\in \mathfrak T$ which is nested in all elements of the family, because otherwise, by the previous observation, we could find an infinite chain $T_1\sqsupsetneq T_2\sqsupsetneq\ldots$ which would contradict the finite complexity of the HHS. Hence, we say that $\bar T$ is \emph{the} weak wedge of $U$ and $V$, and denote it by $U\wedge_\text{min}V$.
\end{rem}

\subsection{Clean containers}
The second main property was first articulated in~\cite{ABD}, and is still very natural.

\begin{property}[Clean containers]\label{property:clean_containers}
The HHS $(\cuco Z,\frakS)$ has \emph{clean containers} if the following holds. Let $T\in\frakS$. Suppose that $U\propnest T$ and 
$$A=\{V\in\frakS:V\propnest T,V\orth U\}\neq \emptyset.$$
Then there exists $U^\orth_T\propnest T$, which we call the \emph{orthogonal complement} of $U$ inside $T$, such that $U^\orth_T\orth U$ and $V\nest U^\orth_T$ for each $V\in A$.
\end{property}

In other words, the clean container property states that there exists a unique container for $U$ inside $T$, as in Definition~\eqref{item:dfs_orthogonal}, and it is actually orthogonal to $U$. When the ambient domain $T$ coincides with $S$ we simply write $U^\orth:=U^\orth_S$.

\subsection{Orthogonals for non-split domains}
For the next property, which is the first real requirement on $(Z,\frakS)$, we first recall a definition from \cite{Asymptotic_dim}:
\begin{defn}[Friendly] Let $V,W\in\frakS$. Then $W$ is \emph{friendly} to $V$ if $W\nest V$ or $W\orth V$. 
\end{defn}
Notice that, as often happens in life, friendship is not always a symmetric relation.

\begin{defn}[Split]\label{defn:split}
    A domain $U\in\frakS$ is \emph{split} if there exists a $\nest$-minimal domain $W\nest U$ such that, for every $V\nest U$, we have $W\nest V$ or $W\orth V$. We say that $W$ is a \emph{Samaritan} for $U$, since it is friendly to every other $V\nest U$.
\end{defn}
An example of a split domain in the usual HHS structure of the mapping class group is as follows. Let $U$ be a subsurface given by the disconnected union of an annulus $W$ and another subsurface. Then $W$ is a Samaritan for $U$, since any subsurface $V\nest U$ which does not contain $W$ must be nested in $U- W$. We postpone the details to Subsection~\ref{sec:mcg}.

\begin{rem}\label{rem:middle_in_split}
Notice that, if $W$ is a Samaritan for $U$ and $W\nest U'\nest U$, then by definition $U'$ is also split with Samaritan $W$.
\end{rem}

\begin{rem}
If $U$ is $\nest$-minimal then it is trivially split, since it coincides with its unique Samaritan.\end{rem}

\begin{property}\label{property:orthogonal_for_non_split} A hierarchically hyperbolic space has \emph{orthogonals for non-split domains} if for every two domains $U\propnest V$ either $U$ is split or there exists $W\propnest V$ such that $W\orth U$.
\end{property}

This property feels like it should be related to the unbounded product property from \cite{ABD}, but in practice the connection seems to not be very strong.

% \begin{rem}
%     \marginpar{\textcolor{red}{The referee suggested adding a comparison with the unbounded products property from [ABD21] (=bounded domain dichotomy+whenever $F_U$ is unbounded then so is $E_U$), but the two properties seem quite different to me. Should we write anything?}}
% \end{rem}

\subsection{(Everywhere) dense product regions}
Finally, for our construction to work, we must require that every coordinate space $\fontact U$ can be reconstructed from the projections coming from the domains nested inside $U$.

\begin{property}[DPR]\label{property:dpr}
    A hierarchically hyperbolic space $(\cuco Z,\frakS)$ has \emph{dense product regions} if there exists a constant $M_0$ such that, whenever $U\in\frakS$ is not $\nest$-minimal, for any $p\in \fontact U$ there exists $V\propnest U$ such that $\dist_U(p,\rho^V_U)\le M_0$.
\end{property}

Notice that, up to choosing a bigger constant $M_0$, we may always find a domain $V\propnest U$ as in the previous property which is also $\nest$-minimal, since if $V\propnest V'\propnest U$ then $\dist_U(\rho^{V}_U, \rho^{V'}_U)\le 10E$ by the consistency axiom~\eqref{item:dfs_transversal}.

\begin{rem}
    In the study of HHS, it is often the case that one can remove minimal domains with bounded coordinate spaces. This is, however, not possible in our setting, since we might need these domains to witness DPR.
\end{rem}

Actually, we will use a seemingly stronger, yet equivalent version of the DPR property, which we now state.

\begin{property}[EDPR]\label{property:EDPR}
A hierarchically hyperbolic space $(\cuco Z,\frakS)$ has \emph{everywhere dense product regions} if there exists a constant $C_0$ such that the following holds. For every $U\in \frakS$ and every $x\in F_U$ there exists a maximal family $V_1,\ldots, V_k\propnest U$ of $\nest$-minimal, pairwise orthogonal domains such that $x$ is $C_0$-close to the relative product region $P^U_{\{V_i\}_{i=1,\ldots, k}}$.
\end{property}

Clearly, property~\eqref{property:EDPR} is stronger than property~\eqref{property:dpr}. The converse also holds:% if one allows a single change in the HHS constant $E$ from Definition~\ref{defn:HHS}, but, remarkably, without changing the rest of the HHS structure, including the domain set and all projections. 

\begin{lemma}\label{lem:dpr_grows}
    A HHS with the DPR property~\eqref{property:dpr} also has the EDPR property~\eqref{property:EDPR}.
\end{lemma}

\begin{proof} Let $(\cuco Z, \frakS)$ be a HHS with the DPR property, with some constant $M_0$, and let $E$ be a HHS constant for $(\cuco Z, \frakS)$. We shall prove the existence of a constant $C_0=C_0(M_0,E)$ such that $(\cuco Z, \frakS)$ also enjoys the EDPR property with respect to $C_0$. 

We will prove the Lemma by induction on the \emph{level} $l$ of $U$, that is, the maximum $k$ such that there exists a chain $U_0\propnest\ldots\propnest U_k=U$. If $l=0$ then $U$ is $\nest$-minimal and the EDPR property~\eqref{property:EDPR} clearly holds. 
    
Now suppose the theorem holds for every domain of level strictly less than $l$, and let $U$ be a domain of level $l$. Before going on with the proof we recall some definitions from \cite{HHS_II}.

We will say that the collection $\mathfrak{U}$ of elements of $\frakS_U$ is \emph{totally orthogonal} if any pair of distinct elements of $\mathfrak{U}$ are orthogonal. Given a totally orthogonal family $\mathfrak{U}$, we say that $W\nest U$ is \emph{$\mathfrak{U}$–generic} if there exists $V\in\mathfrak{U}$ so that $W$ is not orthogonal to $V$. 

Now fix $x\in F_U$ that we want to realise with minimal domains. A totally orthogonal collection $\mathfrak{S}$ is \emph{$C$–good} if any $E$–partial realisation point $y$ for $\mathfrak{U}$, as defined in the partial realisation axiom~\eqref{item:dfs_partial_realisation}, has the property that for each $W\nest U$ we have $\dist_W(x_w,y_w)\le C$. Notice that our goal is to find a maximal family $\mathfrak{U}$ which is made of minimal supports and $C$-good for some uniform constant $C$. Notice moreover that, if $\mathfrak{U}$ is $C$-good but not maximal, then we can add $\nest$-minimal domains to $\mathfrak{U}$ and complete it to a maximal totally orthogonal family. The latter will again be $C$-good, because a partial realisation point for the larger family is also a partial realisation point for $\mathfrak{U}$.

A totally orthogonal collection $\mathfrak{U}$ is \emph{$C$–generically good} if any $E$–partial realisation point $y$ for $\mathfrak{U}$ has the property that for each $\mathfrak{U}$–generic $W$ we have $\dist_W(x_w,y_w)\le E$.

We allow that $\mathfrak U$ can be empty. In this case, we say that a $C$-partial realisation point for $\emptyset$ is simply a point $y$ such that $\dist_U(x_U, y_U)\le C$. Notice that no $W$ is $\emptyset$–generic.

Now the following lemma holds:

\begin{lemma}\label{lem:totally_orth_growing}
    For every $C\ge 100 E^3$ the following holds. Let $\mathfrak U$ be totally orthogonal and $C$–generically good but not $C$–good. Then there exists a totally orthogonal, $10C$–generically good collection $\mathfrak{U}'$ with $\mathfrak{U}\subsetneq \mathfrak{U}'$, obtained by adding $\nest$-minimal domains.
\end{lemma}

This fact is proven exactly as \cite[Lemma 3.3]{HHS_II}, whose proof runs verbatim in our case. The only difference is that our inductive hypothesis, which replaces that of \cite[Theorem 3.1]{HHS_II}, allows us to assume that the additional domains are all $\nest$-minimal.

Now we can prove Lemma~\ref{lem:dpr_grows}. Recall that we want to realise a point $x\in F_U$. If $\mathfrak{U}=\emptyset$ is already $M_0$-good we can choose any maximal family of pairwise orthogonal, minimal domains $V_1,\ldots, V_n\nest U$ such that $\dist_U(x_U, \rho^{V_1}_U)\le M_0$ (whose existence is granted by the DPR property~\eqref{property:dpr}). Then any realisation point $y$ for $\{(V_i, x_{V_i})\}$ is also a realisation point for $\mathfrak{U}=\emptyset$, since $\dist_U(x_U, y_U)\le M_0$ by construction, and therefore $x$ and $y$ are $M_0$-close in every coordinate space.

Otherwise, since no $W$ is $\emptyset$-generic, we can apply Lemma~\ref{lem:totally_orth_growing} and find a larger $\mathfrak U_1$ which is $10 M_0$-generically good. If $\mathfrak U_1$ is $10 M_0$-good we can complete it to a maximal family of pairwise orthogonal, $\nest$-minimal domains which is again $10 M_0$-good, and we are done. Otherwise, we can repeat the process with $\mathfrak U_1$. Since there is a bound on the cardinality of totally orthogonal sets, in finitely many steps we necessarily get a good totally orthogonal set made of minimal supports, and this concludes the proof.
\end{proof}

\subsection{The main theorem}

\begin{rem}[Normalisation preserves our hypotheses]\label{rem:normalised_HHS}
Before stating the main theorem we notice that, if $(\cuco Z,\frakS)$ has one of the properties defined in the previous section, then so does the normalised structure, as defined in Remark~\ref{rem:normalisation}. This is because the new structure has the same set of domains $\frakS$ with the same relations of nesting and orthogonality, thus all combinatorial assumptions on the domain set (wedges, clean containers and property~\eqref{property:orthogonal_for_non_split}) are preserved under the normalisation procedure. Moreover, the DPR property~\eqref{property:dpr} still holds as well, since the coarse closest point projection is coarsely Lipschitz. Hence, our Assumption~\ref{assump:normalised} that the HHS structure is normalised does not lose any generality.
\end{rem}

\begin{thm}\label{thm:main}
Let $(\cuco Z,\frakS)$ be a normalised hierarchically hyperbolic space with weak wedges, clean containers, the orthogonals for non-split domains property~\eqref{property:orthogonal_for_non_split} and the DPR property~\eqref{property:dpr}. There exists a combinatorial HHS $(X,\mathcal{W})$ such that $\cuco Z$ is quasi-isometric to $\mathcal{W}$.

Moreover, let $G$ be a finitely generated group which acts on $(\cuco Z,\frakS)$ by automorphisms. Then $G$ acts on $(X,\mathcal{W})$, and the quasi-isometry $f:\,\mathcal{W}\to \cuco Z$ is coarsely $G$-equivariant. 
\end{thm}

\begin{proof}[Outline of the proof of Theorem~\ref{thm:main}]
The graphs $X$ and $W$ are constructed in Subsections~\ref{subsec:minimal_orth_graph_and_blow_up} and~\ref{subsec:HHS_W}, respectively.  
The four conditions of~\ref{defn:combinatorial_HHS} are verified in Subsection~\ref{sec:verify_axioms}, and more precisely:
\begin{itemize}
\item Condition 1 is Corollary~\ref{cor:finite_complexity_HHS_blow_up};
\item The two parts of Condition 2 are proved in Subsections~\ref{subsec:hyperbolic_links} and~\ref{subsec:qi_embedded_links};
\item Condition 3 is Lemma~\ref{lem:simplicial_wedge_property}, which is implied by Lemma~\ref{lem:intersection_of_links};
\item Condition 4 is Lemma~\ref{lem:edges_in_link}.
\end{itemize}
In Definition~\ref{defn:realisation_map} we define a map $f:\mathcal{W}\to \cuco Z$, which we prove to be a quasi-isometry in Lemma~\ref{lem:W_Z_qi}. Finally, the “moreover” part of the statement is proved in Lemma~\ref{thm:G-equivariance}.
\end{proof}

\begin{rem}
    Towards relative HHS versions of our results we remark the following. The proof will feature $\nest$-minimal domains only in Lemma \ref{lem:hyperbolic_links} and Claim \ref{claim:Y_Delta_for_minimal}. The proof of the former shows that certain augmented links are quasi-isometric to the coordinate spaces of $\nest$-minimal domains, with hyperbolicity of these only used to then conclude that those augmented links are hyperbolic. The proof of the latter does not use hyperbolicity of coordinate spaces of $\nest$-minimal domains.
\end{rem}

\section{Construction of the combinatorial HHS}\label{sec:construction}
In this Section we construct the pair $(X,\mathcal{W})$ and the map $f:\,\mathcal{W}\to \cuco Z$.
\begin{rem}\label{rem:only_clean}
    The construction of $(X,\mathcal{W})$ and $f$ will only require $(\cuco Z,\frakS)$ to have clean containers. The other hypotheses of Theorem~\ref{thm:main} will be needed later, to ensure that $(X,\mathcal{W})$ is actually a combinatorial HHS and that $f$ is a quasi-isometry.
\end{rem}

\subsection{Moral compass}
Before going into the actual details, we explain the idea of the construction, and why it should work (at least morally).

First, we consider the graph $\bar X$ whose vertices are all $\nest$-minimal domains of $\frakS$ (see Definition~\ref{defn:blow-up}). Now, whenever $U\in\frakS$ is not $\nest$-minimal, its coordinate space $\fontact U$ can be reconstructed by just looking at the projections $\rho^V_U$ coming from the $\nest$-minimal domains, by the dense product regions property~\ref{property:dpr}, and such projections are close whenever the $\nest$-minimal domains are orthogonal, by Lemma~\ref{lem:close_proj_of_orthogonals}. Hence, in a sense, the graph $\bar X$ will contain all information about the HHS structure coming from the non-$\nest$-minimal domains. 

However, $\bar X$ does not see the coordinate spaces of $\nest$-minimal domains. Therefore, for every vertex $V$ of $\bar X$ we consider its coordinate space $\fontact V$, which we may assume to be a simplicial graph up to quasi-isometry (see for example \cite[Lemma 3.B.6]{cornulier2014metric}), and we replace $V$ with the cone over the $0$-skeleton $\fontact V^{(0)}$. This way, the apex of the cone, call it $v_V$, will have $\fontact V^{(0)}$ inside its link, and after adding the right $\mathcal{W}$-edges we will be able to see $\fontact V$ as the augmented link of some simplex. Call $X$ the graph obtained after this “blow-up” procedure (again, see Definition~\ref{defn:blow-up}).

Hence, a vertex of $\mathcal{W}$, which corresponds to a maximal simplex of $X$, is the data of a collection $V_1,\ldots, V_k$ of pairwise orthogonal and $\nest$-minimal domains, and a point $x_i\in \fontact V_i$ for all $i=1, \ldots,k$. Such a collection $\{(V_i,x_i)\}$ admits a unique realisation point, in the sense of Theorem~\ref{thm:realisation}, thus we can define a map $f$ between vertices of $\mathcal{W}$ and points in $\cuco Z$ (see Definition~\ref{defn:realisation_map}).

Regarding the edges of $\mathcal{W}$, morally we would like to say that, if $\Sigma=\{(V_i,x_i)\}_{i=1,\ldots, k}$ and $\Delta=\{(W_j,x_j)\}_{j=1,\ldots, l}$ are two maximal simplices, then they are $\mathcal{W}$-adjacent if and only if their realisation points are close, so that $f$ is a quasi-isometry almost by definition. In turn, such realisation points are close if and only if their coordinates are close in every coordinate space. For some technical reasons (mainly appearing in the proof of Lemma~\ref{lem:edges_in_link}), the exact definition of the edges of $\mathcal{W}$ will also take into account the supports $\{V_i\}_{i=1,\ldots, k}$ and $\{W_j\}_{j=1,\ldots, l}$ of the two simplices: the more these supports intersect, the further we allow the coordinates of the realisation points to be (see Definition~\ref{defn:W-edge}).

\subsection{The minimal orthogonality graph and its blow-up}\label{subsec:minimal_orth_graph_and_blow_up}
The first step in the proof of Theorem~\ref{thm:main} is to construct the simplicial complex $X$, which will heavily depend on both $\cuco Z$ and the actual HHS structure. For the purpose of this subsection we do not need to assume any property on $(\cuco Z,\frakS)$.

Let $\frakS_{\text{min}}$ be the set of $\nest$--minimal elements of $\frakS$.  Let $\bar X^{(1)}$ be the graph with vertex-set $\frakS_{\text{min}}$, with $U,V\in\frakS_{\text{min}}$ joined by an edge when $U\orth V$. Let $\bar X$ be the flag complex on $\bar X^{(1)}$.

For each $U\in\frakS_{\min}$ we can assume, up to quasi-isometry, that $\fontact U$ is a graph. Thus let $L(U)$ be the cone on $(\fontact U)^{(0)}$, and denote by $v_U$ the cone-vertex.  Let $X^{(1)}$ be the graph formed from $\bigsqcup_{U\in\frakS_{\text{min}}}L(U)$ by joining each vertex of $L(U)$ to each vertex of $L(V)$ whenever $U\orth V$ (i.e. whenever $U,V$ are adjacent in $\bar X^{(1)}$).  Let $X$ be the flag complex on $X^{(1)}$.

\begin{figure}[htp]
    \centering
    \includegraphics[width=\textwidth]{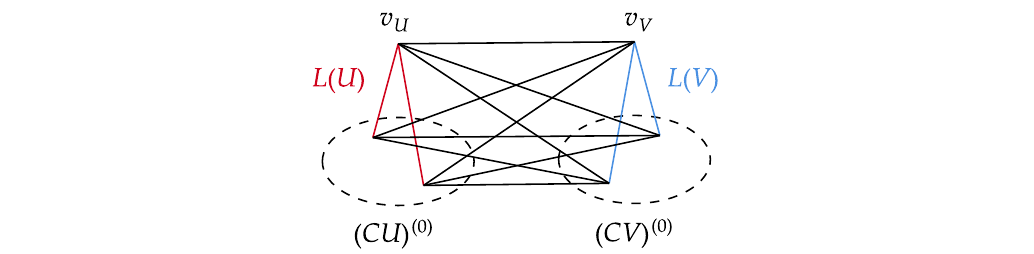}
    \caption{$U$ and $V$ are orthogonal, $\nest$-minimal domains. Therefore, the cone over $(\fontact U)^{(0)}$ (here, in red) and the cone over $(\fontact V)^{(0)}$ (in blue) form a join. Notice that any two points in $(\fontact U)^{(0)}$ are \emph{not} adjacent.}
    \label{fig:UorthV}
\end{figure}

\begin{defn}[Blow-up]\label{defn:blow-up}
Define the retraction $p:\,X\to\bar X$ by collapsing each subcomplex $L(U)$ to the vertex $U$. We will refer to $\bar X$ as the \emph{minimal orthogonality graph} of the structure $(\cuco Z,\frakS)$, and to $X$ as the \emph{blow-up} of $\bar X$.
\end{defn}

For each simplex $\Delta$ of $X$, let $\bar\Delta=p(\Delta)$ be the image simplex in $\bar X$. We will say that $\Delta$ is \emph{supported} in $\bar\Delta$.

Given a simplex $\Delta$ of $X$ and a vertex $U\in\bar \Delta$, let $\Delta_U=\Delta\cap p^{-1}(U)$. Note that $\Delta_U$ is either a vertex of $\Delta$ or an edge of $\Delta$.  Moreover, we have 
$$\Delta=\bigstar_{U\in\bar(\Delta)^{(0)}}\Delta_U.$$
A careful inspection of the construction yields the following (compare with \cite[Lemma 4.12]{ELTAG_HHS}):
\begin{lemma}[Decomposition of links]\label{lem:decomposition_of_links}
Let $\Delta$ be a simplex of $X$. Then
$$\link(\Delta)=p^{-1}(\link_{\bar X}(\bar\Delta))\star (\bigstar_{U\in\bar\Delta^{(0)}}\link_{p^{-1}(U)}(\Delta_U)).$$
\end{lemma}

\begin{cor}\label{cor:bounded_links}
Let $\Delta$ be a simplex of $X$. Then one of the following holds:
\begin{enumerate}
 \item \label{cor:bounded_links1} $\link(\Delta)$ is either a single vertex or a non-trivial join.
 \item \label{cor:bounded_links2} For each $U\in\bar\Delta^{(0)}$, we have that $\Delta_U$ is an edge.
 \item \label{cor:bounded_links3}  The simplex $\bar\Delta$ is maximal inside $\bar X$ and the following holds. There exists $U\in\bar\Delta^{(0)}$ such that $\Delta_U=v_U$, and for every $V\in\bar\Delta^{(0)}-\{U\}$ we have that $\Delta_V$ is an edge.  
\end{enumerate}
\end{cor}

\begin{figure}[htp]
    \centering
    \includegraphics[width=\textwidth]{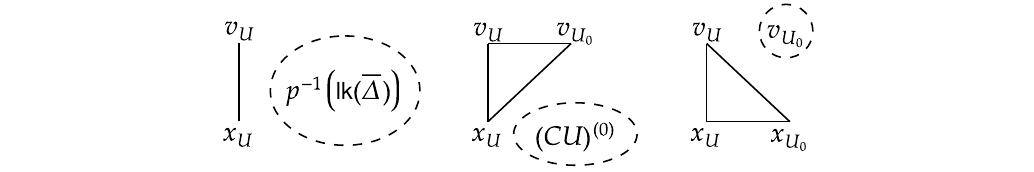}
    \caption{All cases in which $\link(\Delta)$ (represented by the dashed ellipses) is not a non-trivial join.}
    \label{fig:link_of_Delta}
\end{figure}

\begin{proof}
    If at least two of the terms of the join from Lemma~\ref{lem:decomposition_of_links} are non-empty then $\link(\Delta)$ is a non-trivial join, and Item~\eqref{cor:bounded_links1} holds. Thus suppose that exactly one of the terms of the join is non-empty, so that $\link(\Delta)$ coincides with that term.
    
    If $p^{-1}(\link_{\bar X}(\bar\Delta))\neq\emptyset$ then for every $U\in\bar\Delta^{(0)}$ we have $\link_{p^{-1}(U)}(\Delta_U))=\emptyset$. Hence $\Delta_U$ is an edge, and Item~\eqref{cor:bounded_links2} is satisfied (see the picture on the left in Figure~\ref{fig:link_of_Delta}).
    
    Otherwise, suppose that $p^{-1}(\link_{\bar X}(\bar\Delta))=\emptyset$ (that is, that $\bar\Delta$ is a maximal simplex) and that $\link_{p^{-1}(U)}(\Delta_U))=\emptyset$ for all domains except one, call it $U_0$. In particular, we have that 
    $$\link(\Delta)=\link_{p^{-1}(U_0)}(\Delta_{U_0})).$$
    If $\Delta_{U_0}=v_{U_0}$ then we are in the case of Item~\eqref{cor:bounded_links3} (see the central picture in Figure~\ref{fig:link_of_Delta}). Otherwise $\Delta_{U_0}$ is a point in $(\fontact U)^{(0)}$, and therefore $\link_{p^{-1}(U_0)}(\Delta_{U_0}))=v_{U_0}$ is a single point, thus again giving Item~\eqref{cor:bounded_links1} (see the picture on the right in Figure~\ref{fig:link_of_Delta}).
\end{proof}

\begin{rem}\label{rem:finite_dim}
By definition, the maximum cardinality of a simplex in $X$ is twice the maximum cardinality of a family of $\nest$-minimal, pairwise orthogonal domains, which in turn is bounded above by the \emph{complexity} of the HHS structure by \cite[Lemma 2.1]{HHS_II}. Thus $X$ has finite dimension.
\end{rem}

\subsection{Edges between maximal simplices}\label{subsec:HHS_W}
Now we define the graph $\mathcal{W}$ for the combinatorial HHS $(X,\mathcal{W})$. The construction of this subsection will only require $(\cuco Z,\frakS)$ to have clean containers.

Let $\mathfrak M(X)$ be the set of maximal simplices of $X$. The vertex set of the graph $\mathcal{W}$ is $\mathfrak M(X)$. Now, given $\sigma\in\mathfrak M(X)$, let $U_1,\ldots,U_n\in\frakS_{\text{min}}$ be the vertices of $p(\sigma)$ (a maximal collection of pairwise orthogonal elements), so that $\sigma=\bigstar_{i=1}^n\{v_{U_i},x_i\}$, where 
$x_i\in\fontact U_i$ and $\{v_{U_i},x_i\}$ is the edge of $X$ joining $v_{U_i}$ to $x_i$. We will refer to the domains $U_1,\ldots,U_n\in\frakS_{\text{min}}$ as the \emph{support} of $\sigma$, and to the points $x_i\in\fontact U_i$ as the \emph{coordinates} prescribed by $\sigma$.

\begin{defn}\label{defn:b_sigma}
Let $\sigma=\bigstar_{i=1}^n\{v_{U_i},x_i\}$ be a maximal simplex of $X$. Define the tuple $b(\sigma)=(b_V)_{V\in\frakS}\in\prod_{V\in\frakS}2^{\fontact V}$ as follows. Let $V\in\frakS$. Note that $V$ is not properly nested in any $U_i$, since each $U_i\in\frakS_{\text{min}}$. Moreover, $V$ cannot be orthogonal to all of the $U_i$; otherwise, we could choose $V'\nest V$ with $V'\in\frakS_{\text{min}}$ and observe that $V'\orth U_i$ for all $i$, contradicting maximality of $\sigma$. So, for some $i\leq n$, we have either $U_i\transverse V$ or $U_i\propnest V$ or $U_i=V$. In particular, if $V\not\in\{U_1,\ldots,U_n\}$, then for some $i$ we have a set $\rho^{U_i}_V\subset\fontact V$ of diameter at most $E$.  Moreover, by Lemma~\ref{lem:close_proj_of_orthogonals} and the fact that $U_i\orth U_j$ for all $i\neq  j$, we have $\diam(\bigcup_{i\leq n,U_i\not\bot V}\rho^{U_i}_V)\leq 10E$.

Then:
\begin{itemize}
 \item If $V=U_i$ for some $i\leq n$ set $b_V=x_i$. This is why we called $x_i$ the coordinate prescribed by $\sigma$.
 \item Otherwise, let $b_V\in\fontact V$ be $b_V=\bigcup_{i\leq n,U_i\not\bot V}\rho^{U_i}_V$, which is a non-empty set of diameter at most $10E$.
\end{itemize}
\end{defn}

\begin{defn}[Co-level]\label{defn:colevel}
Let $U\in\frakS$. We define the \emph{co-level} of $U$, denoted by $\text{co-lv}(U)$, as the maximum $k$ such that there exists a chain of the form $U=U_0\propnest\ldots\propnest U_k=S$.
\end{defn}

Notice that if $U\propnest V$ then $\text{co-lv}(U)\gneq\text{co-lv}(V)$. Moreover, the maximum co-level is $n-1$, where $n$ is the complexity of the HHS structure.

\begin{defn}[Orthogonal complement]
   Let $\bar\Delta=\{U_1,\ldots, U_k\}$ be a non-empty simplex inside $\bar X$ and let $V$ be a domain such that $U_i\nest V$ for all $i$. If there exists $T\propnest V$ that is orthogonal to all $U_i$ then the \emph{orthogonal complement} of $\bar\Delta$ inside $V$, which we will denote as $\bar\Delta^\orth_V$, is constructed inductively as follows:
   \begin{itemize}
       \item set $\{U_1\}^\orth_V$ as the orthogonal complement of $U_1$ inside $V$, as in Definition~\ref{property:clean_containers};
       \item If $\{U_1,\ldots,U_i\}^\orth_V$ has already been defined, then $\{U_1,\ldots,U_{i+1}\}^\orth_V$ is the orthogonal complement of $U_{i+1}$ inside $\{U_1,\ldots,U_i\}^\orth_V$.
   \end{itemize}
   If $V=S$ we denote the orthogonal complement of $\bar\Delta$ in the maximal domain simply as $\bar\Delta^\orth$.
\end{defn}

Notice that the definition is independent of the order of the vertices, because by construction $\bar\Delta^\orth_V$ is also the unique $\nest$-maximal element $T$ which is nested inside $V$ and is orthogonal to all vertices of $\bar\Delta$.

\begin{notation}\label{notation:exceptional_orth}
With an innocent abuse of notation, we could say that:
\begin{itemize}
    \item the orthogonal complement of the empty simplex is the $\nest$-maximal element $S$, and $\text{co-lv}(S)=0$;
    \item the orthogonal complement of a maximal simplex is empty, and $\text{co-lv}(\emptyset)=n$.
\end{itemize}
\end{notation}

Now we can finally define the edges of $\mathcal{W}$. 
\begin{defn}[$\mathcal{W}$--edges]\label{defn:W-edge}
Let $\Sigma,\Delta$ be two maximal simplices of $X$ and let $\bar\Sigma,\bar\Delta$ be their supports. Let $b(\Sigma)=(b_U)_{U\in\frakS}$ and $b(\Delta)=(c_U)_{U\in\frakS}$. Let $n$ be the complexity of the HHS structure. Let $\lambda\ge 0$ be some constant.

Let $W=(\bar\Sigma\cap\bar\Delta)^\orth$ be the orthogonal complement of the intersection, and let $k$ be the co-level of $W$ (with the Notation~\ref{notation:exceptional_orth} for the exceptional cases). Then $\Sigma$ and $\Delta$ are $\mathcal{W}$--adjacent if and only if, for every $U\in \frakS$, $$\dist_{\fontact  U}(b_U,c_U)\le (k+1)\lambda.$$\end{defn}

In other words, the more two adjacent simplices share their supports, the further away we let them be in the fewer and fewer domains where $b_U$ and $c_U$ may actually differ. Notice that the definition depends on the constant $\lambda$, which we will later choose to be large enough. 

\subsection{The realisation map}\label{subsec:realisation_map}
Finally, we define the map $f:\,\mathcal{W}\to \cuco Z$ that will be the required quasi-isometry.
\begin{defn}[Realisation map]\label{defn:realisation_map}
For every simplex $\sigma=\bigstar_{i=1}^n\{v_{U_i},x_i\}$, the partial realisation axiom~\eqref{item:dfs_partial_realisation} provides the existence of a \emph{realisation point} $z$ such that:
\begin{itemize}
    \item $\dist_{U_i}(x_i, \pi_{U_i}(z))\le E$;
    \item for every $V\in \frakS$ such that $U_i\propnest V$ or $U_i\transverse V$ for some $i\leq n$ we have that $\dist_{V}(\rho^{U_i}_V, \pi_{V}(z))\le E$.
\end{itemize}
In other words, the coordinates of $z$ are $E$-close to the tuple $b(\sigma)$ from Definition~\ref{defn:b_sigma}, and we say that $z$ \emph{realises} $b(\sigma)$. Moreover, by the uniqueness axiom~\eqref{item:dfs_uniqueness} we have that $z$ is uniformly coarsely unique, and the bound only depends on $E$. Hence, setting $f(\sigma)=z$ gives a well-defined coarse map $f:\,\mathcal{W}\to \cuco Z$.
\end{defn}

\begin{rem}[Consistency of $b(\sigma)$]\label{rem:b_sigma_consistent}
    The existence of a $z$ that realises $b(\sigma)$ also shows that the latter is a $20E$-consistent tuple. Indeed, every coordinate of $b(\sigma)$ has diameter at most $10E$ and is $E$-close to the corresponding coordinates of $z$, which satisfies the consistency axiom~\eqref{item:dfs_transversal}.
\end{rem}

\section{Proof of the main Theorem}\label{sec:verify_axioms}
Unless otherwise stated, in this Section we will work under the following assumption:

\begin{ass}\label{ass:standing_assumption}
$(\cuco Z,\frakS)$ is a normalised hierarchically hyperbolic space with weak wedges~\eqref{defn:weak_wedge}, clean containers~\eqref{property:clean_containers}, the orthogonals for non-split domains property~\eqref{property:orthogonal_for_non_split} and the DPR property (we will work with its strong form, which is property~\eqref{property:EDPR}). Let $E$ be the HHS constant for $(\cuco Z,\frakS)$, and for every $U\in \frakS$ let $F_U$ be the space of $20E$-consistent partial tuples. Let $X$ be the graph from Subsection~\ref{subsec:minimal_orth_graph_and_blow_up}, and let $\mathcal{W}$ be the graph from Subsection~\ref{subsec:HHS_W}, whose edges depend on the constant $\lambda$ from Definition~\ref{defn:W-edge}, that we will later choose (see Lemma~\ref{lem:edges_in_link}, Claim~\ref{claim:minimal_product_regions} and Lemma~\ref{lem:W_Z_qi}). Finally, let $f:\,\mathcal{W}\to \cuco Z$ be the realisation map from Definition~\ref{defn:realisation_map}.
\end{ass}
We choose $F_U$ to be the space of $20E$-consistent tuples, in order to include the tuples of the form $b(\sigma)$ from Definition~\ref{defn:b_sigma}.

Our first goal is to prove the following, which readily implies the first half of Theorem~\ref{thm:main}:
\begin{thm}\label{thm:XWf}
Under Assumption~\ref{ass:standing_assumption} there exists $\widetilde{\lambda}\ge 0$ such that, whenever $\lambda\ge \widetilde{\lambda}$, the pair $(X,\mathcal{W})$ is a combinatorial HHS, and $f$ is a quasi-isometry.
\end{thm}

\subsection{Weak orthogonal complements}
Before getting into the proof, we develop some more technical notation and lemmas regarding the interaction between weak wedges and clean containers.

\begin{defn}[Weak orthogonal complement]\label{defn:weak_complement}
    Let $\bar\Delta=\{U_1,\ldots, U_k\}$ be a simplex inside $\bar X$. Then its \emph{weak orthogonal complement} is 
    $$\bar\Delta^\orth_{\text{min}}:=\bar\Delta^\orth\ww \bar\Delta^\orth.$$
\end{defn}

In other words, $\bar\Delta^\orth_{\text{min}}$ is the unique domain $T\in\frakS$ such that:
\begin{itemize}
    \item $T\orth U$ for every $U\in \bar\Delta$;
    \item $T$ contains every $V\in \frakS_{\text{min}}$ which is orthogonal to $\bar\Delta$ (that is, every $V\in \link_{\bar X}(\bar\Delta)$);
    \item $T$ is $\nest$-minimal among domains with the previous properties.
\end{itemize}

The weak orthogonal complement is uniquely determined by the link of $\bar\Delta$, in the following sense:
\begin{lemma}\label{lem:nested_min_orth}
    Let $\bar\Delta, \Bar\Sigma$ be two simplices inside $\bar X$. Then the following are equivalent:
    \begin{enumerate}
        \item $\link(\bar\Delta)\subseteq\link(\bar\Sigma)$;
        \item $\bar\Delta^\orth_{\text{min}}\nest\bar\Sigma^\orth_{\text{min}}$.
    \end{enumerate}
\end{lemma}

\begin{proof}
    (1 $\Rightarrow$ 2) Let $T:=\bar\Delta^\orth_{\text{min}}\ww \bar\Sigma^\orth_{\text{min}}\nest  \bar\Sigma^\orth_{\text{min}}$. Since $T\nest \bar\Delta^\orth_{\text{min}}$, we have that $T\orth U$ for every $U\in\bar\Delta$. Moreover, by definition of weak wedge, $T$ contains every $V\in \frakS_{\text{min}}$ which is nested inside both $\bar\Delta^\orth_{\text{min}}$ and $\bar\Sigma^\orth_{\text{min}}$, that is, every $V\in \link_{\bar X}(\bar\Delta)\cap \link_{\bar X}(\bar\Sigma)=\link_{\bar X}(\bar\Delta)$. Therefore, by minimality of $\bar\Delta^\orth_{\text{min}}$, we must have that $T=\bar\Delta^\orth_{\text{min}}$, hence $\bar\Delta^\orth_{\text{min}}\nest \bar\Sigma^\orth_{\text{min}}$.

    (2 $\Rightarrow$ 1) Simply notice that, if $V$ is $\nest$-minimal, then $V\in \link(\bar\Delta)$ if and only if $V\nest\bar\Delta^\orth_{\text{min}}$, by definition of $\bar\Delta^\orth_{\text{min}}$, and the same holds for $\bar\Sigma$.
\end{proof}

\begin{cor}
    Let $\bar\Delta, \Bar\Sigma$ be two simplices inside $\bar X$. Then $\link(\bar\Delta)=\link(\bar\Sigma)$ if and only if  $\bar\Delta^\orth_{\text{min}}=\bar\Sigma^\orth_{\text{min}}$.
\end{cor}

\begin{lemma}\label{lem:weak_compl=compl}
    Let $\bar\Delta$ be a simplex of $\bar X$. Then exactly one of the following holds:
    \begin{enumerate}
        \item $\bar\Delta^\orth_{\text{min}}=\bar\Delta^\orth$;
        \item $\bar\Delta^\orth_{\text{min}}$ is split. 
    \end{enumerate} 
    Moreover, in the second case $\link(\bar\Delta)$ is the cone with cone point $V\in \link(\bar\Delta)$, where $V$ is any Samaritan for $\bar\Delta^\orth_{\text{min}}$.
\end{lemma}

\begin{proof}
Clearly $\bar\Delta^\orth_{\text{min}}\nest\bar\Delta^\orth$, and if they do not coincide then the orthogonals for non-split domains property~\ref{property:orthogonal_for_non_split} states that either $\bar\Delta^\orth_{\text{min}}$ is split or there exists $V\nest\bar\Delta^\orth$ such that $V\orth \bar\Delta^\orth_{\text{min}}$. But then there exists some $\nest$-minimal domain $V'\nest V$ which lies in $\bar\Delta^\orth$ but not in $\bar\Delta^\orth_{\text{min}}$, contradicting the definition of the latter.

For the “moreover” part, just notice that, if $V$ is a Samaritan for $\bar\Delta^\orth_{\text{min}}$, then $V$ is orthogonal to any other vertex of $\link(\bar\Delta)$.
\end{proof}

\subsection{Intersection of links and finite complexity}
Now we turn to the proof of Theorem~\ref{thm:XWf}. First, we check the parts of Definition~\ref{defn:combinatorial_HHS} that depend on $X$ only.

\begin{lemma}[Verification of Definition~\ref{defn:combinatorial_HHS}.\eqref{item:cHHS_join}]
\label{lem:simplicial_wedge_property}
Let $\Sigma,\Delta$ be non-maximal simplices of $X$ and suppose that there exists a non-maximal simplex $\Gamma$ such that $[\Gamma]\nest[\Sigma]$, $[\Gamma]\nest[\Delta]$ and $\diam(\fontact ([\Gamma]))\ge 3$. Then there exists a non-maximal simplex $\Pi$ which extends $\Sigma$ such that $[\Pi]\nest[\Delta]$ and all $\Gamma$ as above satisfy $[\Gamma]\nest[\Pi]$.
\end{lemma}

Arguing as in the proof of \cite[Theorem 6.4]{BHMS} (more precisely, at the beginning of the paragraph named “\textbf{$(X,W)$ is a combinatorial HHS}”), one sees that Lemma~\ref{lem:simplicial_wedge_property} is implied by the following, which is \cite[Condition 6.4.B]{BHMS} there:

\begin{lemma}\label{lem:intersection_of_links}
    Under Assumption~\ref{ass:standing_assumption}, let $\Sigma,\Delta$ be non-maximal simplices of $X$. Then there exist two (possibly empty or maximal) simplices $\Pi,\Psi\subset X$ such that $\Sigma\subset \Pi$ and 
    $$\link(\Sigma)\cap\link(\Delta)=\link(\Pi)\star\Psi.$$
\end{lemma}

\begin{proof}[Proof of Lemma~\ref{lem:intersection_of_links}]
We subdivide the proof into two major steps.
\par\medskip

\textbf{Finding the support of the extended simplex}: Let $\bar\Sigma$ and $\bar\Delta$ be the supports of $\Sigma,\Delta$, respectively, and let $\bar\Sigma^\orth, \bar\Delta^\orth\in \frakS$ be their orthogonal complements. Let $\bar\Phi=\bar \Delta\cap \link(\bar \Sigma)$, and let $Y_0$ be the orthogonal complement of $\bar\Phi$ inside $\bar\Sigma^\orth$, that is, $Y_0=(\bar\Sigma\star\bar\Phi)^\orth$. Finally, set $\bar\Theta=\bar\Psi=\emptyset$. We will progressively add vertices to these simplices, which will form the supports of the simplices $\Pi$ and $\Psi$ we are looking for.

If $\bar\Sigma^\orth$ and $\bar\Delta^\orth$ have no common nested domain we formally set $W_0=\emptyset$. Otherwise, let $W_0=\bar\Sigma^\orth\wedge_\text{min} \bar\Delta^\orth$ be the weak wedge of the orthogonal complements. Since by construction $W_0$ is nested in  $\bar\Delta^\orth$, every vertex of $\bar\Phi$ is orthogonal to $W_0$, Thus $W_0\nest Y_0$, since $W_0$ is also nested in $\bar\Sigma^\orth$ and $Y_0$ is a clean container. 

Now we do the following procedure, which is divided into three parts.
\par\medskip

\textit{Part 1}: If $W_0=\emptyset$ or $W_0$ is non-split then we can set $W'=W_0$ and $Y_0'=Y_0$ and skip to Part 2. Otherwise, there exists a Samaritan $U_1\nest W_0$ such that $U_1\orth V$ for every other $\nest$-minimal domain $V\propnest W_0$, and we add $U_1$ to $\bar\Psi$. Then let $W_1=\{U_1\}^\orth_{W_0}$ (which might be empty), and similarly let $Y_1=\{U_1\}^\orth_{Y_0}$. Clearly $W_1\nest Y_1$, thus if $W_1$ is again split we can repeat this argument with $Y_1$ and $W_1$. This procedure, which adds one vertex at a time to $\bar\Psi$, must end after at most $n$ steps by the finite complexity axiom~\eqref{item:dfs_complexity}, since every new $W_i$ is properly nested into $W_{i-1}$ for all $i$. Moreover, this procedure stops when $W'=\bar\Psi^\orth_{W_0}$ is either empty or non-split. Set $Y_0'=\bar\Psi^\orth_{Y_0}$ (which again might be empty).
\par\medskip

\textit{Part 2}: Now, if $W'=\emptyset$ we choose a simplex $\bar\Theta=\{V_1, \ldots, V_k\}$ of pairwise orthogonal, $\nest$-minimal domains inside $Y'_0$, and we skip to Part 3. Otherwise $W'\nest W_0$ is a non-split domain which is nested inside $Y_0'\nest Y_0$. If $W'= Y_0'$ we set $\bar\Theta=\emptyset$. Otherwise, by the orthogonals for non-split domains property~\eqref{property:orthogonal_for_non_split} there exists a $\nest$-minimal domain $V_1\nest Y_0'$ such that $V_1\orth W'$. Then let $Y'_1=\{V_1\}^\orth_{Y'_0}$, which contains $W'$ and is properly nested into $Y_0'$. Now we can iterate this construction with $W'$ and $Y_1'$, and the procedure has to stop after at most $n$ steps since $Y_i'$ is properly nested inside $Y_{i-1}'$ for all $i$. Thus, in the end we find a simplex $\bar\Theta=\{V_1, \ldots, V_k\}$ of $\nest$-minimal and pairwise orthogonal domains, which are nested in $Y_0'$ and whose orthogonal complement inside $Y_0'$ is $W'$. 
\par\medskip

\textit{Part 3}: Summing up, we have defined some (possibly formally empty) domains and two simplices $\bar\Psi,\bar\Theta\subset \bar X$ such that
$$\begin{cases}
    Y_0=(\bar\Sigma\star\bar\Phi)^\orth;\\
    Y_0'=\bar\Psi^\orth_{Y_0};\\
    W_0=\bar\Sigma^\orth\wedge_\text{min}\bar\Delta^\orth;\\
    W'=\bar\Psi^\orth_{W_0}=\bar\Theta^\orth_{Y'_0}.
\end{cases}$$
Hence we have that
$$(\bar\Sigma\star\bar\Phi\star\bar\Psi\star\bar\Theta)^\orth=W'\nest W_0=\bar\Sigma^\orth\wedge_\text{min}\bar\Delta^\orth.$$
Since $\link(\bar\Sigma)\cap\link(\bar\Delta)$ is the subgraph of $X$ spanned by all $\nest$-minimal domain that are nested inside $W_0$, the previous nesting can be restated as
$$\link(\bar\Sigma)\cap\link(\bar\Delta)\supseteq\link(\bar\Sigma\star\bar\Phi\star\bar\Psi\star\bar\Theta),$$
and since the domains of $\bar\Psi$ lie in $W_0$ by the construction from Part 1 we also have that
$$\link(\bar\Sigma)\cap\link(\bar\Delta)\supseteq\link(\bar\Sigma\star\bar\Phi\star\bar\Psi\star\bar\Theta)\star\bar\Psi.$$

Moreover, the converse inclusion is also true, since by the construction of Part 1 we have that, if a minimal domain $V$ is nested in $W_0$, then either $V$ is one of the vertices of $\bar\Psi$ or $V\nest\bar\Psi^\orth_{W_0}=W'$. Then we have proved that
\begin{equation}\label{eq:intersection_of_lk_in_barX}
\link(\bar\Sigma)\cap\link(\bar\Delta)=\link(\bar\Sigma\star\bar\Phi\star\bar\Psi\star\bar\Theta)\star\bar\Psi,
\end{equation}
where 
$$\begin{cases}
    \bar\Phi=\link(\bar\Sigma)\cap\bar\Delta;\\
    \bar\Psi\subseteq\link(\bar\Delta)\cap\link(\bar\Sigma);\\
    \bar\Theta\subseteq\link(\bar\Sigma)- \text{Star}(\bar\Delta).
\end{cases}$$ 

\textbf{Finding the extension of $\Sigma$}: Let $\bar\Lambda=\bar\Phi\star\bar\Psi\star\bar\Theta$, and let $\Pi$ be the simplex defined as follows:
\begin{itemize}
    \item $p(\Pi)=\bar\Sigma\star\bar\Lambda$;
    \item If $U\in \bar\Sigma$ does not belong to $ \text{Star}(\bar\Delta)$ then $\Pi_U$ is an edge containing $\Sigma_U$, so that $\link_{p^{-1}(U)}\left(\Pi_U\right)=\emptyset$;
    \item If $U\in \bar\Sigma\cap \link(\bar\Delta)$ then $\Pi_U= \Sigma_U$;
    \item If $U\in \bar\Sigma\cap \bar\Delta$ then $\Pi_U=\Sigma_U$ whenever $\Sigma_U$ and $\Delta_U$ are single vertices “of the same kind” (that is, either they are both the vertex of the cone or they are both points in the base); otherwise $\Pi_U$ is an edge containing $\Sigma_v$. In other words, we choose $\Pi_U$ so that 
    $$\link_{p^{-1}(U)}\left(\Pi_U\right)=\link_{p^{-1}(U)}(\Sigma_U)\cap \link_{p^{-1}(U)}(\Delta_U);$$
    \item If $U\in \bar\Phi$ then $\Pi_U= \Delta_U$;
    \item If $U\in \bar\Psi$ then $\Pi_U$ is the cone point $v_U$.
    \item If $U\in \bar\Theta$ then $\Pi_U$ is an edge, so that $\link_{p^{-1}(U)}\left(\Pi_U\right)=\emptyset$;
\end{itemize}

\begin{figure}[htp]
\centering
\renewcommand{\arraystretch}{1.5}
\begin{tabular}{c|c|c}
$\Pi_U$&$\bar\Sigma$&$\link(\bar\Sigma)$\\
\hline
$\bar\Delta$& extend $\Sigma_U$ if needed&$\Delta_U$ for every $U\in\bar\Phi$\\
$\link(\bar\Delta)$&$\Sigma_U$& $v_U$ for every $U\in \bar\Psi$\\
$\bar X- \text{Star}(\bar\Delta)$&complete $\Sigma_U$ to an edge& choose an edge for every $U\in \bar\Theta$
\end{tabular}
    \caption{Schematic representation of the simplex $\Pi$. Each cell describes how $\Pi_U$ is defined whenever the domain $U$ belongs to the area given by the intersection between the row label and the column label (for example, if $U\in\bar\Sigma\cap \link(\bar\Delta)$ we have that $\Pi_U=\Sigma_U$).}
    \label{tab:Sigma_star_Pi}
\end{figure}

Moreover, we define the simplex $\Psi$ such that $p(\Psi)=\bar\Psi$ and that, for every $U\in\bar\Psi$, $\Psi_U=v_U$ (this is exactly how we defined $\Pi_U$ for $U\in\bar\Psi$).

We are finally ready to prove that $\link(\Delta)\cap\link(\Sigma)=\link(\Pi)\star\Psi$. First, we notice that $\Psi\subseteq\link(\Delta)\cap\link(\Sigma)$ since its support is $\bar\Psi$, whose vertices lie in $W_0$. Next, we argue that $\link(\Pi)\subseteq \link(\Delta)$. Let $u\in \link(\Pi)$ and let $U=p(u)$. If $U\in \bar\Sigma\star\bar\Lambda$ then a careful inspection of how we defined $\Pi$ shows that $u\in \link(\Delta)$. Otherwise $U\in \link(\bar\Sigma\star\bar\Lambda)\subset\link(\bar\Delta)$, and therefore $u\in\link(\Delta)$. Thus we showed that $\link(\Delta)\cap\link(\Sigma)\supseteq \link(\Pi)\star\Psi$.

For the converse inclusion, let $u\in \link(\Delta)\cap\link(\Sigma)$, so that $U=p(u)$ belongs to $\text{Star}( \bar\Sigma)\cap \text{Star}(\bar\Delta)$ and $u\in(\link_{p^{-1}(U)}\Sigma_U)\cap(\link_{p^{-1}(U)}\Delta_U)$. There are four possible cases:
\begin{itemize}
    \item If $U\in \bar\Sigma\cap\bar\Delta$ then $(\link_{p^{-1}(U)}\Sigma_U)\cap(\link_{p^{-1}(U)}\Delta_U)=\link_{p^{-1}(U)}(\Pi_U)$, as we already noticed.
    \item If $U\in \bar\Sigma\cap\link(\bar\Delta)$ then $\Delta_U=\emptyset$, and again by construction we have $(\link_{p^{-1}(U)}\Sigma_U)\cap(\link_{p^{-1}(U)}\Delta_U)=\link_{p^{-1}(U)}\Sigma_U=\link_{p^{-1}(U)}(\Pi_U)$.
    \item Symmetrically, if $U\in \link(\bar\Sigma)\cap\bar\Delta=\bar\Phi$ then $\Sigma_U=\emptyset$, and we have $(\link_{p^{-1}(U)}\Sigma_U)\cap(\link_{p^{-1}(U)}\Delta_U)=\link_{p^{-1}(U)}\Delta_U=\link_{p^{-1}(U)}(\Pi_U)$.
    \item Finally, suppose $U\in \link(\bar\Sigma)\cap\link(\bar\Delta)$, that is, $U$ is nested in $W_0$. Then by construction either $U$ is nested in $W'$ or $U\in\bar\Psi$. In the former case $U\in\link(\bar\Sigma\star\bar\Lambda)$, hence $u\in \link(\Pi)$. In the latter case, either $u=v_U$ is the cone point, which belongs to $\Psi$, or $u\in \link_{p^{-1}(U)}(v_U)$, and since $\Pi_U=v_U$ we have that $u\in \link(\Pi)$.
\end{itemize}
This concludes the proof.
\end{proof}

We point out the following by-product of the proof:
\begin{cor}\label{cor:intersection_of_links_in_barX}
    Under Assumption~\ref{ass:standing_assumption}, let $\bar\Sigma,\bar\Delta$ be two simplices of $\bar X$. Then there exist two simplices $\bar\Lambda,\bar\Psi\subset \link(\bar\Sigma)$ such that 
    $$\link(\bar\Sigma)\cap\link(\bar\Delta)=\link(\bar\Sigma\star\bar\Lambda)\star\bar\Psi.$$
    Furthermore, we can assume that $\link(\bar\Sigma\star\bar\Lambda)$ is not the cone with cone point $V\in \link(\bar\Sigma\star\bar\Lambda)$.
\end{cor}

\begin{proof}
    The first part of the statement is just Equation~\eqref{eq:intersection_of_lk_in_barX}. For the “furthermore” part, if $\link(\bar\Sigma\star\bar\Lambda)$ is the cone with cone point $V\in \link(\bar\Sigma\star\bar\Lambda)$, then we have that 
    $$\link(\bar\Sigma\star\bar\Lambda)\star\bar\Psi=\link(\bar\Sigma\star\bar\Lambda\star V)\star\bar\Psi\star V.$$
    In this case we can set $\bar\Lambda'=\bar\Lambda\star V$ and $\bar\Psi'=\bar\Psi\star V$, and then we check again if $\link(\bar\Sigma\star\bar\Lambda')$ is a cone with cone point $V'$. This process must end after finitely many steps, since $\bar X$ has finite dimension by Remark~\ref{rem:finite_dim}.
\end{proof}

As a consequence of the previous Lemma we can also verify another axiom:
\begin{cor}[Verification of Definition~\ref{defn:combinatorial_HHS}.\eqref{item:cHHS_flag}]\label{cor:finite_complexity_HHS_blow_up}
Under Assumption~\ref{ass:standing_assumption}, $X$ has finite complexity in the sense of Definition~\ref{defn:finite_complexity_cHHS}.
\end{cor}

\begin{proof}
One can argue exactly as in the proof of \cite[Claim 6.9]{BHMS}, which only uses that $X$ has finite dimension, as pointed out in Remark~\ref{rem:finite_dim}, and \cite[Condition 6.4.B]{BHMS}, which is our Lemma~\ref{lem:intersection_of_links}. 
\end{proof}

\subsection{Fullness of links}
\begin{lemma}[Verification of Definition~\ref{defn:combinatorial_HHS}.\eqref{item:C_0=C}]\label{lem:edges_in_link}
Under Assumption~\ref{ass:standing_assumption}, there exists a constant $\lambda_0=\lambda_0(E)$ such that $\mathcal{W}$ has the following property whenever $\lambda\geq \lambda_0$. Let $\Delta$ be a non-maximal simplex of $X$.  Suppose that $v,w\in\link(\Delta)$ are distinct, non-adjacent vertices which are contained in $\mathcal{W}$--adjacent maximal simplices $\sigma_v,\sigma_w$. Then there exist maximal simplices $\Pi_v,\Pi_w$ of $X$ such that $\Delta\star v\subseteq\Pi_v$ and $\Delta\star w\subseteq\Pi_w$.
\end{lemma}

\begin{proof}
Recall that $p:\,X \to \bar X$ is the retraction from Definition~\ref{defn:blow-up} that maps every vertex of the blow-up to its support. Moreover, for every maximal simplex $\sigma=\bigstar_{i=1}^k \{v_{U_i},x_i\}$ of $X$, where $U_i\in \frakS $ and $x_i\in \fontact U_i$, let  $(b(\sigma)_W)_{W\in\frakS}$ be the tuple from Definition~\ref{defn:b_sigma}, which was defined as follows:

$$b(\sigma)_W=\begin{cases}
    x_i\mbox{ if }W=U_i;\\
    \bigcup_{U_i\not\bot W} \rho^{U_i}_W\mbox{ otherwise}.
\end{cases}$$

Set $V=p(v)$, $W=p(w)$ and $\bar\Delta=p(\Delta)$. Let $\sigma_v,\sigma_w$ be the two $\mathcal{W}$--adjacent simplices containing $v$ and $w$, respectively, and let $\bar\Sigma=\bar\sigma_v\cap\bar\sigma_w$ (which is possibly empty). Let $k=\text{co-lv}(\bar\Sigma^\orth_{\text{min}})$. Recall that, by Definition~\ref{defn:W-edge} of the edges of $\mathcal{W}$, we have that, for every $U\in\frakS$, 
$$\dist_{\fontact U}(b(\sigma_v)_U,b(\sigma_w)_U)\le (k+1)\lambda.$$

If $V=W$ then we can complete $\Delta$ to two simplices $\Pi_v,\Pi_w$ with the same support (so that, in the sense of Notation~\ref{notation:exceptional_orth}, the co-level is $n$) and such that $v\in\Pi_v, w\in\Pi_w$ and these simplices coincide away from $V$. Since $$\dist_V(b(\Pi_v)_V,b(\Pi_w)_V)=\dist_V(b(\sigma_v)_V,b(\sigma_w)_V)\le (k+1)\lambda\le (n+1)\lambda,$$
while by construction $\dist_U(b(\Pi_v)_U,b(\Pi_w)_U)=0$ whenever $U\neq V$, we are done.

Moreover, notice that $V \not \bot W$, otherwise $v$ and $w$ would be adjacent in $X$ by how we constructed the blow-up graph.
\par\medskip

Thus we are left with the case when $V\neq W$ and $V\not\bot W$. In particular none of them lies inside either $\bar\Delta$ or $\bar\Sigma$. Therefore $V,W\in \link(\bar\Sigma)\cap \link(\bar\Delta)$. By Corollary~\ref{cor:intersection_of_links_in_barX} we can find two simplices $\bar\Phi, \bar\Phi'\subset \link(\bar \Sigma)$ such that 
$$\link(\bar\Sigma)\cap \link(\bar\Delta)=\link(\bar\Delta\star\bar\Phi)\star\bar\Phi'.$$
Moreover, we can assume that $\link(\bar\Delta\star\bar\Phi)$ is not a cone with any of its vertices as cone point, which implies that $(\bar\Delta\star\bar\Phi)^\orth_{\text{min}}=(\bar\Delta\star\bar\Phi)^\orth$ by Lemma~\ref{lem:weak_compl=compl}.

Now, notice that $V$ cannot lie in $\bar\Phi'$, since otherwise it would be orthogonal to every other vertex of $\link(\bar\Sigma)\cap\link(\bar\Delta)$, including $W$. Hence $V$ must lie in $\link(\bar\Delta\star\bar\Phi)$, and the same holds for $W$. Now set $\bar \Psi_v=\bar\sigma_v\cap \link(\bar\Delta\star\bar\Phi)$, which contains $V$ as we just argued, and similarly $\bar \Psi_w=\bar\sigma_w\cap \link(\bar\Delta\star\bar\Phi)$. The situation in $\bar X$ is therefore as  in Figure~\ref{fig:W_edge_1}.

\begin{figure}[htp]
    \centering
    \begin{tikzcd}
&V\in\bar\Psi_v \ar[dl, no head]\ar[d, no head]\ar[dr, no head]&\\
\bar\Phi\ar[r, no head]&\bar\Delta&\bar\Sigma\\
&W\in\bar\Psi_w \ar[ul, no head]\ar[u, no head]\ar[ur, no head]&
\end{tikzcd}
    \caption{The simplices involved in the construction of Lemma~\ref{lem:edges_in_link}, where edges denote joins in $\bar X$. Actually, $\bar\Sigma$ and $\bar\Delta\star\bar\Phi$ need not be disjoint, but none of them can contain $V$ or $W$.}
    \label{fig:W_edge_1}
\end{figure}
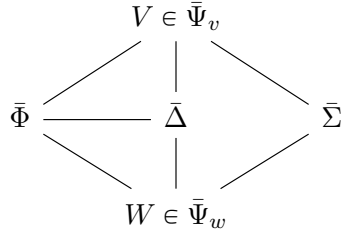

Now, $\link(\bar\Delta\star\bar\Phi)\subseteq \link(\bar\Sigma)$ by construction, or equivalently $(\bar\Delta\star\bar\Phi)^\orth_{\text{min}}\nest \bar\Sigma^\orth_{\text{min}}$ by Lemma~\ref{lem:nested_min_orth}. Hence
\begin{equation}\label{eq:deltaorth_and_sigmaorth}
    (\bar\Delta\star\bar\Phi)^\orth=(\bar\Delta\star\bar\Phi)^\orth_{\text{min}}\nest \bar\Sigma^\orth_{\text{min}}\nest \bar\Sigma^\orth.
\end{equation}

There are two possible cases, depending on whether the two orthogonal complements coincide or not.
\par\medskip

\textbf{Case 1}. Suppose first that $(\bar\Delta\star\bar\Phi)^\orth=\bar\Sigma^\orth$. This means that $\bar\sigma_v=\bar\Sigma\star\bar\Psi_v$, because a support of $\bar\sigma_v$ which is orthogonal to $\bar\Sigma$ must also be orthogonal to $\bar\Delta\star\bar\Phi$. Therefore $\bar\Delta\star\bar\Phi\star\bar\Psi_v$ is already a maximal simplex, and the same is true for $\bar\Delta\star\bar\Phi \star\bar\Psi_w$. Now complete $\Delta$ to maximal simplices $\Pi_v$, $\Pi_w$ as follows:

\begin{itemize}
    \item $\Pi_v$ is supported on $\bar\Delta\star\bar\Phi \star\bar\Psi_v$, and similarly $\Pi_w$ is supported on $\bar\Delta\star\bar\Phi \star\bar\Psi_w$;
    \item if $U\in\bar\Delta$ and $\Delta_U\cap (\fontact U)^{(0)}=\emptyset$ (that is, if $\Delta$ does not prescribe the coordinate for $U$) then choose the same coordinate both for $\Pi_v$ and $\Pi_w$;
    \item if $U\in\bar\Delta$ and $\Delta_U\cap (\fontact U)^{(0)}=\{p_U\}$ (that is, if $\Delta$ already prescribes the coordinate for $U$) then set $\Pi_v=\Pi_w=\{v_U,p_U\}$;
    \item if $U\in\bar\Phi$, choose the same coordinate both for $\Pi_v$ and $\Pi_w$;
    \item if $U\in \bar\Psi_v$ choose for $\Pi_v$ the coordinate prescribed by $\sigma_v$, and similarly if $U\in \bar\Psi_w$ choose for $\Pi_w$ the coordinate prescribed by $\sigma_w$.
\end{itemize}

Now, since $\bar\Sigma^\orth=(\bar\Delta\star\bar\Phi)^\orth$ they have the same co-level $k$. Thus, in order to show that $\Pi_v$ and $\Pi_w$ are $\mathcal{W}$--adjacent, it is enough to prove that, for any $U\in\frakS$, we have $\dist_{\fontact U}(b(\Pi_v)_U, b(\Pi_w)_U)\le (k+1)\lambda$, because $p(\Pi_v)$ and $p(\Pi_w)$ coincide at least on $\bar\Delta\star\bar\Phi$.

If $U\not\bot(\bar\Delta\star\bar\Phi)$ then clearly $\dist_{\fontact U}(b(\Pi_v)_U, b(\Pi_w)_U)=0$. Otherwise $U$ is also orthogonal to $\bar\Sigma$, and by maximality of $\sigma_v$ it cannot be orthogonal to every vertex of $\bar\sigma_v-\bar\Sigma=\bar\Psi_v$. This means that $b(\Pi_v)_U=b(\sigma_v)_U$, since they both depend only on the coordinates over $\bar\Psi_v$, which are the same for both $\Pi_v$ and $\sigma_v$ by construction. For the same reason $b(\Pi_w)_U=b(\sigma_w)_U$. This in turn means that 
$$\dist_{\fontact U}(b(\Pi_v)_U, b(\Pi_w)_U)=\dist_{\fontact U}(b(\sigma_v)_U, b(\sigma_w)_U)\le (k+1)\lambda.$$
\par\medskip

\textbf{Case 2}. Now we are left to deal with the case when $(\bar\Delta\star\bar\Phi)^\orth\sqsubsetneq \bar\Sigma^\orth$. We will find two maximal simplices $\Pi_v,\Pi_w$ whose supports extend $\bar\Delta\star \bar\Phi\star \bar\Psi_v$ and $\bar\Delta\star \bar\Phi\star \bar\Psi_w$, respectively, and then it will suffice to prove that, for every domain $U\in\frakS$, we have
$$\dist_{\fontact U}(b(\Pi_v)_U, b(\Pi_w)_U)\le (k+2)\lambda.$$
In other words, it will be enough to loosen the threshold just by adding a single $\lambda$. Since $\text{co-lv}((\bar\Delta\star \bar\Phi\star \bar\Psi_v)^\orth)\lneq\text{co-lv}(\bar\Sigma^\orth)=k$ we will then have that $\Pi_v$ and $\Pi_w$ are $\mathcal{W}$--adjacent.

Let $\bar\Theta_v=\bar\sigma_v- (\bar\Sigma\star\bar\Psi_v)$ be the simplex spanned by all remaining vertices of $\bar\sigma_v$. Moreover, let $R_v=(\bar\Delta\star\bar \Phi\star \bar \Psi_v)^\orth$, if it exists (if not, the following construction is unnecessary because our simplex is already maximal). Notice that $R_v$ is also orthogonal to $\bar\Sigma$, hence every domain $U\nest R_v$ is orthogonal to every domain in $\bar\Sigma\star\bar\Psi_v$, and therefore it cannot be orthogonal to every vertex of $\bar\Theta_v$ by maximality of $\bar\sigma_v=\bar\Sigma\star\bar\Psi_v\star\bar\Theta_v$.

Then let $r^v=(r^v_U)_{U\nest R_v}\in F_{R_v}$ be the tuple defined as follows:
\begin{itemize}
\item if $U\in \bar\Theta_v$ then $r^v_U$ is the coordinate prescribed by $\sigma_v$;
\item otherwise $r^v_U=\bigcup\rho^{U'}_U$ where the union varies among all $U'\in \bar\Theta_v$ whose projection to $U$ is defined.
\end{itemize}
By the previous argument, $r^v_U$ is well-defined for any $U\nest R_v$. Moreover, arguing exactly as in Remark~\ref{rem:b_sigma_consistent}, one sees that $r^v$ is indeed a $20E$-consistent tuple, and therefore an element of $F_{R_v}$. 

Now, if $R_v$ is $\nest$-minimal then $r^v$ is just a point in $\fontact R_v$, and we set $\bar\Omega_v=R_v$ and $\Omega_v$ as the edge $\{v_{R_v},r^v\}$. Otherwise, by the EDPR property~\eqref{property:EDPR} there exist a maximal family of pairwise orthogonal, $\nest$-minimal domains $\bar\Omega_v=\{O_1,\ldots, O_l\}$ whose realisation point is $C_0$-close to the realisation of $r$ in $F_{R_v}$. This means that, if we define the simplex $\Omega_v$, supported in $\bar\Omega_v$, by choosing for every $I\in\bar\Omega_v$ the coordinate $r^v_I$, then for every $U\nest R_v$ and every maximal simplex $\Omega_v'$ containing $\Omega_v$, the $U$-coordinate of the realisation tuple $b(\Omega_v')$ is $M$-close to $r^v_U$, where $M=M(C_0,E)$ is a constant coming from the distance formula, Theorem~\ref{thm:distance_formula}.

Define $\bar\Theta_w, R_w, \bar\Omega_w$ analogously, so that the situation looks like in Figure~\ref{fig:W_edge_2}.

\begin{figure}[htp]
\centering
\begin{tikzcd}
\bar\Omega_v \ar[d, no head] \ar[dr, no head] \ar[drr, bend left=120, looseness=1.3, no head]\ar[r, no head]&V\in\bar\Psi_v \ar[dl, no head]\ar[d, no head]\ar[dr, no head]&\bar\Theta_v \ar[d, no head] \ar[l, no head]\\
\bar\Phi\ar[r, no head]&\bar\Delta&\bar\Sigma\\
\bar\Omega_w \ar[u, no head] \ar[ur, no head] \ar[urr, bend right=120, looseness=1.3, no head]\ar[r, no head]&W\in\bar\Psi_w \ar[ul, no head]\ar[u, no head]\ar[ur, no head]&\bar\Theta_w \ar[u, no head] \ar[l, no head]
\end{tikzcd}
\caption{The simplices of $\bar X$ involved in the second case of the proof of Lemma~\ref{lem:edges_in_link}, where edges represent joins. By construction $\bar\sigma_v=\bar\Sigma\star\bar\Psi_v\star\bar\Theta_v$, and similarly for $\bar\sigma_w=\bar\Sigma\star\bar\Psi_w\star\bar\Theta_w$.}
\label{fig:W_edge_2}
\end{figure}
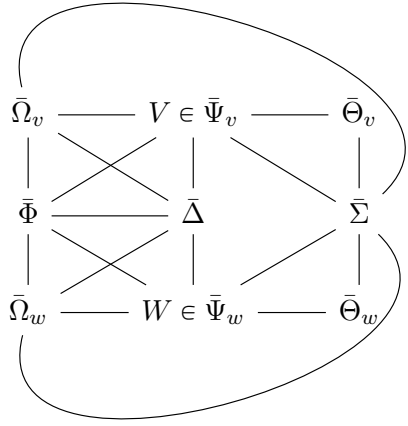

Then complete $\Delta$ to maximal simplices $\Pi_v$, $\Pi_w$ as follows:

\begin{itemize}
    \item $\Pi_v$ is supported on $\bar\Delta\star\bar\Phi \star\bar\Psi_v\star\bar\Omega_v$, and similarly $\Pi_w$ is supported on $\bar\Delta\star\bar\Phi \star\bar\Psi_w\star\bar\Omega_w$;
    \item if $U\in\bar\Delta\star\bar\Phi$, choose the same coordinates both for $\Pi_v$ and $\Pi_w$;
    \item if $U\in \bar\Psi_v$ choose for $\Pi_v$ the coordinate coming from $\sigma_v$, and similarly if $U\in \bar\Psi_w$ choose for $\Pi_w$ the coordinate coming from $\sigma_w$
    \item if $U\in \bar\Omega_v$ choose the coordinate coming from $\Omega_v$, and similarly if $U\in \bar\Omega_v$ choose the coordinate coming from $\Omega_v$.
\end{itemize}

Now we show that for every $U\in\frakS$ we have $\dist_{\fontact U}(b(\Pi_v)_U, b(\Pi_w)_U)\le (k+2)\lambda$. 

If $U\not\bot(\bar\Delta\star\bar\Phi)$ then clearly $\dist_{\fontact U}(b(\Pi_v)_U, b(\Pi_w)_U)=0$. Otherwise $U$ is also orthogonal to $\bar\Sigma$, since $(\bar\Delta\star\bar\Phi)^\orth\nest \bar\Sigma^\orth$.

\begin{claim}\label{claim:m20e} If $U\orth(\bar\Delta\star\bar\Phi)$ then
    $\dist_{\fontact U}(b(\Pi_v)_U,b(\sigma_v)_U)\le M$.
\end{claim}

\begin{claimproof}[Proof of Claim~\ref{claim:m20e}]
First suppose that $U\not\bot \bar\Psi_v$. Then $b(\Pi_v)_U$ and $b(\sigma_v)_U$ both contain the set
$$\alpha^v_U=\begin{cases}
    (\Psi_v)_U\mbox{ if }U\in \bar\Psi_v;\\
    \bigcup_{T\in\bar\Psi_v,\, T\not\bot U}\rho^T_U\mbox{ otherwise}.
\end{cases}$$
We are left with the case when $U\orth \bar\Delta\star\bar\Phi\star\bar \Phi_v$. In this case $b(\Pi_v)_U$ contains the set
$$\beta^v_U=\begin{cases}
    (\Psi_v)_U\mbox{ if }U\in \bar\Omega_v;\\
    \bigcup_{T\in\bar\Omega_v,\, T\not\bot U}\rho^T_U\mbox{ otherwise}.
\end{cases}$$
But by our definition of $\Omega_v$, $\beta_v$ is $M$-close to the coordinate 
$$r^v_U=\begin{cases}
    (\Psi_v)_U\mbox{ if }U\in \bar\Theta_v;\\
    \bigcup_{T\in\bar\Theta_v,\, T\not\bot U}\rho^T_U\mbox{ otherwise},
\end{cases}$$
which is a subset of $b(\sigma_v)_U$.
\end{claimproof}

The proof of the Claim also applies to $\Pi_w$ and $\sigma_w$. Hence
$$\dist_{\fontact U}(b(\Pi_v)_U,b(\Pi_w)_U)\le \dist_{\fontact U}(b(\sigma_v)_U,b(\sigma_w)_U)+ 2M\le (k+1)\lambda+2M.$$
Therefore it is enough to choose $\lambda_0= 2M$, so that 
$(k+1)\lambda+2M\le (k+2)\lambda$ whenever $\lambda\ge \lambda_0$.
\end{proof}

\subsection{Hyperbolic links}\label{subsec:hyperbolic_links}
Recall that, for any domain $U$, we defined $F_U$ as the space of $20E$-consistent tuples for $U$, which is a sub-HHS of $(\cuco Z,\frakS)$ with maximal domain $U$. For convenience, up to quasi-isometry we can assume that $F_U$ is a graph (again, by \cite[Lemma 3.B.6]{cornulier2014metric}). Moreover, for every $V\propnest U$ one has the relative product region $P^U_V\subset F_U$, which we will often refer to as $P_V$ when the ambient domain will be clear.

\begin{defn}\label{defn:factored_space}
    The  \emph{factored space} $\hat F_U$ associated to $U$ is the graph obtained from $F_U$ by coning off the product region $P_V$ for every $V\propnest U$.
\end{defn}
We will denote by $H_V$ the vertex of the cone over $P_V$. By \cite[Corollary 2.9 and Remark 2.10]{Asymptotic_dim}, if $\cuco Z$ is a normalised HHS then $\hat F_U$ is uniformly quasi-isometric to $\fontact U$, and therefore uniformly hyperbolic. More precisely, the quasi-isometry is induced by the projection $\pi_U:\,F_U\to \fontact U$. Then the strategy to prove that links of simplices inside $X$ are hyperbolic will be to show that each of them is either bounded, quasi-isometric to a coordinate space or to a factored space.

\begin{lemma}\label{lem:hyperbolic_links}
Under Assumption~\ref{ass:standing_assumption}, there exists $\lambda_1\ge \lambda_0$ such that the following holds whenever $\lambda\ge \lambda_1$. There exists $\delta$, depending both on $\lambda$ and on the HHS structure, such that, for every non-maximal simplex $\Delta\subset X$, the associated coordinate space $\fontact ([\Delta])$ is $\delta$-hyperbolic.
\end{lemma}

\begin{proof} 
We consider all possible shapes of $\link(\Delta)$, according to Corollary~\ref{cor:bounded_links}.
\par\medskip

\textbf{If $\link(\Delta)$ is a point or a non-trivial join}: In this case $\fontact ([\Delta])$, which is obtained by adding edges to $\link(\Delta)$, has diameter at most $2$, and therefore it is $2$-hyperbolic.
\par\medskip

\textbf{If $\Delta$ is almost-maximal and $\Delta_U=v_U$ for some $U\in\bar\Delta$}: In this case $\link(\Delta)$ is the base of the cone over $\fontact U$. Now, by construction two points $p,q\in\link(\Delta)$ belong to $\mathcal{W}$--adjacent maximal simplices if and only if $$\dist_{\fontact U}(p,q)\le (n+1)\lambda.$$
This shows that $\fontact ([\Delta])$ and $\fontact U$ are quasi-isometric with uniform constants, and since $\fontact U$ is $E$-hyperbolic then $\fontact ([\Delta])$ is $\delta$-hyperbolic for some constant $\delta=\delta(\lambda,E,n)$.
\par\medskip

\textbf{If $\Delta_W$ is an edge for every $W\in\bar\Delta$}: 
In this case $\fontact ([\Delta])$ is $(2,2)$-quasi-isometric to the subgraph of $\bar X^{+\mathcal{W}}$ spanned by $\link(\bar\Delta)$, via the retraction $p:\,X\to \bar X$. Let $\link(\bar\Delta)^{+\mathcal{W}}$ be this subgraph, and let $U=\bar\Delta^\orth$. 

Now, if $U$ is $\nest$-minimal then $\link(\bar\Delta)^{+\mathcal{W}}$ consists only of $U$, hence is uniformly bounded. Then suppose that $U$ is non-$\nest$-minimal. At the level of vertices we can define a map $\psi:\,\link(\bar\Delta)^{+\mathcal{W}}\to \hat F_U$ by sending every $\nest$-minimal domain $V\propnest U$ to the cone point $H_V$ over the corresponding product region $P_V$. Our goal is to prove that $\psi$ is a quasi-isometry with uniform constants, and therefore $\link(\bar\Delta)^{+\mathcal{W}}$ is uniformly hyperbolic.

First we show that $\psi$ is coarsely surjective. For every $x\in F_U$, by the EDPR property~\eqref{property:EDPR} we can find a maximal family  $V_1,\ldots, V_k\propnest U$ of $\nest$-minimal domains whose product region is $C_0$-close to $x$. Then in particular the product region $P_{V_1}$ is $C_0$-close to $x$ in $F_U$, which means that the corresponding cone point $H_{V_1}=\psi(V_1)$ is $(C_0+1)$-close to $x$ in $\hat F_U$. 

Next, we prove that the map $\psi$ is uniformly Lipschitz, by showing that adjacent vertices in $\link(\bar\Delta)^{+\mathcal{W}}$ map to uniformly close points inside $\hat F_U$. If $V,V'\propnest U$ are $\nest$-minimal domains which are joined by an edge of $\mathcal{W}$ then one of the following must hold:
\begin{itemize}
    \item $V\orth V'$: in this case, the product regions $P_V$ and $P_{V'}$ are already uniformly close inside $F_U$. Indeed, if one chooses two coordinates $p\in \fontact V$ and $p'\in \fontact V'$ then, by the partial realisation axiom~\eqref{item:dfs_partial_realisation} for the sub-HHS $F_U$, one can find an element $x\in F_U$ whose coordinates must satisfy the following properties:
    $$\begin{cases}
        \dist_V(x_V, p)\le E;\\
        \dist_{V'}(x_{V'}, p')\le E;\\
        \dist_{U}(x_{U}, \rho^V_U)\le E \mbox{ whenever }U\not \bot V;\\
        \dist_{U}(x_{U}, \rho^{V'}_U)\le E \mbox{ whenever }U\not \bot V'.\\
    \end{cases}$$
    Therefore, by definition, $x$ belongs to both $P_V$ and $P_V'$. This means that the cone points $H_V$ and $H_{V'}$ are at distance at most $2$.
    \item There exist two maximal simplices $\Pi$, $\Pi'$ which extend $\Delta$ and such that $V\in\bar\Pi$ and $V'\in\bar\Pi'$, and for every $W\in\frakS$ we have 
    $$\dist_{\fontact W}(b(\Pi)_W, b(\Pi')_W)\le (1+\text{co-lv}(U))\lambda\le (1+n)\lambda.$$
    In this case let $x=(b(\Pi)_W)_{W\nest U}$ and $x'=(b(\Pi')_W)_{W\nest U}$ be the corresponding tuples inside $F_U$. By the Distance Formula (Theorem~\ref{thm:distance_formula}) for the sub-HHS $F_U$, the distance of these points is bounded by some constant $D$ depending only on $\lambda,n,E$. Moreover $x\in P_V$ and $x'\in P_{V'}$ by construction, and therefore $x$ is adjacent to $H_V$ inside $\hat F_U$ and similarly for $x'$ and $H_{V'}$. Therefore $\dist_{\hat F_U}(H_V,H_{V'})\le 2+D$.
\end{itemize}

In order to complete the proof that $\psi$ is a quasi-isometry we need the following, which is the only spot where we have to choose $\lambda_1$ carefully:

\begin{claim}\label{claim:minimal_product_regions}
    Under Assumption~\ref{ass:standing_assumption}, and with the notation of Lemma~\ref{lem:hyperbolic_links}, there exists $\lambda_1\ge \lambda_0$ such that the following holds whenever $\lambda\ge \lambda_1$. Let $V,W\propnest U$ be two $\nest$-minimal domains. If the product regions $P_V$ and $P_W$ lie within distance at most $2C_0+1$ in $F_U$, where $C_0$ is the constant from the EDPR property~\eqref{property:EDPR}, then $V$ and $W$ are $\mathcal W$--adjacent in $\link(\bar\Delta)^{+\mathcal{W}}$.
\end{claim}

\begin{claimproof}[Proof of Claim~\ref{claim:minimal_product_regions}]
Let $y\in P_V, z\in P_W$ be two points that are $(2C_0+1)$-close in $F_U$. Moreover, let $\bar\Theta=\{T_1,\ldots, T_l\}$ be a family of pairwise orthogonal, $\nest$-minimal domains such that $T_i\orth U$ for all $i=1,\ldots, l$. Our goal is to complete $\bar\Theta$ to a maximal family $V=V_0,V_1,\ldots, V_k$ of pairwise orthogonal, $\nest$-minimal elements whose product region is uniformly close to $y$, and similarly for $W$ and $z$.

If $V=V_0$ has no orthogonal inside $U$ then we have nothing to do. Otherwise, define $a=(y_I)_{I\nest V^\orth}\in F_{V^\orth_U}$ as the coordinates of $y$ in the domains of the orthogonal complement. If $V^\orth_U$ is not itself minimal, by the EDPR property~\eqref{property:EDPR} there exist a maximal family $V_1,\ldots, V_k\propnest V^\orth_U$ of $\nest$-minimal and pairwise orthogonal domains whose product region inside $F_{V^\orth_U}$ is $C_0$-close to $a$. Either way, there exist $p_i\in \fontact V_i$, for $i=0,\ldots,k$, such that the realisation point $b(\Sigma)$ of the simplex $\Sigma=\bigstar_{i=0}^k\{v_{V_i},p_i\}\subset X$ is $C_0$-close to $y$. Arguing the same way for $W$ we can find a simplex $\Sigma'$, whose support contains $W$ and whose realisation point is $C_0$-close to $w$. Therefore the realisation points of these two simplices are $(4C_0+1)$-close, and since projections to coordinate spaces are uniformly coarsely Lipschitz there exists $M=M(C_0)$ such that, for every domain $I\in \frakS$, $\dist_{\fontact I}(b(\Sigma), b(\Sigma'))\le M$. Therefore, if we set $\lambda_1=M$ we have that, whenever $\lambda\ge \lambda_1$, for every $I\in \frakS$
$$\dist_{\fontact I}(b(\Sigma), b(\Sigma'))\le \lambda.$$
Hence $\Sigma$ and $\Sigma'$ are $\mathcal{W}$--adjacent by definition, and therefore $V$ and $W$ are $\mathcal{W}$--adjacent in $\link(\bar\Delta)^{+\mathcal{W}}$.
\end{claimproof}

Now, we claim that $\dist_{\link(\bar\Delta)^{+\mathcal{W}}}(V,V')\le 4\dist_{\hat F_U}(H_V,H_{V'})+2$ for every two $\nest$-minimal domains $V,V'\propnest U$, and this will complete the proof that $\psi$ is a quasi-isometric embedding. Consider a geodesic $\gamma\subset \hat F_U$ from $H_{V}$ to $H_{V'}$. The vertices of $\gamma$ can either be tuples of $F_U$ or cone points associated to product regions. Let $x_1,\ldots, x_k$ be the vertices of $\gamma$ belonging to $F_U$, and for every $x_i$ consider a $\nest$-minimal domain $I_i$ whose product region is $C_0$-close to $x_i$ in $F_U$, which exists by the EDPR property~\eqref{property:EDPR}. Then the situation in $\hat F_U$ is as shown in Figure~\ref{fig:hatF}.

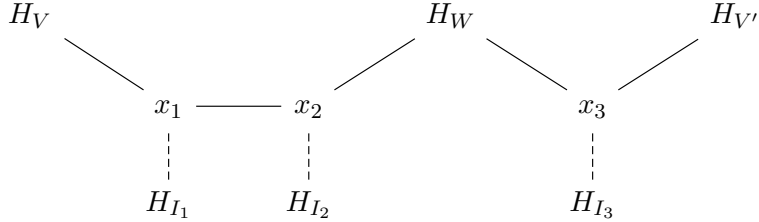
\begin{figure}[htp]
    \centering
    \begin{tikzcd}
        H_V\ar[dr, no head]&&&H_W\ar[dr, no head]&&H_{V'}\\
        &x_1\ar[r, no head]&x_2\ar[ur, no head]&&x_3\ar[ur, no head]&\\
       & H_{I_1}\ar[u, no head, dashed]&H_{I_2}\ar[u, no head, dashed]&&H_{I_3}\ar[u, no head, dashed]&
    \end{tikzcd}
    \caption{The continuous segments represent the path $\gamma\subset \hat F_U$. For each vertex $x_i$ of $\gamma$ which lies inside $F_U$ we choose a $\nest$-minimal domain $I_i\in \frakS$ whose product region is $C_0$-close to $x_i$ in $F_U$.}
    \label{fig:hatF}
\end{figure}

Now, Claim~\ref{claim:minimal_product_regions} implies that $V$ is $\mathcal{W}$--adjacent to $I_1$ (since $x_1\in P_V$ by construction), and similarly that $V'$ is $\mathcal{W}$--adjacent to $I_k$. Then the proof is complete if we show that, for every $i=1,\ldots, k-1$, $\dist_{\link(\bar\Delta)^{+\mathcal{W}}}(I_i,I_{i+1})\le 4$, because then $\dist_{\link(\bar\Delta)^{+\mathcal{W}}}(V, V')\le 4L(\gamma)+2$ where $L(\gamma)$ is the length of $\gamma$.

Thus let $x_i,x_{i+1}$ be two consecutive vertices of $\gamma\cap F_U$. If they are joined by an edge of $F_U$ (in Figure~\ref{fig:hatF} this happens to $x_1$ and $x_2$) we have that $P_{I_i}$ and $P_{I_{i+1}}$ are at most $(2C_0+1)$-close, thus $I_i$ and $I_{i+1}$ are $\mathcal{W}$--adjacent by Claim~\ref{claim:minimal_product_regions}.

Otherwise $\gamma$ might contain a segment of the form $\{x_i, H_W, x_{i+1}\}$, where $W\propnest U$ is a domain such that $x_i,x_{i+1}\in P_W$ (in Figure~\ref{fig:hatF} this happens to $x_2$ and $x_3$).

If $W$ is $\nest$-minimal then by Claim~\ref{claim:minimal_product_regions} we see that $I_i$ and $W$ are $\mathcal{W}$--adjacent, and similarly for $I_{i+1}$ and $W$. Therefore $\dist_{\link(\bar\Delta)^{+\mathcal{W}}}(I_i, I_{i+1})\le 2$.

Otherwise suppose that $W$ is not $\nest$-minimal, and let $y_i=((x_i)_T)_{T\nest W}\in F_W$, which is again $20E$-consistent since it is a sub-tuple of a consistent tuple. By the EDPR property~\eqref{property:EDPR} there exists a minimal domain $R_i\propnest W$ such that $y_i$ is $C_0$-close to the relative product region $P_{R_i}^W$ inside $F_W$. This also means that the whole tuple $x_i$ is $C_0$-close to the relative product region $P_{R_i}= P_{R_i}^U$ inside $F_U$, since any tuple in $P_{R_i}^W$ can be completed to a tuple in $P_{R_i}^U$ by choosing some coordinates for the domains $T\nest W^\orth_U$. Therefore $I_i$ and $R_i$ are $\mathcal{W}$--adjacent, by Claim~\ref{claim:minimal_product_regions}, since their product regions in $F_U$ are $(2C_0+1)$--close. 

Now, if $W$ is split then pick a Samaritan $Q\nest W$. Hence either $R_i=Q$, or $R_i\orth Q$ and therefore they are $\mathcal{W}$--adjacent. Then $\dist_{\link(\bar\Delta)^{+\mathcal{W}}}(I_i, Q)\le 2$, and by repeating this procedure with $I_{i+1}$ and the same Samaritan $Q$ we also get that $\dist_{\link(\bar\Delta)^{+\mathcal{W}}}(I_{i+1}, Q)\le 2$. Therefore $\dist_{\link(\bar\Delta)^{+\mathcal{W}}}(I_i, I_{i+1})\le 4$.

If otherwise $W$ is non-split, by property~\eqref{property:orthogonal_for_non_split} there exists a minimal domain $Q'\propnest U$ which is orthogonal to $W$. Then $Q'$ and $R_i$ are $\mathcal{W}$--adjacent since they are orthogonal, and arguing as before we get that $$\dist_{\link(\bar\Delta)^{+\mathcal{W}}}(I_i, I_{i+1})\le \dist_{\link(\bar\Delta)^{+\mathcal{W}}}(I_i, Q')+\dist_{\link(\bar\Delta)^{+\mathcal{W}}}(Q', I_{i+1}) \le 4.$$
This concludes the proof.\end{proof}

\subsection{QI-embedding of coordinate spaces}\label{subsec:qi_embedded_links}
\begin{lemma}
    Under Assumption~\ref{ass:standing_assumption}, whenever $\lambda\ge \lambda_1$ there exists a constant $\delta'$, depending both on $\lambda$ and on the HHS structure, such that, for every non-maximal simplex $\Delta\subset X$, the associated coordinate space $\fontact ([\Delta])$ is $(\delta',\delta')$-quasi-isometrically embedded inside $Y_\Delta$. 
\end{lemma}

\begin{proof}
    Again, we look at all possible shapes of $\link(\Delta)$. If it is a single point or a non-trivial join then of course it is $(2,2)$-quasi-isometrically embedded inside $Y_\Delta$. In all other cases, since $\fontact ([\Delta])$ is a subgraph of $Y_\Delta$, it will be enough to construct a coarsely Lipschitz retraction from $Y_\Delta$ to $\link(\Delta)$.
    \par\medskip
    
    \textbf{If $\Delta$ is almost maximal}: Let $U$ be the $\nest$-minimal domain such that $\Delta_U=v_U$. Then $\link(\Delta)^{(0)}$ is the copy of $(\fontact U)^{(0)}$ inside $X$, and since in Lemma~\ref{lem:hyperbolic_links} we showed that $\fontact ([\Delta])$ is uniformly quasi-isometric to $\fontact (U)$, we will often replace distances in $\fontact ([\Delta])$ with distances in $\fontact (U)$ without explicit mention. 

    \begin{claim}\label{claim:Y_Delta_for_minimal} $Y_\Delta$ is the subgraph spanned by $\link(\Delta)^{(0)}=(\fontact U)^{(0)}$ and the cones over $\fontact V$ for all $\nest$-minimal domains $V\transverse U$. In other words
    $$Y_\Delta=\text{span}_{X^{+\mathcal{W}}}\left\{\fontact U^{(0)}\cup\bigcup_{p(v)\transverse U}\{v\}\right\}.$$
    \end{claim}

\begin{claimproof}[Proof of Claim~\ref{claim:Y_Delta_for_minimal}]
    Since, by definition $Y_\Delta=\text{span}_{X^{+\mathcal{W}}}\left\{X-\Sat(\Delta)\right\}^{(0)}$, we will equivalently prove that 
    $$\Sat(\Delta)^{(0)}=\{v_U\}\cup\bigcup_{p(v)\orth U} \{v\}.$$
    Indeed, clearly $v_U\in \Sat(\Delta)^{(0)}$, while $\fontact U^{(0)}\cap\Sat(\Delta) =\emptyset$. Moreover, whenever $V=p(v)$ is orthogonal to $U$, it is possible to complete $\{U,V\}$ to a maximal family of pairwise orthogonal domains $\bar\Sigma=\{U,V,V_1,\ldots,V_k\}$, and we can find an almost-maximal simplex $\Sigma$ supported in $\bar\Sigma$ which contains $v$ and such that $\link(\Sigma)=\link(\Delta)=\fontact U^{(0)}$. Hence $v\in \Sat(\Delta)$. Conversely, if $v\in \Sat(\Delta)$ then $\link(v)\supseteq \link(\Delta)$ by definition, therefore $V=p(v)$ either coincides with, or is orthogonal to $U$.
\end{claimproof}
    
    Now, at the level of vertices, we can define a retraction $r:\,Y_\Delta\to \fontact ([\Delta])$ as follows:
    $$r(v)=\begin{cases}
        v \mbox{ if }v\in (\fontact ([\Delta]))^{(0)};\\
        \rho^{p(v)}_U \mbox{ otherwise.}
    \end{cases}$$
    Notice that the retraction is well-defined, as a consequence of Claim~\ref{claim:Y_Delta_for_minimal}. We are left to prove that this retraction is coarsely Lipschitz. Let $v,v'$ be two adjacent vertices in $Y_\Delta$, and we will show that $\dist_{\fontact (\Delta)}(r(v),r(v'))$ is uniformly bounded from above.
    \begin{itemize}
        \item If both $v,v'$ belong to $(\fontact U)^{(0)}$ then they are adjacent in $\fontact ([\Delta])$, by how the latter is defined.
        \item If $v\in (\fontact U)^{(0)}$ but $v'\not\in (\fontact U)^{(0)}$ then, setting $V'=p(v')$ we must have that $U\transverse V'$. Moreover, by definition of $\mathcal{W}$--edges, there must be two simplices $\sigma, \sigma'$ such that $v\in \sigma$, $v'\in \sigma'$ and the corresponding realisation tuples $b(\sigma),b(\sigma')$ are at least $(n+1)\lambda$-close in every coordinate space. Now $b(\sigma)_U=v$ by construction, while $b(\sigma')_U$ is a set of diameter $10E$ which contains that $\rho^{V'}_U$. Hence 
        $$\dist_U(r(v),r(v'))=\dist_U(v,\rho^{V'}_U)=\dist_U(b(\sigma)_U,b(\sigma')_U)+10E\le (n+1)\lambda+10E.$$
        \item We are left with the case where both $v,v'$ do not belong to $(\fontact U)^{(0)}$, and we want to find an upper bound for $\dist_U(r(v),r(v'))=\dist_U(\rho^{V}_U,\rho^{V'}_U)$, where $V=p(v)$ and $V'=p(v')$. If $V=V'$ then we have nothing to prove. Otherwise $V$ and $V'$ may be adjacent in $\bar X^{+\mathcal{W}}$ for two different reasons. If $V\perp V'$, then by Lemma~\ref{lem:close_proj_of_orthogonals} we see that $\dist_U(\rho^{V}_U,\rho^{V'}_U)\le 2E$. Otherwise there are two simplices $\sigma, \sigma'$ such that $v\in \sigma$, $v'\in \sigma'$ and the corresponding realisation tuples $b(\sigma),b(\sigma')$ are at least $(n+1)\lambda$-close in every coordinate space. Then, similarly to the previous case, we have that
        $$\dist_U(\rho^{V}_U,\rho^{V'}_U)\le\dist_U(b(\sigma),b(\sigma'))+20E\le (n+1)\lambda+20E.$$
    \end{itemize}

    \textbf{If $\Delta_V$ is an edge for every $V\in \bar\Delta$}: Let $\bar\Delta=\{V_1,\ldots, V_l\}$ and let $U$ be the orthogonal complement of $\bar\Delta$. If $U$ is split then $\link(\Delta)$ is a join, and therefore it is quasi-isometrically embedded into $Y_\Delta$ since it is uniformly bounded. 
    
    Thus we can assume that $U$ is non-split. The next step is the following:
    
    \begin{claim}\label{claim:ydelta_not_orth}
        If a vertex $v\in X$ belongs to $Y_\Delta$ then $V=p(v)$ is not orthogonal to $U$, and therefore $\rho^{p(v)}_U$ is well-defined. 
    \end{claim}
    \begin{claimproof}[Proof of Claim~\ref{claim:ydelta_not_orth}]
        We prove the contrapositive of the statement, that is, we show that if $V=p(v)$ is orthogonal to $U$ then 
        $v\in\Sat(\Delta)$. We have that $U\nest V^\orth$, and either they coincide or there exists a $\nest$-minimal domain $V_1$ inside $V^\orth$ such that $V_1\orth U$, by property~\eqref{property:orthogonal_for_non_split} (which applies since we are in the case when $U$ is non-split). Then after finitely many steps we can find a simplex $\bar\Sigma=\{V=V_0,\ldots, V_k\}$ containing $V$ and whose orthogonal complement is $U$. Then $U=\bar\Sigma^\orth=\bar\Delta^\orth$, and therefore $\link_{X^{+\mathcal{W}}}(\bar\Sigma)=\link_{X^{+\mathcal{W}}}(\bar\Delta)$ because they are both spanned by the $\nest$-minimal domains inside $U$. Hence, since $p(v)=V\in \bar\Sigma$ we have that $v\in\Sat(\Delta)$. 
    \end{claimproof}

    Now, the proof of Lemma~\ref{lem:hyperbolic_links} gives a uniform quasi-isometry $\fontact ([\Delta])\mapsto\hat F_U$, mapping the cone over a $\nest$-minimal domain $V\propnest U$ to the cone point $H_V$. In turn, $F_U$ is uniformly quasi-isometric to $\fontact U$ via the projection map. Hence the composition of these maps is a quasi-isometry $\psi:\,\fontact ([\Delta])\to\fontact U$, which maps the cone over a $\nest$-minimal domain $V\propnest U$ to $\rho^V_U$.
    
    Now define $r:\,Y_\Delta\to \fontact ([\Delta])$ by mapping every vertex $v\in Y_\Delta$ to $\rho^{p(v)}_U\in \fontact U$ (which is well-defined by Claim~\ref{claim:ydelta_not_orth}) and then applying the inverse quasi-isometry $\psi^{-1}:\,\fontact U \to\fontact ([\Delta])$. Notice that if $p(v)\nest U$ then $r(v)=p(v)$, and therefore $r$ is a coarse retraction onto $\fontact ([\Delta])$. 
    
    We are left to prove that, whenever $v,v'\in Y_\Delta$ are adjacent vertices, then $\dist_U(\rho^{p(v)}_U,\rho^{p(v')}_U)$ is uniformly bounded from above, and therefore $r$ is coarsely Lipschitz as it is the composition of a Lipschitz map and the uniform quasi-isometry $\psi^{-1}$. There are three possibilities:
    \begin{itemize}
        \item If $p(v)=p(v')$ then we have nothing to prove.
        \item If $p(v)\orth p(v')$ then by Lemma~\ref{lem:close_proj_of_orthogonals} we have $\dist_U(\rho^{p(v)}_U,\rho^{p(v')}_U)\le 2E$;
        \item If $v,v'$ lie in $\mathcal{W}$-adjacent simplices $\sigma, \sigma'$, respectively, then we have that
        $$\dist_U(\rho^{p(v)}_U,\rho^{p(v')}_U)\le \dist_U(b(\sigma),b(\sigma'))+20E\le (n+1)T+20E$$
    \end{itemize}
    This concludes the proof.
\end{proof}

\subsection{The realisation map is a quasi-isometry}
We are left to prove the following Lemma, which concludes the proof of Theorem~\ref{thm:XWf}:
\begin{lemma}\label{lem:W_Z_qi}
Under Assumption~\ref{ass:standing_assumption} there exists $\lambda_2\geq \lambda_1$ such that, whenever  $\lambda\geq \lambda_2$, the map $f:\mathcal{W}\to \cuco Z$ from Definition~\ref{defn:realisation_map} is a 
quasi-isometry.
\end{lemma}

\begin{proof}
First we show that $f$ is Lipschitz. Given two $\mathcal{W}$--adjacent simplices $\Sigma,\Delta$, we have that, for every $U\in\frakS$, $\dist_{\fontact U}(b(\Sigma)_U,b(\Delta)_U)\le (n+1)\lambda$, because $n$ is the maximum co-level, so gives the highest threshold in the Definition~\ref{defn:W-edge} of $\mathcal{W}$--edges. Moreover, since $f(\Sigma)$ realises $b(\Sigma)$ and $f(\Delta)$ realises $b(\Delta)$, for every $U\in \frakS$ we have that
$$\dist_{\fontact U}(f(\Sigma),f(\Delta))\le2E+\dist_{\fontact U}(b(\Sigma)_U,b(\Delta)_U)\le 2E+(n+1)\lambda.$$
Thus, by the distance formula (Theorem~\ref{thm:distance_formula}), there exists $M=M(E, \lambda,n)$ such that $\dist_{\cuco Z}(f(\Sigma), f(\Delta))\le M$. This proves that $f$ is $M$-Lipschitz.

Furthermore, $f$ is coarsely surjective. Indeed, the whole HHS $\cuco Z$ coincides with $F_S$, where $S\in \frakS$ is the $\nest$-maximal element. Hence, by the EDPR property~\ref{property:EDPR} there exists a constant $C_0$ such that every $z\in \cuco Z$ is $C_0$-close to the product region $P_{\{V_i\}_{i=1,\ldots, k}}$ associated to a maximal family $V_1,\ldots, V_k$ of $\nest$-minimal, pairwise orthogonal domains. In particular, $z$ is $C_0$-close to some realisation point $z'$ for some simplex $\Delta$, and such a point uniformly coarsely coincides with $f(\Delta)$ by the uniqueness axiom~\eqref{item:dfs_uniqueness}.

We are left to prove that, for every two maximal simplices $\Sigma,\Delta$ of $X$, their distance in $\mathcal{W}$ is bounded above in terms of $\dist_{\cuco Z}(f(\Sigma), f(\Delta))$. Now, ${\cuco Z}$ is a $K$-quasigeodesic metric space for some $K\ge0$, therefore it is possible to find a $(K,K)$-quasigeodesic path 
$$\gamma=\{x_0=f(\Sigma),x_1,\ldots,x_{l-1}, x_l=f(\Delta)\}$$
from $f(\Sigma)$ to $f(\Delta)$. In particular, the number $l$ of vertices of this path is bounded above by $K \dist_{\cuco Z}(f(\Sigma), f(\Delta))+K$, and for every $i=0,\ldots, l-1$ we have that $\dist_{\cuco Z}(x_i, x_{i+1})\le 2K$.

Moreover, by coarse surjectivity of $f$, for every $i=1,\ldots, l-1$ we can find a simplex $\Sigma_i$ such that $f(\Sigma_i)$ is $C_0$-close to $x_i$. Hence, setting $\Sigma_0=\Sigma$ and $\Sigma_l=\Delta$, for every $i=0,\ldots, l-1$ we have that $\dist_{\cuco Z}(f(\Sigma_i),f(\Sigma_{i+1}))\le D$, where $D=2C_0+2K$. 
\begin{claim}
    If $\lambda_2$ is large enough, each two consecutive simplices $\Sigma_i$ and $\Sigma_{i+1}$ must be $\mathcal{W}$--adjacent.
\end{claim} 
If this is true we are done, since then
$$\dist_{\mathcal{W}}(\Sigma, \Delta)\le \sum_{i=0}^{l-1}\dist_\mathcal{W}(\Sigma_i, \Sigma_{i+1})=l\le K \dist_{\cuco Z}(f(\Sigma), f(\Delta))+K.$$

To prove the claim notice that, by the distance formula Theorem~\ref{thm:distance_formula}, there exists a constant $M_0=M_0(D)$ such that, for every $U\in \frakS$, we have that $\dist_{\fontact U}(f(\Sigma_i),f(\Sigma_{i+1}))\le M_0$. Then in turn $\dist_{\fontact U}(b(\Sigma_i),b(\Sigma_{i+1}))\le M_0+2E$, and if we choose $\lambda_2\ge M_0+2E$ we have that $\Sigma_i$ and $\Sigma_{i+1}$ are $\mathcal{W}$--adjacent whenever $\lambda\ge \lambda_2$ (regardless of how their supports intersect, because $\lambda$ is the tightest threshold for the definition of a $\mathcal{W}$--edge).
\end{proof}

\section{Adding a group action}\label{sec:converse_group}
Here we make some remarks on the construction of the combinatorial HHS for the case when ${\cuco Z}$ is acted on by some finitely generated group $G$, in the sense of Subsection~\ref{subsec:HHG}. First, we want to show that the construction from Section~\ref{sec:construction} is $G$-equivariant, that is, the action of $G$ on ${\cuco Z}$ induces a “compatible” action on $(X,\mathcal{W})$. This will prove the “moreover” statement of Theorem~\ref{thm:main}. Then, in Theorem~\ref{thm:comb_hhg} we will prove that, if $G$ has a HHG structure coming from the action on ${\cuco Z}$ then it will also have a HHG structure coming from the action on $(X,\mathcal{W})$.

The following is the combinatorial counterpart of Definition~\ref{defn:action_on_HHS}:
\begin{defn}
    We say that a group $G$ \emph{acts} on the pair $(X,\mathcal{W})$, where $\mathcal{W}$ is an $X$-graph, if $G$ acts on $X$ by simplicial automorphisms, and the $G$–action on the set of maximal simplices of $X$ extends to an action on $\mathcal{W}$.
    \end{defn}
Here we show that if we start with a $G$-action on $({\cuco Z},\frakS)$ then the pair $(X,\mathcal{W})$ inherits a $G$-action:

\begin{thm}\label{thm:G-equivariance}
Let $({\cuco Z},\frakS)$ be a HHS with clean containers. Let $G$ be a finitely generated group acting on an HHS $({\cuco Z},\frakS)$ by automorphisms. Let $(X,\mathcal{W})$ be the graphs constructed in Subsections~\ref{subsec:minimal_orth_graph_and_blow_up} and~\ref{subsec:HHS_W}. Then $G$ acts on $(X,\mathcal{W})$. Moreover, the realisation map $f:\,\mathcal{W}\to {\cuco Z}$ from Definition~\ref{defn:realisation_map} is coarsely $G$-equivariant, meaning that for every $g\in G$ the following diagram coarsely commutes, with constants independent of $g$:
$$\begin{tikzcd}
\mathcal{W}\ar{r}{f}\ar{d}{g}&{\cuco Z}\ar{d}{g}\\
\mathcal{W}\ar{r}{f}&{\cuco Z}\\
\end{tikzcd}$$
\end{thm}

\begin{rem}
The space $({\cuco Z},\frakS)$ is only required to have clean containers since this is the only requirement to build the two graphs $(X,\mathcal{W})$, as pointed out in Remark~\ref{rem:only_clean}.
\end{rem}

\begin{proof}[Proof of Theorem~\ref{thm:G-equivariance}]
    The group $G$ acts on the domain set $\frakS$, preserving nesting and orthogonality, therefore there is an induced action on the minimal orthogonality graph $\bar X$, mapping every vertex $U$ to $g^\sharp(U)$. This action extends to the blow-up graph $X$ if, for every $\nest$-minimal domain $U$, every point $p\in\fontact U$ and every $g\in G$ we set 
    $$g(p):=g^\diamond(U)(p)\in \fontact (g^\sharp (U)).$$
    With a slight abuse of notation, we will denote the image of a vertex $v\in X$ under the action of an element $g\in G$ simply by $gv$. Similarly, we will drop the superscript for the action on $\bar X$ and set $gU:=g^\sharp (U)$.

    Now we show that the induced action on $\mathfrak M(X)$ preserves $\mathcal{W}$--edges. Let $\Sigma, \Delta$ be two maximal simplices and let $\bar\Sigma, \bar\Delta$ be their supports. Let $W=(\bar\Sigma\cap\bar\Delta)^\orth$, in the sense of Notation~\ref{notation:exceptional_orth} for the exceptional cases, and for every $g\in G$ let $W'=(g\bar\Sigma\cap g\bar\Delta)^\orth$. Clearly we have that $gW=W'$ since the action preserves orthogonality, and the co-level $k=\text{co-lv}(W)$ coincides with the co-level of $gW$ since the action preserves nesting.
    
    Now, if $\Sigma$ and $\Delta$ are $\mathcal{W}$--adjacent then the realisation tuples $b(\Sigma), b(\Delta)$ are $(k+1)\lambda$-close in every coordinate space. Thus, in order to prove that $g(\Sigma)$ and $g(\Delta)$ are again $\mathcal{W}$--adjacent, one needs to show that also $b(g\Sigma), b(g\Delta)$ are $(k+1)\lambda$-close in every coordinate space. Notice that we need the same threshold on distances, since the co-levels coincide.
    
    Recall that, if $\bar\Delta=\{U_1,\ldots, U_k\}$ then the realisation tuple of $\Delta$ is
    $$b(\Delta)_V=
    \begin{cases}
    x_i \quad\mbox{if }V=U_i\in \bar\Delta;\\
    \bigcup_{V\not\bot U_i}\rho^{U_i}_{V} \quad\mbox{otherwise}.
    \end{cases}$$
    If we apply $g$ to all coordinates we get
    $$g(b(\Delta)_V)=
    \begin{cases}
    g(x_i) \quad\mbox{if }gV=gU_i\in \bar\Delta;\\
    \bigcup_{gV\not\bot gU_i}g\rho^{U_i}_{V}=\bigcup_{gV\not\bot gU_i}\rho^{gU_i}_{gV} \quad\mbox{otherwise}.
    \end{cases}$$
    where we used that, as discussed in Remark~\ref{rem:commutativity_of_diag_in_hhg}, we can assume that $\rho^{gU}_{gV}=g\rho^{U}_{V}$ for every $g,U,V$. On the other hand the latter expression is the realisation tuple of $g\Delta$ by definition, thus we just showed that for every $V\in \frakS$ we have $b(g\Delta)_{gV}=g(b(\Delta)_V)$. Therefore
    $$\dist_{\fontact gV}(b(g\Delta)_{gV},b(g\Sigma)_{gV})=\dist_{\fontact gV}(g(b(\Delta)_V),g(b(\Sigma)_V))=\dist_{\fontact V}(b(\Delta)_{V},b(\Sigma)_{V}),$$
    where we used that the map $\fontact V\to \fontact (gV)$ induced by $g$ is an isometry. This means that if $b(\Sigma)$ and $ b(\Delta)$ are $(k+1)\lambda$-close in every coordinate space then so are $b(g\Sigma)$ and $b(g\Delta)$, and therefore $g$ preserves $\mathcal{W}$--adjacency.

    Finally, in order to prove that the realisation map $f:\,\mathcal{W}\to {\cuco Z}$ is coarsely $G$-equivariant we just note that, as proved above, for every $g\in G$ and every maximal simplex $\Delta\in \mathcal{W}^{(0)}$, we have that the tuple $b(g\Delta)$, which coarsely coincides with the coordinates of $f(g\Delta)$, is equal to $gb(\Delta)$. Moreover, as discussed in Remark~\ref{rem:commutativity_of_diag_in_hhg}, we can assume that for every $g\in G$, $V\in\frakS$ and $p\in {\cuco Z}$ we have that $\pi_{gV}(gp)=g\pi_{V}(p)$. This, applied to $p=f(\Delta)$, tells us that
    $$\pi_{gV}(gf(\Delta))=g\pi_{V}(f(\Delta))\sim g (b(\Delta)_V)=b(g\Delta)_{gV}\sim\pi_{gV}(f(g\Delta)),$$
    where $\sim$ denotes equality up to a bounded error. Thus the coordinates of $gf(\Delta)$ and $f(g\Delta)$ coarsely coincide in every coordinate space, and by the uniqueness axiom we have that $\dist_{{\cuco Z}}(gf(\Delta),f(g\Delta))\le E$.
\end{proof}

Next we turn our attention to actions with more structure:

\begin{thm}[{\cite{BHMS}}]\label{thm:CHHS_moreover}
    Let $(X,\mathcal{W})$ be a combinatorial HHS, and let $G$ be a group acting on $X$ with finitely many orbits of subcomplexes of the form $\link(\Delta)$, where $\Delta$ is a simplex of $X$. Suppose moreover that the action on maximal simplices of $X$ extends to an action on $\mathcal{W}$, which is metrically proper and cobounded. Then $G$ acts metrically properly and coboundedly on $\mathcal{W}$ and with finitely many $G$-orbits of domains, and therefore it is a HHG.
\end{thm}
\begin{proof}
    This is the “moreover” part of \cite[Theorem 1.18]{BHMS}. As stated there, the theorem requires the $G$-action on $X$ to be cocompact, but as discussed in \cite[Remark 1.19]{BHMS} the proof only uses that there are finitely many $G$-orbits of links of simplices.
\end{proof}

\begin{defn}
    We will say that a group $G$ satisfying the hypotheses of Theorem~\ref{thm:CHHS_moreover} is a \emph{combinatorial} hierarchically hyperbolic group, meaning that the HHS structure from Definition~\ref{defn:HHG} is inherited from the action on a combinatorial HHS.
\end{defn} 
Here we show that, under a mild assumption on the action of $G$ on $\frakS$, every HHG whose underlying HHS satisfies the hypotheses of Theorem~\ref{thm:XWf} is a combinatorial HHG.
\begin{thm}\label{thm:comb_hhg}
    Let $G$ be a hierarchically hyperbolic group, and let $({\cuco Z},\frakS)$ be a hierarchically hyperbolic space on which $G$ acts metrically properly and coboundedly. Suppose that $({\cuco Z},\frakS)$ has weak wedges, clean containers, the orthogonals for non-split domains property~\eqref{property:orthogonal_for_non_split}, and the DPR property~\eqref{property:dpr}. Moreover, suppose that the action $G\circlearrowleft \frakS$ has finitely many orbits of unordered tuples $\{V_1,\ldots, V_k\}$ of pairwise orthogonal elements, for every $k\le n$. Then $G$ acts on the pair $(X,\mathcal{W})$ constructed in Subsections~\ref{subsec:minimal_orth_graph_and_blow_up} and~\ref{subsec:HHS_W}, and the action satisfies the hypotheses of Theorem~\ref{thm:CHHS_moreover}. Hence $G$ is a combinatorial HHG.
\end{thm}

\begin{proof}
    Combining Theorem~\ref{thm:G-equivariance} and Theorem~\ref{thm:XWf} we get that $(X,\mathcal{W})$ is a combinatorial HHS, that $G$ acts on $(X,\mathcal{W})$, and that the realisation map $f:\mathcal{W}\to {\cuco Z}$ is a coarsely $G$-equivariant quasi-isometry. The latter fact already implies that $G$ acts metrically properly and coboundedly on $\mathcal{W}$, because the same properties hold for the $G$-action on ${\cuco Z}$.

    Then we are left to prove that the $G$-action on $X$ has a finite number of orbits of the form $\link(\Delta)$, where $\Delta$ is a simplex of $X$. Recall that, by Lemma~\ref{lem:decomposition_of_links}, the link of $\Delta$ is given by
    $$\link(\Delta)=p^{-1}(\link_{\bar X}(\bar\Delta))\star (\bigstar_{U\in\bar\Delta^{(0)}}\link_{p^{-1}(U)}(\Delta_U)).$$
    Let $\bar\Delta=\{U_1,\ldots, U_k\}$. Then $\link(\Delta)$ is uniquely determined by the tuple $\{U_1,\ldots, U_k\}$ and by the choice of $\link_{p^{-1}(U_i)}(\Delta_{U_i})$ for every $i$. The assumption on the action of $G$ on $\frakS$ tells us that there is a finite number of $G$-orbits of tuples of the form $\{U_1,\ldots, U_k\}$, because these elements are pairwise orthogonal by definition. Moreover, given such a tuple, for every $i$ we have that
    $\link_{p^{-1}(U_i)}(\Delta_{U_i})$ can only be one of the following:
    \begin{itemize}
        \item If $\Delta_{U_i}=v_{U_i}$ is the tip of the cone then $\link_{p^{-1}(U_i)}(\Delta_{U_i})=(\fontact U_i)^{(0)}$ is the base of the cone;
        \item If $\Delta_{U_i}$ is a point in the base then $\link_{p^{-1}(U_i)}(\Delta_{U_i})=v_{U_i}$;
        \item If $\Delta_{U_i}$ is an edge then $\link_{p^{-1}(U_i)}(\Delta_{U_i})=\emptyset$
    \end{itemize}
    Therefore there are three possible choices for every $U_i$, and this concludes the proof that $G\circlearrowleft X$ has finitely many orbits of the form $\link(\Delta)$.
\end{proof}

\begin{rem}
    The requirement that there are finitely many $G$-orbit of pairwise orthogonal domains (which we just refer to as \emph{the requirement} in this Remark) is necessary for our combinatorial HHS structure to be a hierarchically hyperbolic \emph{group} structure. To see this, let $U_1,\ldots, U_n$ be a collection of pairwise orthogonal, $\nest$-minimal domains. Then, by how the combinatorial HHS $(X,\mathcal W)$ was defined in Section~\ref{sec:construction}, we see that the simplex $\Sigma=\{v_{U_1},\ldots, v_{U_n}\}$ is the link of any simplex $\Delta$ such that $\Sigma\star\Delta$ is maximal, and in particular $\Sigma$ is a domain in the combinatorial HHS structure. Moreover, by definition a hierarchically hyperbolic group structure for $G$ has finitely many $G$-orbits of domains, so there must be finitely many orbits of simplices corresponding to $\nest$-minimal, pairwise orthogonal domains. 

    The requirement is ``easily satisfied'', in the sense that, by \cite[Lem. 6.6]{BGGH:free-by-Z}, it holds when $(G,\frakS)$ has the natural property that $\stabilizer_G(V)$ acts on $P_V$ coboundedly for all $V\in \frakS$. On the other hand, the latter condition is stronger and is not a consequence of the definition of an HHG since, for example, one can put exotic HHG structures on a free group where this fails.

    If we remove the requirement, the same proof endows a HHG $(G,\frakS)$ with a combinatorial \emph{$G$-HHS} structure, which is defined as a combinatorial HHG structure except that we only require finitely many orbits of \emph{unbounded} domains. Indeed, by inspection of Corollary~\ref{cor:bounded_links}, we see that unbounded augmented links of simplices of $X$ correspond to either $\nest$-minimal domains in $\frakS$, or to the orthogonal complements of certain collections of domains in $\frakS$. In particular, unbounded domains in $(X,\mathcal W)$ correspond to certain domains in $\frakS$, of which there are finitely many $G$-orbits since $(G,\frakS)$ is a HHG.
\end{rem}

\section{Some other hypotheses}\label{sec:additional_hyp}
The hypotheses of Theorem~\ref{thm:main} are rather technical, so in this section we present some more “natural” ones, and we describe how they relate to each other and to the original ones. We also show that, in certain cases, one can match the requirements of Theorem~\ref{thm:main} by adding (quite a lot of) bounded coordinate spaces.

We start with some possible requirements on the domain set.

\begin{property}\label{property:strong_orth} A hierarchically hyperbolic space has the \emph{strong orthogonal property} if for every two domains $U\propnest V$ there exists $W\propnest V$ such that $U\orth W$.
\end{property} 

\begin{property}\label{property:weak_orth} A hierarchically hyperbolic space has the \emph{weak orthogonal property} if for every two domains $U\propnest V$, if $U$ is non-$\nest$-minimal there exists $W\propnest V$ such that $U\orth W$.
\end{property} 

Both these properties imply the orthogonals for non-split domains property~\eqref{property:orthogonal_for_non_split}, since every $\nest$-minimal domain is trivially split (it coincides with its unique Samaritan).

\begin{property}\label{property:bounded_C_for_split}
A hierarchically hyperbolic space $({\cuco Z},\frakS)$ has \emph{bounded split coordinate spaces} if there exists $c\ge 0$ such that, for every $U\in\frakS-{S}$ which is split and non-$\nest$-minimal, the corresponding coordinate space $\fontact U$ has diameter at most $c$.
\end{property}

\begin{lemma}[Comparison between properties]\label{lem:comparison_between_prop}
Let $({\cuco Z},\frakS)$ be a normalised hierarchically hyperbolic space with wedges and clean containers. Then each of the following properties implies the lower ones, meaning that if $({\cuco Z},\frakS)$ has some property (i) from the list, and $j > i$, then there exists a normalised HHS structure $({\cuco Z},\frakS')$, with $\frakS\subset\frakS'$ and satisfying wedges, clean containers and property (j):

\begin{enumerate}
    \item $({\cuco Z},\frakS)$ has the strong orthogonal property~\eqref{property:strong_orth};
    \item $({\cuco Z},\frakS)$ has the weak orthogonal property ~\eqref{property:weak_orth};
    \item $({\cuco Z},\frakS)$ has the weak orthogonal property ~\eqref{property:weak_orth} and dense product regions~\eqref{property:dpr};
    \item\label{item:ons_bsc} $({\cuco Z},\frakS)$ has orthogonals for non-split domains~\eqref{property:orthogonal_for_non_split} and bounded split coordinate spaces~\eqref{property:bounded_C_for_split};
    \item $({\cuco Z},\frakS)$ has orthogonals for non-split domains~\eqref{property:orthogonal_for_non_split} and dense product regions~\eqref{property:dpr}.
\end{enumerate}
\end{lemma}

\begin{proof}[Proof of Lemma~\ref{lem:comparison_between_prop}]
    The implication $1\Rightarrow 2$ is trivial. $3\Rightarrow 4$ follows from Lemma~\ref{lem:dpr_to_BSC} and the fact that $\nest$-minimal domains are split. The implication $2\Rightarrow 3$ is Lemma~\ref{lem:ensure_dpr_for_weak_orth}, while Lemma~\ref{lem:ensure_dpr_for_bscg} is $4\Rightarrow 5$. 
\end{proof}

\begin{lemma}\label{lem:dpr_to_BSC}
If a normalised HHS $({\cuco Z},\frakS)$ has the DPR property~\eqref{property:dpr} then it has the bounded split coordinate spaces property~\eqref{property:bounded_C_for_split}.
\end{lemma}
We should think of the bounded split coordinate space property as the fact that product regions are dense in split domains, and this is morally why Lemma~\ref{lem:dpr_to_BSC} holds.

\begin{proof}
    Let $U$ be a non-minimal split domain and let $W$ be one of its Samaritans. By the DPR property, for every $q\in \fontact U$ there exists a minimal domain $V\propnest U$ such that $\dist_{\fontact U} \left(q,\rho^V_U\right)\le M_0$. Now it suffices to notice that, since $W$ is a Samaritan, we must have that either $V=W$ or $V\orth W$, and by Lemma~\ref{lem:close_proj_of_orthogonals} we have that $\diam_{\fontact U}(\rho^W_U,\rho^{V}_U)\le 10E$. Therefore $q$ is $(M_0+10E)$-close to $\rho^W_U$, and therefore $\fontact U$ has diameter at most $2(M_0+10E)$.
    \end{proof}

The following lemmas show that, if one starts with more general hypotheses, there is often a way to modify the HHS structure in order to ensure the DPR property.

\begin{lemma}\label{lem:ensure_dpr_for_bscg}
Let $({\cuco Z},\frakS)$ be a normalised HHS with wedges, clean containers and satisfying the bounded split coordinate spaces property~\eqref{property:bounded_C_for_split} and the orthogonals for non-split domains property~\eqref{property:orthogonal_for_non_split}. Then there exists a normalised HHS structure $({\cuco Z},\frakS')$ such that $\frakS\subseteq\frakS'$ and $({\cuco Z},\frakS')$ has wedges, clean containers, the DPR property~\eqref{property:dpr} and the orthogonals for non-split domains property~\eqref{property:orthogonal_for_non_split}. 
\end{lemma}

\begin{lemma}\label{lem:ensure_dpr_for_weak_orth}
Let $({\cuco Z},\frakS)$ be a normalised HHS with wedges, clean containers and the weak orthogonal property~\eqref{property:weak_orth}. Then there exists a normalised HHS structure $({\cuco Z},\frakS')$ such that $\frakS\subseteq\frakS'$ and $({\cuco Z},\frakS')$ has wedges, clean containers, the weak orthogonal property~\eqref{property:weak_orth} and the DPR property~\eqref{property:dpr}. 
\end{lemma}

The strategy for proving these two lemmas is the same, so we present an extensive proof only of Lemma~\ref{lem:ensure_dpr_for_bscg}, which is more complicated. Here is a list of the changes which are needed to prove Lemma~\ref{lem:ensure_dpr_for_weak_orth}:
\begin{itemize}
    \item A domain $T^U_x$ must be added for every non-minimal $U$ (thus not only if $U$ is non-split);
    \item The argument below to show the DPR property for wide domains will apply to all non-$\nest$-minimal domains; 
    \item Since $\frakS$ and $\frakS'$ will have the same non-$\nest$-minimal domains, the weak orthogonal property will be preserved.
\end{itemize}

\begin{proof}[Proof of Lemma~\ref{lem:ensure_dpr_for_bscg}]
What we will actually prove is that, if a normalised HHS $({\cuco Z},\frakS)$ has the bounded split coordinate space property~\eqref{property:bounded_C_for_split} then we can find a structure $({\cuco Z},\frakS')$ with the DPR property. Moreover, if $({\cuco Z},\frakS)$ has wedges, clean containers or the orthogonals for non-split domains property then the procedure will preserve these properties.

We will say that a domain $U\in \frakS$ is \emph{wide} if $U$ is non-split (and in particular non-$\nest$-minimal) or $U=S$. For every wide domain $U$ we do the following. Recall that we defined $F_U$ as the space of all $20E$-consistent tuples. For every $x=(x_V)_{V\nest U}\in F_U$ we define a domain $T^U_x$ whose coordinate space $\fontact T^U_x$ is a point, and we let $\frakS'$ be the union of $\frakS$ and these new domains. Now we show that $({\cuco Z},\frakS')$ is an HHS, defining new projections, relative projections and relations when needed.
\par\medskip

\textbf{Projections}: The projection $\pi_{T^U_x}:\,{\cuco Z}\to \fontact T^U_x$ is just the constant map, while all other projections are inherited from the original structure.
\par\medskip

\textbf{Nesting}: The domains $T^U_x$ are $\nest$-minimal, and $T^U_x\nest V$ if and only if $U\nest V$. If $U\propnest V$ then we define the projection $\rho^{T^U_x}_V=\rho^U_V$; moreover we set $\rho^{T^U_x}_U=x_U$. The projections in the opposite direction (namely, $\rho^V_{T^U_x}$ whenever $U\nest V$) are just the constant maps.
\par\medskip

\textbf{Finite complexity}: Since we just added $\nest$-minimal domains inside non-$\nest$-minimal ones, the complexity of the HHS structure remains the same.
\par\medskip

\textbf{Orthogonality and clean containers}: If $U,V\in \frakS$, $U$ is wide and $x\in F_U$, we say that $V\orth T^U_x$ if and only if $U\orth V$. Moreover, if $V$ is also wide and $y\in F_V$ we say that $ T^U_x\orth T^V_y$ if and only if $U\orth V$.

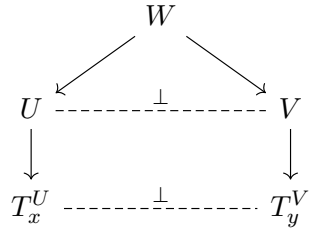
\begin{figure}[htp]
    \centering
\begin{tikzcd}
    &W\ar{dr}\ar{dl}&\\
    U\ar{d} \ar[rr, "\orth", dashed, no head]&&V\ar{d}\\
    T^U_x\ar[rr, "\orth", dashed, no head]&&T^V_y
\end{tikzcd}
    \caption{The diagram, where arrows denote nesting, shows the only way two newly added domains $T^U_x$ and $T^V_y$ can be orthogonal inside some $W\in \frakS$.}
    \label{fig:orthogonal_new_domains}
\end{figure}

Notice that containers already exist for any situation involving elements of $\frakS'- \frakS$. Indeed, suppose $T^U_x\nest W$ for some $W\in\frakS$, which implies that either $U=W$ or $U\propnest W$. In the first case, no container is needed, since $T^U_x$ is transverse to every $V\propnest U$ (and therefore also to every $T^V_y$ if $V$ is wide). In the second case, every $V\in\frakS$ which is properly nested inside $W$ and orthogonal to $T^U_x$ is also orthogonal to $U$, and therefore already nested inside the container for $U$ inside $W$. Moreover, if $T^V_y\propnest W$ and $T^U_x\orth T^V_y$ then $V\propnest W$ and $U\orth W$, which means that $T^V_y$ is already nested inside the container for $U$ inside $W$.

Conversely, if $V\in\frakS$ and $T^U_x$ are orthogonal and properly nested in $W$ then $U\propnest W$ and $U\orth V$, thus $T^U_x$ is already nested inside the container for $V$ inside $W$.

Now, if $({\cuco Z},\frakS)$ has clean containers then so does $({\cuco Z},\frakS')$. This is because, as argued above, the container for $T^U_x$ inside some $W$ is the container for $U$ inside $W$, and if this container is orthogonal to $U$ then it is also orthogonal to $T^U_x$.
\par\medskip

\textbf{Transversality}: If $U,V\in \frakS$, $U$ is wide and $x\in F_U$, then by construction $T^U_x\transverse V$ if and only if one of the following holds:
\begin{itemize}
    \item $U\transverse V$: in this case we define $\rho^{T^U_x}_V=\rho^U_V$;
    \item $V\propnest U$: in this case we set $\rho^{T^U_x}_V= x_V$.
\end{itemize}
Moreover, whenever $Y\in\frakS'$ is transverse to $T^U_x$ we set $\rho^V_{T^U_x}=\fontact T^U_x$.
\par\medskip

\textbf{Consistency}: Since the only elements of $\frakS'$ whose coordinate spaces are not points are in $\frakS$, the first two consistency inequalities are trivial. Moreover, the final clause of the consistency axiom holds, since whenever $U\propnest V$ or $U\transverse V$ we defined $\rho^{T^U_x}_V$ to be $\rho^U_V$.
\par\medskip

\textbf{Uniqueness, BGI, large links}: Since the only elements of $\frakS'$ whose coordinate spaces are not points are in $\frakS$, these axioms for $({\cuco Z},\frakS')$ follow from the corresponding ones for $({\cuco Z},\frakS)$.
\par\medskip

\textbf{Partial realisation}: Let $V_1,\ldots, V_k\in \frakS'$ be pairwise orthogonal elements, and let $p_i\in\fontact V_i$. We show that we can find a partial realisation point for $\{(V_i, p_i)\}$. Up to permutation we can assume that $V_1,\ldots, V_l\in \frakS$, for some $l\le k$, and $V_i=T^{U_i}_{x_i}$ for every $l<i\le k$. Moreover, we can assume that the family $V_1,\ldots, V_k$ is maximal, up to adding domains belonging to $\frakS$, since a realisation point for a bigger family is also a realisation point for the original one.

Now, since $\cuco Z$ is normalised, for every $i\le l$ we can find $y_i\in F_{V_i}$ such that $\dist_{V_i}(\pi_{V_i}(y_i),p_i)$ is uniformly bounded. Then, for any $U\in \frakS'$ set

\begin{equation}\label{eq:qu}
q_U=\begin{cases}
    (y_i)_U \mbox{ if }U\nest V_i,\,i\le l;\\
    (x_i)_U \mbox{ if }U\nest U_i,\,l<i\le k;\\
    \bigcup_{V_i\transverse U\mbox{ or }V_i\propnest U} \rho^{V_i}_U\mbox{ otherwise}.
\end{cases}\end{equation}
Notice that, since $V_1,\ldots, V_k$ is a maximal family, the coordinate $q_U$ is always well-defined. Moreover, it is easy to see that $(q_U)_{U\in \frakS'}$ is $20E$-consistent, thus by the realisation Theorem~\ref{thm:realisation} there exists $z\in \cuco Z$ whose coordinates are uniformly close to $(q_U)_{U\in \frakS'}$. It is also easy to verify that $z$ is a partial realisation point for $\{(V_i, p_i)\}$, since by Equation~\ref{eq:qu} it has the right coordinates whenever $U=V_i$, $V_i\transverse U$ or $V_i\propnest U$ for some $i\le k$.
\par\medskip

\textbf{DPR}: Let $U\in\frakS$ be a non-minimal domain and let $p\in \fontact U$. If $U$ is non-split then there is $z\in {\cuco Z}$ such that $\pi_U(z)$ is uniformly close to $p$, by our normalisation assumption. Thus, if we set $x=(\pi_V(z))_{V\nest U}$, we have that $T^U_x$ projects uniformly close to $p$ in $\fontact U$, and the DPR property holds with the same constant as the one coming from the normalisation assumption. On the other hand, if $U$ is split then $\fontact U$ is uniformly bounded by the bounded split coordinate space property~\eqref{property:bounded_C_for_split}. Hence $p\in \fontact U$ is uniformly close to $\rho^V_U$ for any $V\propnest U$.
\par\medskip

\textbf{Additional properties}: We now check that if $({\cuco Z},\frakS)$ has one of the properties below then so does $({\cuco Z},\frakS')$.
\begin{itemize}
\item \textbf{Wedges}: The new domains are all $\nest$-minimal, therefore for every $T^U_x,V'\in\frakS'$ we can set
$$T^U_x\wedge V'=\begin{cases}
    T^U_x\mbox{ if }T^U_x\nest V';\\
    \emptyset\mbox{ otherwise}.
\end{cases}$$
Then we just need to verify that, if $V,W\in \frakS$ and there exists $R\in\frakS'$ which is nested inside both, then the wedge of $V$ and $W$ exists in $\frakS'$. What we will actually show is that $V$ and $W$ already have a well-defined wedge inside $\frakS$, call this wedge $Q$, and that whenever $R\in\frakS'$ is nested inside both $V$ and $W$ then $R\nest Q$. Therefore the wedge in $\frakS'$ will coincide with the wedge in $\frakS$.

Now let $R\in \frakS'$ as above. If $R\in \frakS$ then $Q$ exists and $R\nest Q$, by the wedge property for $\frakS$. Otherwise $R=T^U_x\propnest U$ for some $U$ and some $x\in F_U$, and by definition of the new nesting relations $U$ must be already nested in both $V$ and $W$. Hence again $Q$ exists and $R\nest U\nest Q$.
\item \textbf{Orthogonals for non-split domains property~\eqref{property:orthogonal_for_non_split}}: 
Let $U\propnest V$, with $V\in \frakS$ since the new domains are all minimal. If $U$ is one of the new domains then $U$ is $\nest$-minimal, and therefore split. Hence, suppose that $U\in \frakS$. If $U$ is non-split in $\frakS$ then property~\eqref{property:orthogonal_for_non_split}, which holds for $\frakS$, ensures the existence of some $W\orth U$ inside $V$. Otherwise $U$ is split in $\frakS$ with some Samaritan $W$, and we claim that $U$ is again split in $\frakS'$ with the same Samaritan. Let $Q\in\frakS'$ be such that $Q\nest U$. If $Q\in \frakS$ then $W\nest Q$ or $W\orth Q$, by definition of Samaritan, and we have nothing to prove. Otherwise we have that $Q=T^{R}_x\nest R\nest U$ for some non-split domain $R\in \frakS$ and some $x\in F_R$. We cannot have that $W\nest R$, since then $R$ would be split in $\frakS$ by Remark~\ref{rem:middle_in_split}; thus $R$ (and therefore $Q$) is orthogonal to $W$.
\end{itemize}
This concludes the proof.\end{proof}

\begin{rem}
    The weak wedge property is not preserved by the procedure of Lemma~\ref{lem:ensure_dpr_for_bscg} to ensure the DPR property, since we are adding many “loose” $\nest$-minimal domains that will not be nested in the original weak wedges. This is why in this section we are assuming the “strong” wedge property.
\end{rem}

\begin{rem}[DPR for HHG with cobounded actions]\label{rem:dpr_for_HHG}
    If $G$ is a HHG, but the underlying space $({\cuco Z},\frakS)$ does not have the DPR property already, we cannot argue as in Lemma~\ref{lem:ensure_dpr_for_bscg} to enforce it. This is because, if we add a domain of the form $T^U_x$ whenever $x\in F_U$ and we define the $G$-action on $\frakS'$ in the obvious way (that is, by setting $gT^U_x=T^{gU}_{gx}$), then this action cannot have finitely many $G$-orbits of domains, since $G$ is countable while $F_U$ might be uncountable.
    
    Nonetheless, the DPR property is fairly common for HHG. For instance, it is implied by the fact that, for every $U\in \frakS$, the stabiliser $\text{Stab}_G(U)$ acts coboundedly on $\fontact U$. Indeed, if $V$ is $\nest$-minimal and $V\propnest U$, then the collection $\{\rho^{gV}_U\}_{g\in \text{Stab}_G(U)}$ coincides with the $\text{Stab}_G(U)$-orbit of $\rho^V_U$, by definition of a HHG, and the latter is coarsely dense in $\fontact U$ by coboundedness.

    In turn, if the HHS structure $({\cuco Z},\frakS)$ is normalised, we can further reduce to checking the weaker requirement that, for every $U\in \frakS$, the stabiliser $\text{Stab}_G(U)$ acts coboundedly on the product region $P_U$. Indeed, the projection $\pi_U:\,P_U\to \fontact U$ is coarsely surjective (by normalisation), coarsely Lipschitz (by Definition~\ref{defn:HHS} of a HHS) and $\text{Stab}_G(U)$-equivariant (by Definition~\ref{defn:action_on_HHS} of the $G$-action), therefore if $\text{Stab}_G(U)$ acts coboundedly on $P_U$ then it also acts coboundedly on $\fontact U$.

    In practice, the latter requirement is not particularly restrictive, since all "reasonable" HHGs have cobounded actions on their product regions, and all known methods of producing new HHGs tend to preserve this property (unless one comes up with some very artificial structures).
\end{rem}

\subsection{Orthogonal sets}\label{subsec:orthogonal_sets}
Another property that one could require on the index set is that orthogonal complementation is an involution. Such a property is satisfied by CAT(0) cube complexes with (weak) factor systems (see Subsection~\ref{subsec:CCC}) and is implied by the strong orthogonality property~\eqref{property:strong_orth} (Lemma~\ref{lem:strong_orth_implies_orthdetnest}). However, it goes in a somewhat different direction than the orthogonals for non-split domains property~\eqref{property:orthogonal_for_non_split}, in the sense that is not enough to prove Theorem~\ref{thm:main}. Indeed, in Section~\ref{subsec:counterex_CCC} we will provide an example of an HHS $\cuco Z$ with wedges, clean containers, dense product regions and where orthogonal complementation is an involution, but such that the graph $X$ from Definition~\ref{defn:blow-up} cannot be the support of a CHHS structure for $\cuco Z$.

\begin{defn}\label{defn:orth_set}
A partially ordered set $(\mathfrak F,\nest)$ is called \emph{orthogonal} if there exists a symmetric relation $\orth$ on $\mathfrak F$ such that the following hold for all $U,V,W\in \mathfrak F$:
\begin{itemize}
\item $U\not\bot U$;
\item $\mathfrak F$ has a unique $\nest$-maximal element $S$;
\item if $U\nest V$ and $V\orth W$ then $U\orth W$;
\item (\textbf{Wedges.}) if $W\nest U,V$, then there exists $U\wedge V\in \mathfrak F$ such that $U\wedge V\nest U,V$, and for all $W\nest U,V$ we have that $W\nest U\wedge V$;
\item (\textbf{Clean containers.}) for all $U$ such that there exists $V\orth U$, there exists $U^\orth\in \mathfrak F$ such that, for all $V\orth U$ we have that $V\nest U^\orth$, and $W\orth U^\orth$ if and only if $W\nest U$;
\item (\textbf{Orthogonality determines nesting.}) $U\nest V$ (resp. $U\propnest V$) if and only if the set of $W$ for which $V\orth W$ is contained (resp. properly contained) in the set of $W'$ for which $U\orth W'$. In particular, if nothing is $\orth$-related to $V$ then $V$ is the unique $\nest$-maximal element, while if there exists $W\orth V$ then $V^\orth\nest U^\orth$ (resp. $V^\orth\propnest U^\orth$) if and only if $U\nest V$ (resp. $U\propnest V$).
\end{itemize}
\end{defn}

\begin{lemma}[Orthogonality determines nesting for HHS]\label{lem:orth_implies_nest_equiv}
    Let $(\cuco Z, \frakS)$ be an HHS with wedges and clean containers. Then the following are equivalent:
    \begin{itemize}
        \item (\textbf{Complementation is an involution.}) For all $U\in \frakS-\{S\}$, $U^\orth$ is defined and $U^{\orth\orth}=U$.
        \item (\textbf{Orthogonality determines nesting.}) For all $U,V\in \frakS-\{S\}$, we have that $U\propnest V$ if and only if $V^\orth\propnest U^\orth$.
    \end{itemize}
    If one (hence both) of the previous holds, then $\frakS$ is an orthogonal set.
\end{lemma}

\begin{proof}
    The first part of the Lemma is proven exactly as \cite[Proposition 6.1]{Rcubes} (the statement there is for real cubings, but as pointed out in \cite[Remark 6.2]{Rcubes} the same argument works for HHSs). The second part follows from the properties of the domain set $\frakS$ of a HHS (see Definition~\ref{defn:HHS}).
\end{proof}

\begin{lemma}\label{lem:strong_orth_implies_orthdetnest}
    Let $(\cuco Z, \frakS)$ be an HHS with wedges, clean containers and the strong orthogonal property~\eqref{property:strong_orth}. Then $\frakS$ is an orthogonal set.
\end{lemma}

\begin{proof}
    Just notice that, for every $U\in \frakS-\{S\}$, the strong orthogonal property grants the existence of $U^\orth$; moreover $U= U^{\orth\orth}$, because otherwise we could find a $V\nest U^{\orth\orth}$ which is orthogonal to $U$, and this would contradict the definition of $U^\orth$. Now the conclusion follows from Lemma~\ref{lem:orth_implies_nest_equiv}.
\end{proof}

\section{Near equivalence of HHS and combinatorial HHS}\label{sec:HHS_iff_CHHS}
In this section we show that, if the hypotheses on $({\cuco Z},\frakS)$ are the strongest possible, then the combinatorial HHS $(X,\mathcal{W})$ arising from the construction has the following two nice properties:

\begin{defn}\label{defn:simp_cont} A combinatorial HHS $(X,\mathcal{W})$ has \emph{simplicial containers} if for any simplex $\Delta\subset X$ there exists a simplex $\Phi\subset X$ such that $$\link(\link(\Delta))=\link(\Phi).$$
\end{defn}

\begin{defn}\label{defn:simp_wedge} A combinatorial HHS $(X,\mathcal{W})$ has \emph{simplicial wedges} if for any two simplices $\Delta, \Sigma\subset X$ there exists a simplex $\Pi$ which extend $\Sigma$ such that $$\link(\Delta)\cap\link(\Sigma)=\link(\Pi).$$
\end{defn}

\begin{thm}\label{thm:strong_iff}
Let $({\cuco Z},\frakS)$ be a normalised hierarchically hyperbolic space. Then ${\cuco Z}$ has wedges, clean containers and the strong orthogonal property~\eqref{property:strong_orth} if and only if there exists a CHHS $(X,\mathcal{W})$ with simplicial wedges and simplicial containers such that $\mathcal{W}$ is quasi-isometric to ${\cuco Z}$.
\end{thm}

\begin{proof}
First we show that, if $(X,\mathcal{W})$ has simplicial wedges and simplicial containers then $\mathcal{W}$, with the HHS structure described in Subsection~\ref{subsec:W_structure}, has the following properties:
\begin{itemize}
    \item \textbf{Wedges}: given two non-maximal simplices $\Sigma,\Delta$, if there exists a simplex $\Gamma$ such that $[\Gamma]\nest[\Sigma]$ and $[\Gamma]\nest[\Delta]$, then $$\link(\Gamma)\subseteq \link(\Sigma)\cap\link(\Delta)=\link(\Pi)$$
    for some $\Pi$ depending only on $\Sigma, \Delta$. Therefore $[\Sigma]\wedge[\Delta]=[\Pi]$.
    \item \textbf{Strong orthogonal property}: Let $[\Delta]\propnest[\Delta']$. Now 
    $$ \link(\link(\Delta))\cap\link(\Delta')=\link(\Phi)\cap\link(\Delta')=\link(\Pi)$$
    for some simplices $\Phi,\Pi$ whose existence is granted by simplicial containers and wedges, respectively. Thus $[\Pi]\nest[\Delta']$, and since $\link(\Pi)\subseteq\link(\link(\Delta))$ we have that $[\Pi]\orth[\Delta]$.
    \item \textbf{Clean containers}: Let $[\Delta]\propnest[\Delta']$ and suppose that there exists $[\Sigma]\propnest[\Delta']$ which is orthogonal to $[\Delta]$. By definition
    $$\link(\Sigma)\subseteq \link(\link(\Delta))\cap\link(\Delta')=\link(\Pi),$$
    where $\Pi$ is the simplex defined above. Then $[\Sigma]\propnest[\Pi]$, which means that $[\Pi]$ is the container for $[\Delta]$ inside $[\Delta']$. Moreover, since $[\Pi]\orth[\Delta]$, this container is also clean.
\end{itemize}

Now we turn our attention to the converse statement. Let $({\cuco Z},\frakS)$ be a HHS with wedges, clean containers and the strong orthogonal property. By Lemma~\ref{lem:ensure_dpr_for_weak_orth} we can find a structure $({\cuco Z},\frakS')$ with wedges, clean, containers, the weak orthogonal property and the DPR property. Then Theorem~\ref{thm:main} applies to $({\cuco Z},\frakS')$ and outputs a combinatorial HHS $(X,\mathcal{W})$ where $\mathcal{W}$ is quasi-isometric to ${\cuco Z}$.
\par\medskip

\textbf{$(X,\mathcal{W})$ has simplicial wedges}: Let $\bar\Sigma$ and $\bar\Delta$ be the supports of $\Sigma$ and $\Delta$, respectively, and let $\bar\Sigma^\orth, \bar\Delta^\orth\in \frakS$ be their orthogonal complements. Let $\bar\Phi=\bar \Delta\cap \link(\bar \Sigma)$, and let $Y_0$ be the orthogonal complement of $\bar\Phi$ inside $\bar\Sigma^\orth$, that is, $Y_0=(\bar\Sigma\star\bar\Phi)^\orth$. Notice that $Y_0$ cannot be one of the minimal domains $T^U_x$ that were added in Lemma~\ref{lem:ensure_dpr_for_weak_orth} to ensure the DPR property. This is because, if $Y_0=T^U_x$, then $T^U_x$ is orthogonal to $\bar\Sigma\star\bar\Phi$, and by construction $U$ is orthogonal to the same simplex. But by definition of $Y_0$ we must have that $U\nest Y_0=T^U_x\propnest U$, which is a contradiction. 
\par\medskip

\textit{Part 1}: If $(\bar\Sigma)^\orth$ and $(\bar\Delta)^\orth$ don't have any nested domain in common, we formally set $W_0=\emptyset$ and we skip to Part 2. Otherwise we can consider the wedge $W_0=(\bar\Sigma)^\orth\wedge(\bar\Delta)^\orth$. Notice that $W_0$ cannot be one of the minimal domains from Lemma~\ref{lem:ensure_dpr_for_weak_orth}. This is because, if $W_0=T^U_x$, then $T^U_x$ is orthogonal to both $\bar\Sigma$ and $\bar\Delta$, and by construction $U$ is also orthogonal to the two simplices. But by definition of $W_0$ as a wedge we must have that $U\nest W_0=T^U_x\propnest U$, which is a contradiction. This means that $W_0\in \frakS$ as well. 
\par\medskip

\textit{Part 2}: Now suppose that $W_0$ has been defined as in Part 1. If $W_0=Y_0$ then we set $\bar\Theta=\emptyset$. Otherwise, by the strong orthogonal property, which holds for every two elements of the original domain set $\frakS$, we can find a $\nest$-minimal domain $V_0\nest Y_0$ such that $V_0\orth W_0$ (if $W_0\neq\emptyset$). Now let $Y_1=\{V_0\}^\orth_{Y_0}$, which is again in $\frakS$ and contains $W_0$, and we can argue as above. 

In both cases  we can find a (possibly empty) simplex $\bar\Theta=\{V_0,\ldots, V_k\}$ of $\bar X$ such that $W_0=(\bar\Sigma\star\bar\Phi\star\bar\Theta)^\orth$, and this readily implies that
\begin{equation}\label{eq:lklk_strong}
    \link(\bar\Sigma)\cap\link(\bar\Delta)=\link(\bar\Sigma\star\bar\Phi\star\bar\Theta).
\end{equation}

From now on we can argue exactly as in Lemma~\ref{lem:simplicial_wedge_property} (more precisely, we can repeat the arguments of the paragraph “\textbf{Finding the extension of $\Sigma$}”), and find a simplex supported in $\bar\Sigma\star\bar\Phi\star\bar\Theta$ which extends $\Sigma$ and whose link is $\link(\Sigma)\cap\link(\Delta)$.
\par\medskip

\textbf{$(X,\mathcal{W})$ has simplicial containers}: Let $\Delta\subseteq X$ be a simplex and let $\bar\Delta$ be its support. Let $\bar\Delta_1 \subset\bar\Delta$ be the domains $U$ such that $\Delta_U$ is a single point, and let $\bar\Delta_2\subset\bar\Delta$ be the domains such that $\Delta_U$ is an edge. Taking the link of the expression in Lemma~\ref{lem:decomposition_of_links}, which described the shape of the link of $\Delta$ inside $X$, we have that
$$\link(\link(\Delta))=p^{-1}(\link(\link(\bar\Delta)))\cap\bigcap_{U\in\bar\Delta_1}\link_X(\link_{p^{-1}(U)}(\Delta_U)).$$
Thus $v\in X$ belongs to $\link(\link(\Delta))$ if its support $V$ lies in $(\link(\link(\bar\Delta)))$, and either $V\orth \bar\Delta_1$, or $V\in\bar\Delta_1$ and $v\in \link_{p^{-1}(V)}(\link_{p^{-1}(V)}(\Delta_V)) $. Therefore, let $W=(\bar\Delta^\orth)^\orth$. Again, $W$ does not coincide with any $T^U_x$, because if $T^U_x= W$ then $T^U_x\orth \bar\Delta^\orth$, and therefore also $U\orth \bar\Delta^\orth$. Then again $W\in\frakS$, and the strong orthogonal property implies that there exists a simplex $\bar\Theta$ inside $\bar X$ such that $W=\bar\Theta^\orth$. Notice that every $U\in\bar\Delta$ is nested in $W$ by construction. Now define a simplex $\Phi$ with support $\bar\Theta\star\bar\Delta_1$ by choosing an edge for every domain $U$ which is a vertex of $\bar\Theta$, and a point $q_U\in\link_{p^{-1}(U)}(\Delta_U)$ for every $U\in \bar\Delta_1$. Thus by construction
$$\link(\Phi)=p^{-1}(\link(\bar\Theta\star\bar\Delta_1)))\star\bigcup_{U\in\bar\Delta_1}\link_{p^{-1}(U)}(q_U)=$$
$$=p^{-1}(\link(\link(\bar\Delta)))\cap\bigcap_{U\in\bar\Delta_1}\link_X(\link_{p^{-1}(U)}(\Delta_U))=\link(\link(\Delta)),$$
and we are done. 
\end{proof}

\begin{rem}\label{rem:already_has_DPR}
    If in the previous proof the original structure $({\cuco Z},\frakS)$ already has the DPR property~\eqref{property:dpr} then there is no need to invoke Lemma~\ref{lem:ensure_dpr_for_weak_orth}, and the whole proof of Theorem~\ref{thm:strong_iff} works with $\frakS$ instead of $\frakS'$. In other words, if $({\cuco Z},\frakS)$ already has the DPR property then the combinatorial HHS $(X, \mathcal{W})$ is exactly the one constructed in Section~\ref{sec:construction}.
\end{rem}

\section{Mapping class groups are combinatorial HHS}\label{sec:mcg}
Throughout this section, let $S$ be a surface obtained from a closed, connected, oriented surface after removing a finite number of points and open disks; we call such an $S$ a surface of finite-type with boundary. 

It was proven in e.g. \cite[Theorem 11.1]{HHS_II} that, if $S$ has no boundary, then it admits a HHG structure. In this Section we first extend this result to surfaces of finite-type with boundary (see Remark~\ref{rem:surface_with_boundary_HHS}); then we apply our main Theorems to produce two combinatorial HHG structures, one whose underlying graph is a blow-up of the curve graph (Theorem~\ref{thm:mcg_cHHS}), and one with combinatorial wedges and combinatorial containers (Theorem~\ref{thm:strong_mcg}).

\subsection{On the meaning of subsurface}
The “usual” HHG structure for a mapping class group involves \emph{open} subsurfaces, but it will be convenient here to consider a more general type of subsurfaces; here we present the two notions and compare them.

\begin{defn}[{\cite[Section 2.1.3]{BKMM}}]\label{defn:essential_sub}
A subsurface $Y\subset S$ is \emph{essential} if it is the disjoint union of some components of the complement of a collection of disjoint simple closed curves, so that no component is a pair of pants and no two annular components are isotopic.
\end{defn}

Recall that nesting of subsurfaces is defined as follows: $U$ is nested in $V$ if $U$ is contained in $V$ (up to isotopy) and no isotopy class representative of $U$ is disjoint from an isotopy class representative of $V$. (This last clause is only relevant for annuli and unions of annuli; an annulus might be isotopic to a non-essential annulus of another subsurface.)

\begin{thm}[{\cite[Theorem 11.1]{HHS_II}}]\label{thm:original_HHS_mcg}
    Let $S$ be a surface of finite-type with boundary. Its mapping class group $\mathcal{MCG}(S)$ is an HHG with the following structure:
\begin{itemize}
\item $\frakS$ is the collection of isotopy classes of essential subsurfaces.
\item For each $U\in \frakS$ the space $\fontact U$ is its curve graph. 
\item The relation $\nest$ is nesting, $\orth$ is disjointness and $\transverse$ is overlapping.
\item For each $U\in\frakS$, the projection $\pi_U:\,\mathcal{MCG}(S)\to \fontact U$ is constructed using the subsurface projection.
\item For $U,V\in \frakS$ satisfying either $U\propnest V$ or $U\transverse V$, the projection is $ \rho^U_V=\pi_V(\partial U)\subset \fontact V$, while for $V\propnest U$ the map $\rho^U_V:\,\fontact U\to 2^{\fontact V}$ is the subsurface projection.
\end{itemize}
\end{thm}

\begin{rem}\label{rem:surface_with_boundary_HHS}
    In previous literature, the HHG structure above is only considered for surfaces without boundary, but everything goes through as above for surfaces with boundary (including braid groups). There are a few ways to see this, besides inspecting \cite[Section 11]{HHS_II}. One is to regard $\mathcal{MCG}(S)$ as above as a subgroup of the mapping class group of the double of $S$ along all boundary components, by extending mapping classes to be the identity on the complement of $S$. In this setting, $\mathcal{MCG}(S)$ acts properly and coboundedly on  $F_S$ (with finitely many orbits of subsurfaces nested into $S$), giving the required structure.
\end{rem}

We now describe a different HHS structure, whose index set is made of subsurfaces which might include some of their boundary components, and then discuss how it relates to the one from Theorem~\ref{thm:original_HHS_mcg}.

\subsubsection{The new index set.} 
\begin{defn}[Block]\label{defn:block} A \emph{block} is (the isotopy class of) a subsurface of one of the following types:
\begin{enumerate}[label=(\alph*)]
    \item \label{block_annulus} a closed annulus which does not bound a disk or a single puncture;
    \item \label{block_other} a connected, non-annular subsurface of complexity at least $1$, with some (possibly none) of its boundary components included.
\end{enumerate}
The \emph{included boundary} of a block $U$ of type~\ref{block_other} is the set of curves in its topological boundary (relative to $S$) which belong to $U$. With a slight abuse of notation, we say that the included boundary of an annulus is its core curve.
\end{defn}

Given two blocks $U$ and $V$, we say $U\nest V$ if $U$ can be isotoped inside $V$, and every included boundary component of $U$ which is homotopic to a component of the topological boundary of $V$ is also contained in the included boundary of $V$. Moreover, let $\orth$ be disjointness of blocks (up to isotopy), but with the convention that:
\begin{itemize}
    \item if $\gamma$ is a curve in the topological boundary of $U$, then $\gamma\orth U$ if and only if $\gamma$ does not belong to the included boundary of $U$;
    \item if $U$ and $V$ are blocks of type~\ref{block_other}, and they share a component of the included boundary, then they are not orthogonal.
\end{itemize} This way, nesting and orthogonality are mutually exclusive, meaning that if two blocks are orthogonal then they are not $\nest$-related.

\begin{figure}[htp]
    \centering
    \includegraphics[width=\textwidth]{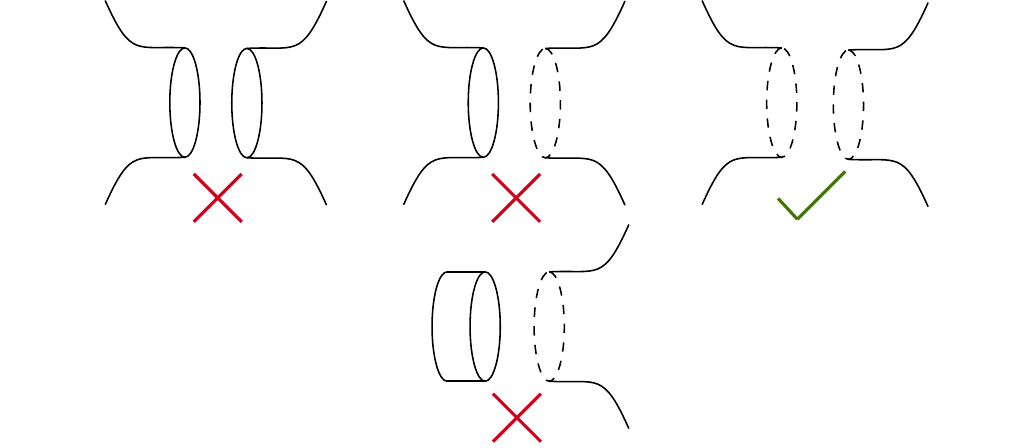}
    \caption{Representation of the forbidden and allowed collections of blocks. Dashed lines represent open boundaries (those which do not contain the boundary curve), while full lines represent included boundaries. The forbidden cases are as follows: no two blocks (which might coincide) have isotopic boundary components, unless both of them are not included (top row); no annulus is isotopic to a boundary component of another block (bottom row).}
    \label{fig:bad_blocks}
\end{figure}

\begin{defn}[Admissible collection]\label{defn:admissible_blocks}
A collection of pairwise orthogonal blocks $U_1,\ldots, U_k$ is \emph{admissible} if, whenever two blocks $U_i$ and $U_j$ (which might coincide) have two topological boundary components which are isotopic (and are distinct if $U_i=U_j$), then none of these components belongs to the included boundary of the respective block. See Figure~\ref{fig:bad_blocks} to understand the forbidden cases.
\end{defn}

Let $\frakS'$ be the set of all admissible collections of blocks, which we see as subsurfaces of $S$ (up to isotopy). Extend nesting and orthogonality to $\frakS'$, with the same conventions about boundary curves. Note that $\nest$-minimal elements of $\frakS'$ are exactly annuli of type~\ref{block_annulus}.

Define the \emph{interior} of an admissible collection as the union of the interiors of its blocks of type~\ref{block_other}. Moreover, define the \emph{included boundary} of an admissible collection as the union of the included boundaries of its blocks. Notice that if $U,V\in \frakS'$ then $U\nest V$ if and only if the interior of $U$ is nested in the interior of $V$, and the included boundary of $U$ is nested in the included boundary of $V$.

\par\medskip
\textbf{Coordinate spaces.} Now, let $V$ be a block. Define $\fontact V$ as follows:
\begin{itemize}
    \item If $V$ is a closed, essential annulus, let $\fontact V$ be its annular curve graph;
    \item If $V$ is a connected, open subsurface of complexity at least $1$, let $\fontact V$ be its curve graph;
    \item If $V$ is of type~\ref{block_other} and its included boundary is non-empty, let $\fontact V$ be the join of the curve graph of its interior and the included boundary.
\end{itemize}
The coordinate space of an admissible collection is, by definition, the join of the coordinate spaces of its blocks. 
\begin{rem}\label{rem:unbounded_CU_for_blocks}
    Notice that, if $U\in \frakS'$, then either $U$ is an annulus or a connected, open subsurface of complexity at least $1$, or $\fontact U$ is uniformly bounded since it is a join.
\end{rem}

\par\medskip
\textbf{Projections.} For every $U\in \frakS'$, let $\pi_U:\,\mathcal{MCG}(S)\to \fontact U$ be the subsurface projection. Analogously, for every $U,V\in \frakS'$ such that $U\propnest V$, let $\rho^V_U:\,\fontact V\to \fontact U$ be the subsurface projection. Moreover, if $\fontact V$ is bounded then define $\rho^U_V$ by choosing any point in $\fontact V$; otherwise $V$ must be a connected, essential subsurface, and therefore at least one of the topological boundary curves of $U$ is nested in $V$, so one can set $\rho^U_V=\partial U\cap V$. Finally, for every $U,V\in \frakS'$ such that $U\transverse V$ let $\rho^U_V$ be the subsurface projection of the topological boundary of $U$ inside $\fontact V$.

\begin{rem}[{Comparison with the index set from \cite{BKMM}}]
 To pass from our new index set $\frakS'$ to $\frakS$, one must exchange every subsurface $U$ with included boundary for the disjoint union of the interior of $U$ and open annuli corresponding to the included boundary curves. This procedure preserves nesting, orthogonality, curve graphs, and projections. Therefore, from the fact that $(\mathcal{MCG}(S), \frakS)$ is a HHG one can deduce that $(\mathcal{MCG}(S), \frakS')$ is a HHG. The main point here is that the two structures have the same unbounded coordinate spaces, namely the curve graphs of connected open subsurfaces that are not pairs of pants, which takes care of most of the axioms (e.g. uniqueness, which constitutes the majority of the work in \cite[Section 11]{HHS_II}).
 \end{rem}

\subsection{First CHHS structure}
Our next goal is to present the first CHHS structures for $\mathcal{MCG}(S)$, whose underlying graph $X$ is a blow-up of the coordinate space $\fontact S$. This will answer a question from \cite[Subsection 1.6]{BHMS}. As mentioned in Remark~\ref{rem:clean_markings_for_mcg}, the graph $\mathcal{W}$ of this structure will be very similar to the graph of complete clean markings from \cite{MasurMinsky2}.
\begin{thm}\label{thm:mcg_cHHS}
Let $S$ be a surface of finite-type with boundary. There exists a combinatorial HHG structure $(X,\mathcal{W})$ for $\mathcal{MCG}(S)$, where $X$ is the blow-up of the curve graph of $S$, obtained by replacing every curve with the cone over its annular curve graph.
\end{thm}

\begin{proof}
It is enough to show that $(\mathcal{MCG}(S), \frakS')$, with the coordinate spaces and the projections defined as above, satisfies the hypotheses of Theorem~\ref{thm:comb_hhg}.

\begin{itemize}
\item \textbf{Wedges}: Let $U,V\in\frakS'$, and let $W$ be the intersection of their interiors, which is a disjoint union of open pairs of pants and open subsurfaces of complexity at least $1$. First, for every curve $\gamma$ in the topological boundary of $W$ which is nested in both $U$ and $V$, glue to $W$ the closed annulus with core $\gamma$, so that one (resp. both) of the boundary curves of the annulus is identified with one (resp. two) curves in the topological boundary of $W$. Then replace every component of $W$ which is a pair of pants with its included boundary (in particular, one has to remove open pairs of pants). Let $W'$ the surface obtained after this procedure, and let $W''$ be the disjoint union of $W'$ and every annulus which is a block of both $U$ and $V$. Now, by construction $W''$ is in $\frakS'$ (since the gluing procedure prevents the forbidden cases from Figure~\ref{fig:bad_blocks} from appearing), and it is nested in both $U$ and $V$. Moreover $W''$ is the wedge of $U$ and $V$, since $T\in \frakS'$ is nested in both $U$ and $V$ if and only if its interior is nested in the interior of $W$ (which is the intersection of the interiors, without the pants components), and its included boundary is nested in the intersection of the included boundaries.
\par\medskip
\item \textbf{Clean containers}: whenever $U\propnest V\in\frakS'$ and there exists $W\in \frakS'$ which is nested in $V$ and orthogonal to $U$, consider the subsurface obtained from $V-U$ after replacing every pair of pants with its included boundary (in particular, one has to remove open pairs of pants). Let $Y$ be the resulting subsurface, which is non-empty since it contains $W$. First we prove that the connected components of $Y$ are blocks. Indeed, the connected surfaces which are not blocks are annuli with at most one boundary included, or pairs of pants with some boundary components included. Now, one of the connected components of $Y$ is an open annulus if and only if two blocks of $U$ share a common curve in their included boundary, and this would mean that these blocks are not orthogonal by our convention. Moreover, one of the connected components of $Y$ is an annulus with exactly one boundary curve included if and only if two blocks of $U$ fall in one of the forbidden cases from Figure~\ref{fig:bad_blocks}. Finally, we manually replaced every pair of pants with its included boundary. Hence $Y$ is a disjoint union of blocks. Similar arguments show that the blocks of $Y$ cannot share curves in their included boundaries (otherwise $U$ would contain an open annulus), and cannot fall in one of the forbidden cases from Figure~\ref{fig:bad_blocks} (otherwise $U$ would contain an annulus with exactly one boundary component included). This shows that $Y$ is an admissible collection, i.e. an element of $\frakS'$, and by construction it is orthogonal to $U$ and nested in $V$. Being the maximal subsurface with these properties, $Y$ is also the clean container for $U$ inside $V$.
\par\medskip
\item \textbf{Orthogonals for non-split domains property~\eqref{property:orthogonal_for_non_split}}: Let $U\propnest V$ be two subsurfaces. If $U$ has a connected component which is an annulus then $U$ is split, and such annulus is one of its Samaritans. Otherwise $U$ is the disjoint union of finitely many subsurfaces of complexity at least $1$, possibly with boundary. If there exists a connected component of $V$ whose intersection with $U$ is trivial, then such component is orthogonal to $U$. Otherwise there must be a curve $\gamma$ in the boundary of $U$ relative to $V$, which must therefore be essential in $V$. Then the associated annular domain $A_\gamma$ is nested in $V$; moreover, $A_\gamma$ is either nested in $U$ or disjoint from $U$, depending on whether the boundary curve is included in $U$ or not. In the former case, $U$ is a split domain, and $A_\gamma$ is one of its Samaritans; in the latter, $A_\gamma$ is orthogonal to $U$.
\par\medskip
\item \textbf{Dense product regions~\eqref{property:dpr}}: let $U\in \frakS'$ be a non-minimal domain. Then either $U$ is a connected, open subsurface of complexity at least $1$, and its curve graph is covered by the projections of the annuli it contains, or $\fontact U$ is uniformly bounded, as pointed out in Remark~\ref{rem:unbounded_CU_for_blocks}.
\par\medskip
\item \textbf{Cofinite action}: The action of the mapping class group on $\frakS$ has finitely many orbits of pairwise disjoint subsurfaces. This is a consequence of the change of coordinate principle (see e.g. \cite[Section 1.3.3]{FarbMargalit}).
\end{itemize}
This proves the Theorem.
\end{proof}

\subsection{Second CHHS structure, with better properties}
Notice that the HHS structure from Theorem~\ref{thm:mcg_cHHS} does not satisfy the weak orthogonal property~\eqref{property:weak_orth}. To see this, let $V$ be an open, connected subsurface of complexity $2$ and let $U$ be the union of two disjoint, essential annuli inside $V$. Then $U$ is non-minimal, but $V-U$ is a disjoint union of open pairs of pants and therefore cannot contain any element of $\frakS'$.

However, we can find a larger index set to ensure even the strong orthogonal property~\eqref{property:strong_orth}:

\begin{thm}\label{thm:strong_mcg}
    Let $S$ be a surface of finite-type with boundary. There exists a combinatorial HHG structure for $\mathcal{MCG}(S)$ with simplicial wedges (Definition~\ref{defn:simp_wedge}) and simplicial containers (Definition~\ref{defn:simp_cont}).
\end{thm}

\begin{proof}
The proof is very similar to that of Theorem~\ref{thm:mcg_cHHS}. First, we weaken Definition~\ref{defn:block}, by allowing a block to be also a pair of pants, with some (possibly none) of its boundary components included. Define nesting and orthogonality between blocks as before, with the same conventions about the included boundary components. This allows one to define the collection $\frakS''$ of admissible blocks, as in Definition~\ref{defn:admissible_blocks}. Notice that now the $\nest$-minimal elements of $\frakS''$ are all closed annuli of type~\ref{block_annulus}, and all open pants. 

Define the interior and the included boundary of an admissible collection as before, but now the interior also includes the interior of the pants components. Again, notice that if $U,V\in \frakS'$ then $U\nest V$ if and only if the interior of $U$ is nested in the interior of $V$, and the included boundary of $U$ is nested in the included boundary of $V$.

Define the coordinate spaces as before, and set the coordinate space of an open pair of pants to be a point. Again, the only elements of $ \frakS'$ with unbounded coordinate spaces are annuli and connected, open subsurfaces of complexity at least $1$, because all other coordinate spaces are either points or joins.

Finally, define the projections as above, using subsurface projections. With the same techniques of \cite[Theorem 11.1]{HHS_II}, one can then show that $(\mathcal{MCG}(S), \frakS'')$ is a HHG.

Next, we observe that $(\mathcal{MCG}(S), \frakS'')$ has the following properties. The proofs are very similar to those which appear in Theorem~\ref{thm:mcg_cHHS} (and even easier since we do not have to remove pairs of pants), but we put them here for clarity:

\begin{itemize}
\item \textbf{Wedges}: Let $U,V\in\frakS''$, and let $W$ be the intersection of their interiors, which is a disjoint union of open pants and open subsurfaces of complexity at least $1$. Then, for every curve $\gamma$ in the topological boundary of $W$ which is nested in both $U$ and $V$, glue to $W$ the closed annulus with core $\gamma$, so that one (resp. both) of the boundary curves of the annulus is identified with one (resp. two) curves in the topological boundary of $W$. Let $W'$ the surface obtained after the gluing, and let $W''$ be the disjoint union of $W'$ and every annulus which is a block of both $U$ and $V$. Now, by construction $W''$ is in $\frakS''$ (since the gluing procedure prevents the forbidden cases from Figure~\ref{fig:bad_blocks} from appearing), and it is nested in both $U$ and $V$. Moreover, the interior of $W''$ is the intersection of the interiors, and the included boundary of $W''$ is the intersection of the included boundaries. This shows that $W''$ is the wedge of $U$ and $V$, since $T\in \frakS'$ is nested in both $U$ and $V$ if and only if its interior is nested in the intersection of the interiors, and its included boundary is nested in the intersection of the included boundaries.
\par\medskip
\item \textbf{Clean containers and the strong orthogonal property~\eqref{property:strong_orth}}: whenever $U\propnest V\in\frakS''$, consider the subsurface $V-U$. First notice that the connected components of $V-U$ are blocks. Indeed, the only connected subsurfaces which are not blocks are annuli with at most one of the two boundary curves included. Now, one of the connected components of $V-U$ is an open annulus if and only if two blocks of $U$ share a common curve in their included boundary, and this would mean that these blocks are not orthogonal by our convention. Moreover, one of the connected components of $V-U$ is an annulus with exactly one boundary curve included if and only if two blocks of $U$ fall in one of the forbidden cases from Figure~\ref{fig:bad_blocks}. Hence $V-U$ is a disjoint union of blocks. The same argument with $U$ and $V-U$ swapped shows that the blocks of $V-U$ cannot share curves in their included boundaries, and cannot fall in one of the forbidden cases from Figure~\ref{fig:bad_blocks}. This shows that $V-U$ is an admissible collection, i.e. an element of $\frakS'$, and by construction it is orthogonal to $U$ and nested in $V$. Being the maximal subsurface with these properties, $V-U$ is also the clean container for $U$ inside $V$.
\par\medskip
\item \textbf{Dense product regions}: let $U\in \frakS''$ be a non-minimal domain. Then either $U$ is an essential open subsurface, and its curve graph is covered by the projections of the annuli it contains, or $\fontact U$ is uniformly bounded, as pointed out above.
\par\medskip
\item \textbf{Cofinite action}: the $\mathcal{MCG}(S)$-action on $\frakS''$ has finitely many orbits of tuples of pairwise orthogonal domains, again as a consequence of the change of coordinates principle.
\end{itemize}
Now the hypotheses of Theorem~\ref{thm:strong_iff} are satisfied, and we can find a combinatorial HHS $(X,\mathcal{W})$ with simplicial wedges and simplicial containers. Since $(\mathcal{MCG}(S), \frakS')$ already has dense product region, the pair $(X,\mathcal{W})$ is exactly the one constructed in Section~\ref{sec:construction} from $(\mathcal{MCG}(S), \frakS')$, as pointed out in Remark~\ref{rem:already_has_DPR}. Hence, by Theorem~\ref{thm:comb_hhg} we have that $(X,\mathcal{W})$ inherit an action of $\mathcal{MCG}(S)$ that makes the latter into a combinatorial HHG, and we are done.
\end{proof}

\section{Why orthogonals for non-split domains?}\label{sec:CHHS_for_orth_sets}
The main goal in this section is to use a simple example of an HHS --- a CAT(0) cube complex with a factor system --- to illustrate the necessity of the orthogonals for non-split domains hypothesis~\eqref{property:orthogonal_for_non_split}.  

Factor systems yield examples of HHS structures where the index set is an \emph{orthogonal set} in the sense of Definition~\ref{defn:orth_set}. It would be illuminating to find conditions on an orthogonal set allowing one to modify the HHS/HHG structure so that some version of the orthogonality properties (\eqref{property:orthogonal_for_non_split},~\eqref{property:weak_orth} or~\eqref{property:strong_orth}) hold. We speculate on this below, and in particular on intriguing relations with problems in lattice theory, namely embedding a complete ortho-lattice inside an orthomodular one, see Remark~\ref{rem:lattice}.

In this section, we use notation from \cite{Rcubes}; see also \cite{HagenSusse} and \cite[Section 8]{HHS_I}.

\subsection{Background on factor systems for CAT(0) cube complexes}\label{subsec:CCC}
For the rest of the Section, let $\cuco Z$ be a CAT(0)-cube complex. 
\begin{defn}[Hyperplane, carrier, combinatorial hyperplane]\label{defn:hyperplane}
A \emph{midcube} in the unit cube $c=[-\frac12,\frac12]^n$ is a subspace obtained by restricting exactly one coordinate to $0$. A \emph{hyperplane} in $\cuco Z$ is a connected subspace $H$ with the property that, for all cubes $c$ of $\cuco Z$, either $H\cap c=\emptyset$ or $H\cap c$ consists of a single midcube of $c$. The \emph{carrier} $\neb(H)$ of the hyperplane $H$ is the union of all closed cubes $c$ of $\cuco Z$ with $H\cap c\neq\emptyset$. The inclusion $H\to\cuco Z$ extends to a combinatorial embedding $H\times[-\frac12,\frac12]\stackrel{\cong}{\longrightarrow}\neb(H)\hookrightarrow\cuco X$ identifying $H\times\{0\}$ with $H$.  Now, $H$ is isomorphic to a CAT(0) cube complex whose cubes are the midcubes of the cubes in $\neb(H)$.  The subcomplexes $H^\pm$ of $\neb(H)$ which are the images of $H\times\{\pm\frac12\}$ under the above map are isomorphic as cube complexes to $H$, and are \emph{combinatorial hyperplanes} in $\cuco Z$.  Thus each hyperplane of $\cuco Z$ is associated to two combinatorial hyperplanes in $\neb(H)$.
\end{defn}

\begin{defn}[Gate maps]
    For any convex subcomplex $\cuco Y\subseteq \cuco Z$ there is a \emph{gate map} $\gate_{\cuco Y}:\cuco Z\to\cuco Y$ such that, for any other convex subcomplex $\cuco Y'\subseteq \cuco Z$, the hyperplanes crossing $\gate_{\cuco Y}(\cuco Y')$ are precisely the hyperplanes which cross both $\cuco Y$ and $\cuco Y'$. 
\end{defn}
 Gate maps are fundamental in the study of cube complexes and median spaces; see, for instance, \cite[Section 2]{HHS_I} for additional background.

\begin{defn}[Parallelism]\label{defn:parallel_subcomplexes}
The convex subcomplexes $F$ and $F'$ are \emph{parallel}, written $F\parallel F'$, if for each hyperplane $H$ of $\cuco Z$, we have $H\cap F\neq\emptyset$ if and only if $H\cap F'\neq\emptyset$.
\end{defn}

\begin{defn}[Orthogonality, orthogonal complement]\label{defn:orthogonal_subcomplexes}
The convex subcomplexes $F$ and $F'$ are \emph{orthogonal}, written $F\orth F'$, if the inclusions $F\to \cuco Z$ and $F'\to \cuco Z$ extend to a convex embedding $F\times F'\to \cuco Z$.

Given a convex subcomplex $F$, let $P_F$ be the smallest subcomplex containing the union of all subcomplexes in the parallelism class of $F$.  By e.g. \cite[Lem. 1.7]{HagenSusse} there is a cubical isomorphism $P_F\to F\times F^\orth$, where $F^\orth$ is a CAT(0) cube complex which we call the \emph{abstract orthogonal complement of $F$}.  For any $f\in F^{(0)}$, the inclusion $P_F\to \cuco Z$ induces an isometric embedding $\{f\}\times F^\orth\to \cuco Z$ whose image is a convex subcomplex that we call the \emph{orthogonal complement of $F$ at $f$} and denote $\{f\}\times F^\orth$.  Observe that $\{f\}\times F^\orth$ and $\{f'\}\times F^\orth$ are parallel for all $f,f'\in F^{(0)}$ (see \cite[Lem. 1.11]{HagenSusse}).  When the base point is unimportant, we sometimes abuse notation and write $F^\orth$ to refer to one of these parallel copies.  
\end{defn}

\begin{lemma}\label{lem:container-supporting}
Let $F\subset \cuco Z$ be a convex subcomplex, let $f\in F^{(0)}$, and write $F^\orth=\{f\}\times F^\orth$.  Suppose that $F'$ is a convex subcomplex such that $F\orth F'$.  Then $F'$ is parallel to a subcomplex of $F^\orth$.  Conversely, if $F'$ is a convex subcomplex of $F^\orth$, then there exists $F''$ parallel to $F'$ with $F''\orth F$.
\end{lemma}

\begin{proof}
This follows easily from \cite[Lem. 1.11]{HagenSusse}; see for instance the proof of \cite[Theorem C.]{HagenSusse}. 
\end{proof}

\begin{defn}[Candidate factor system]
A \emph{candidate factor system} $\mathfrak h$ is a collection of non-trivial convex subcomplexes of $\cuco Z$ (where “non-trivial” means that we exclude singletons) that satisfy the following properties:
\begin{enumerate}
 \item \label{item:hyper_1} $\cuco Z\in\mathfrak h$, and for all combinatorial hyperplanes $H$ of $\cuco Z$, we have $H\in\mathfrak h$;
 \item\label{item:hyper_2} if $F,F'\in\mathfrak h$ then $\gate_F(F')\in\mathfrak h$;
 \item\label{item:hyper_3} if $F\in\mathfrak h$ and $F'$ is parallel to $F$, then $F'\in\mathfrak h$.
\end{enumerate}
\end{defn}

\begin{defn}[Hyperclosure]\label{defn:hypclosure}
The \emph{hyperclosure} $\bar{\mathfrak h}$ of $\cuco Z$ is the intersection of all candidate factor systems, and therefore the unique candidate factor system which is minimal by inclusion. In other words
$$\bar{\mathfrak h}=\left(\bigcup_{i=1}^\infty \mathfrak h_i\right)-\{\mbox{singletons}\},$$
where:
\begin{itemize}
    \item $\mathfrak h_1$ is the collection of all subcomplexes which are parallel to combinatorial hyperplanes, together with the whole space $\cuco Z$;
    \item For all $i\ge 1$, $\mathfrak h_{i+1}=\left\{F\,\mid \,F\parallel \mathfrak g_{F_1}(F_2),\, F_1, F_2\in \mathfrak h_i\right\}$.
\end{itemize}
\end{defn}

\begin{lemma}[{Characterisation of $\bar{\mathfrak h}$, \cite[Theorem 3.3]{HagenSusse}}]\label{lem:F=Corth}
    Let $\cuco Z$ be a locally finite CAT(0) cube complex, and let $\bar{\mathfrak h}$ be its hyperclosure. Then a convex subcomplex $F$ belongs to $\bar{\mathfrak h}$ if and only if there exists a compact, convex subcomplex $C$ such that $F=C^\orth$.
\end{lemma}

\begin{defn}[Weak factor system]
    Let $\mathfrak h$ be a candidate factor system, and let $\mathfrak h_{_{/\sim}}$ the set of parallelism classes of subcomplexes in $\mathfrak h$. If there exists $N\in \mathbb{N}$ such that $N$ bounds the length of chains in the partial order of $\mathfrak h$ given by inclusion, then $\mathfrak h_{_{/\sim}}$ is a \emph{weak factor system}. 
    
    We denote the parallelism class of the subcomplex $F\in \mathfrak h$ by $[F]\in \mathfrak h_{/\sim}$. Two elements $[F],[F']\in \mathfrak h_{/\sim}$ are \emph{nested} (resp. \emph{orthogonal}), and we write $[F]\nest[F']$ (resp. $[F]\orth[F']$) if there exists two representatives $F,F'$ such that $F\subseteq F'$ (resp. $F\orth F'$).
\end{defn}

The class of CAT(0) cube complexes which admit weak factor systems is quite large. For example, virtually special groups, in the sense of Haglund-Wise \cite{HaglundWise}, act geometrically on CAT(0) cube complexes with weak factor systems by \cite[Proposition B]{HHS_I}.  The more general class of CAT(0) cube complexes with geometric group actions and weak factor systems is characterised in \cite{HagenSusse}, and includes some notable non-special examples, like irreducible lattices in products of trees, and certain amalgams of these \cite{Hagen:non-colourable}. There are other amalgams of such lattices that provide the first examples of proper cocompact CAT(0) cube complexes \emph{not} admitting any weak factor system \cite{Shepherd:counterexamples}.

Following \cite{HHS_I}, one can endow a CAT(0) cube complex with an HHS structure when the hyperclosure gives a weak factor system; this is analysed in more detail in \cite[Section 20]{Rcubes}.  We summarise this here in order to connect the existing results more explicitly to the hypotheses of Theorem~\ref{thm:main}.

The first lemma is needed to verify the clean containers property.  It is proved in \cite[Prop. 5.1]{HagenSusse} as part of a more complicated statement whose other parts rely on the presence of a cocompact group action; see \cite[Rem. 20.7]{Rcubes}.  So, for the avoidance of doubt, we extract the exact statement (with the same proof as in \cite{HagenSusse}) here:

\begin{lemma}\label{lem:hyperclosure-closed-under-complements}
Let $\cuco Z$ be a CAT(0) cube complex and let  $\bar{\mathfrak h}$ be its hyperclosure. If $\bar{\mathfrak h}_{/\sim}$ is a weak factor system then, for every $F\in\bar{\mathfrak h}-\{\cuco Z\}$ and every $x\in F^{(0)}$, the subcomplex $F^\orth\coloneq\{x\}\times F^\orth$ belongs to $\bar{\mathfrak h}$.
\end{lemma}

\begin{proof}
We will use \cite[Lem. 5.2]{HagenSusse} after some preliminary setup.  Let $\{H_i\}_{i\in I}$ be the set of hyperplanes that are dual to edges of $F$ incident to $x$.  For $i\in I$, let $H_i^+$ be the combinatorial hyperplane in which the carrier of $N(H_i)$ intersects the $H_i$--halfspace of $\cuco Z$ containing $x$; in particular, $x$ lies in $ H_i^+$.  Let $Y=\bigcap_{i\in I}H_i^+$.  Observe that for any finite $I'\subset I$, we have that $Y(I'):=\bigcap_{i\in I'}H_i^+\in \bar{\mathfrak h}$, unless $Y(I')$ is a single point (recall that we do not allow single points in the hyperclosure).  Assume the former.

By the assumption that the set of parallelism classes represented in the hyperclosure is a weak factor system, together with the observation that $[Y(I')]\nest [Y(I'')]$ when $I''\subset I'$, we see that there exists a finite subset $I_0\subset I$ such that $Y=Y(I_0)$.  Hence $Y\in \bar{\mathfrak h}$, or $Y$ consists of a single vertex.

Now let $\mathcal S$ be the set of all combinatorial hyperplanes $H^\pm$ such that the associated hyperplane $H$ (i.e. the hyperplane $H$ such that the usual identification of $N(H)$ with $H\times[-\frac12,\frac12]$ identifies $H\times\{\pm\frac12\}$ with $H^\pm$) crosses $F$.  By \cite[Lem. 5.2]{HagenSusse}, $F^\orth = \bigcap_{H^\pm\in\mathcal S}\gate_Y(H^\pm)$.  We now argue as before. First, note that if $\mathcal S'\subset\mathcal S$ is finite, then $A(\mathcal S')=\bigcap_{H^\pm\in\mathcal S'}\gate_Y(H^\pm)$ belongs to $\bar{\mathfrak h}$, or it is a single point.

As before, if $\mathcal S'\subset\mathcal S''$ are finite subsets of $\mathcal S$, we have $A(\mathcal S'')\subset A(\mathcal S')$, so our assumption that $\bar{\mathfrak h}$ gives a weak factor system again implies that there is a finite $\mathcal S'$ such that $F^\orth = A(\mathcal S')$.  Hence, either $F^\orth\in \overline{\mathfrak h}$, or $F^\orth$ is a single point.

To complete the proof, we rule out the latter possibility as follows.  By hypothesis, $F\neq \cuco Z$.  By Lemma~\ref{lem:F=Corth}, there exists a (compact) convex subcomplex $C\subset\cuco Z$ with $C^\orth=F$, since $F\in\bar{\mathfrak h}$.  This means that $F\orth C$, so Lemma~\ref{lem:container-supporting} implies that, up to parallelism, $C\subset F^\orth$, so it remains to show that $C$ is non-trivial.  But if $C$ is trivial, then $C\orth\cuco Z$, so by Lemma~\ref{lem:container-supporting}, $\cuco Z\nest C^\orth=F$, contradicting that $F$ is a proper subcomplex.
\end{proof}

\begin{rem}[Why are singletons excluded?]\label{rem:excluded-singletons}
Let $\bar{\mathfrak h}_{/\sim}$ denote the set of parallelism classes in the hyperclosure $\bar{\mathfrak h}$.  In general, one gets from a weak factor system to an HHS structure using Theorem~\ref{thm:HHS_from_WFS} below.  The additional property that the index set is an orthogonal set can be arranged by using the weak factor system provided by the hyperclosure, when it exists --- see Theorem~\ref{thm:hyperclosure_is_orthogonal}.

Recall that we have defined $\bar{\mathfrak h}$ so as to exclude subcomplexes consisting of a single vertex.  Note that any two such subcomplexes $\{x\},\{y\}$ are parallel.  Moreover, $[\{x\}]\nest [F]$, for any $F\in\bar{\mathfrak h}$.  On the other hand, we also have $[\{x\}]\orth[F]$.  This would not be allowed in an HHS structure, so we avoid the issue by excluding $\{x\}$ from the hyperclosure.  For non-trivial subcomplexes, this problem does not occur, by, for instance, \cite[Prop. 2.5]{CapraceSageev} and the fact that hyperplanes do not cross themselves.

In \cite{HHS_I}, something more radical is allowed: one can exclude all subcomplexes below some fixed diameter, and still get an HHS structure.  But this may cease to satisfy the conclusion of Lemma~\ref{lem:hyperclosure-closed-under-complements}.  This is not a problem from the point of view of clean containers, but it can break “orthogonality determines nesting” by creating non-$\nest$--maximal $[F]$ that are not orthogonal to anything in the index set.
\end{rem}

\begin{thm}[{see e.g. \cite[Proposition 20.4]{Rcubes}}]\label{thm:HHS_from_WFS}
    Let $\cuco Z$ be a CAT(0) cube complex with a weak factor system $\mathfrak h_{/\sim}$. Then $(\cuco Z, \mathfrak h_{/\sim})$ is a hierarchically hyperbolic space with wedges, where the coordinate spaces $\fontact [F]$, $[F]\in \mathfrak h_{/\sim}$ and the projections $\pi_{[F]}:\,\cuco Z\to \fontact [F]$ are as in \cite[Remark 13.2]{HHS_I}.
\end{thm}

The next theorem refines the previous one by strengthening the conclusions about the index set.  It should be compared to the very similar \cite[Prop. 20.6]{Rcubes}, where, however, it is not explicit that the HHS structure comes from the hyperclosure.  Note, also, that the clean containers property of the HHS structure from Theorem~\ref{thm:HHS_from_WFS} is observed (in the presence of a group action) in \cite{HagenSusse}.

\begin{thm}\label{thm:hyperclosure_is_orthogonal}
    Let $\cuco Z$ be a CAT(0) cube complex which admits a weak factor system, and let $\bar{\mathfrak h}$ be the hyperclosure of $\cuco Z$. Then $\bar{\mathfrak h}_{/\sim}$ is a weak factor system and an orthogonal set. In particular, $(\cuco Z, \bar{\mathfrak h}_{/\sim})$ is a HHS with wedges, clean containers, and where orthogonality implies nesting.
\end{thm}

\begin{proof}
    In the proof of \cite[Proposition 20.6]{Rcubes} it is shown that, if $\cuco Z$ admits a weak factor system, then there exists a weak factor system $\mathfrak h'_{/\sim}$ which consists of all equivalence classes of subcomplexes of the form
    $$\mathfrak g_{H_1}(\ldots(\mathfrak g_{H_{n-1}}(H_n))\ldots),$$
    where $H_1,\ldots, H_n$ are combinatorial hyperplanes and $n \ge 0$. By Definition~\ref{defn:hypclosure} of the hyperclosure $\bar{\mathfrak h}$, we see that $\mathfrak h'\subseteq \bar{\mathfrak h}$, and they must coincide since the hyperclosure is the minimal candidate factor system. Hence, the quotient of the hyperclosure by parallelism is a weak factor system, and therefore $(\cuco Z, \bar{\mathfrak h}_{/\sim})$ is a HHS with wedges by Theorem~\ref{thm:HHS_from_WFS}.
    
    Now, by Lemma~\ref{lem:hyperclosure-closed-under-complements}, if $\bar{\mathfrak h}_{/\sim}$ is a weak factor system then, for every $F\in \bar{\mathfrak h}$, we have that $F^\orth \in \bar{\mathfrak h}$. In particular, this shows that $\bar{\mathfrak h}_{/\sim}$ has clean containers, since if $[F],[F'],[C]\in \bar{\mathfrak h}_{/\sim}$ are such that $[F],[F']\nest [C]$ and $[F]\orth [F']$, then the clean container for $[F]$ inside $[C]$ is $[C]\wedge [F]^\orth$, where $[F]^\orth=[F^\orth]$, as provided by Lemma~\ref{lem:container-supporting}.
    
    Moreover, by \cite[Corollary 3.4]{HagenSusse} we have that $F^{\orth\orth}=F$ whenever $F\in \bar{\mathfrak h}$, and by Lemma~\ref{lem:orth_implies_nest_equiv} this is equivalent to the fact that orthogonality determines nesting in $\bar{\mathfrak h}_{/\sim}$.
\end{proof}

\begin{rem}\label{rem:easy-DPR-cubes}
 The HHS structure $(\cuco Z,\bar{\mathfrak h}_{/\sim})$ need not have the DPR property.  When it does, this can be verified as follows.  For each $[F]\in\bar{\mathfrak h}_{/\sim}$, the HHS structure from \cite{HHS_I} has as the coordinate space $\fontact[F]$ the \emph{factored contact graph}, which can be viewed as a copy of $F$ with additional edges added to cone off subcomplexes of the form $\gate_F(F')$, where $\gate_F:\cuco Z\to F$ is the gate map and $F'\in\bar{\mathfrak h}$.  Now, the hyperplanes of $F$ have the form $H\cap F$, where $H$ is a hyperplane of $\cuco Z$ crossing $F$.  This shows that subcomplexes of $F$ of the form $\gate_F(F')$ cover $F$.  This implies DPR provided sufficiently many of those subcomplexes to coarsely cover $F$ are actually in $\bar{\mathfrak h}$, i.e. they are not singletons.  So, for example, $(\cuco Z,\bar{\mathfrak h}_{/\sim})$ has the DPR property provided there exists a constant $K$ such that for all non $\nest$--minimal $F\in\bar{\mathfrak h}$ and all $x\in F$, there exists $F'\in\bar{\mathfrak h}$ such that $\dist_{\cuco Z}(x,\gate_F(F'))\leq K$ and $|\gate_F(F')|>1$.
\end{rem}

\subsection{The counterexample}\label{subsec:counterex_CCC}
In this Subsection we present a HHS $(\cuco Z, \frakS)$ which satisfies all hypotheses of the main Theorem~\ref{thm:main} except the orthogonals for non-split domains property~\eqref{property:orthogonal_for_non_split}, and we prove that the graph $X$ from Definition~\ref{defn:blow-up}, constructed using the coordinate spaces in the HHS structure, cannot be the underlying graph of any combinatorial HHS structure for $\cuco Z$.

Remarkably, $\cuco Z$ is a CAT(0) cube complex admitting a weak factor system, and $\frakS=\mathfrak h_{/\sim}$ is the quotient of its hyperclosure by parallelism, which is also an orthogonal set by Theorem~\ref{thm:hyperclosure_is_orthogonal}. This shows that the orthogonals for non-split domains property, which is the most obscure among the properties of our main Theorem, is essential for the construction from Section~\ref{sec:construction} to work, even if we assume that orthogonality determines nesting.

\begin{figure}[htp]
    \centering
    \includegraphics[width=\textwidth]{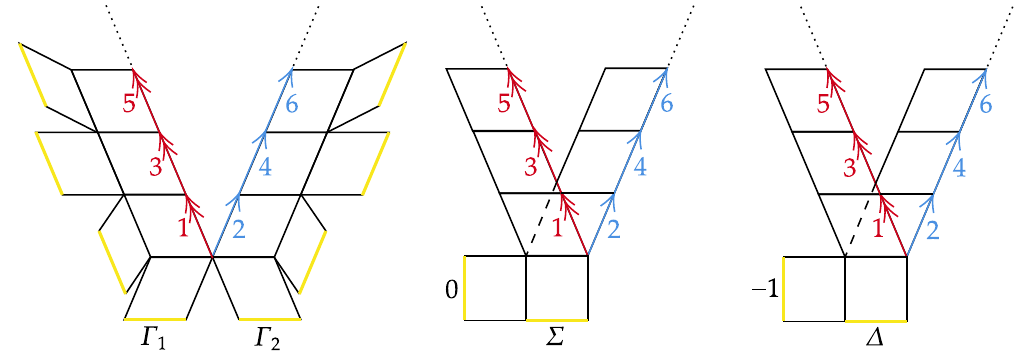}
    \caption{The space $\cuco Z$ is obtained by gluing the three pieces along the red and blue arrows. The vertical edges are labelled by numbers, while the horizontal edges are labelled by Greek capital letters. Notice that $\Gamma_1$ is orthogonal to all edges with positive odd labels, which form an infinite ray, and similarly $\Gamma_2$ is orthogonal to all edges with positive even labels. The yellow edges are representatives for the parallelism classes of edges which are in the hyperclosure.}
    \label{fig:counterexample_image}
\end{figure}

Consider the two-dimensional cube complex $\cuco Z$ obtained by gluing the three infinite complexes in Figure~\ref{fig:counterexample_image} along the red and blue arrows, respecting the numerical labels.

\begin{lemma}\label{lem:counterexample-cat0}
The complex $\cuco Z$ is CAT(0).
\end{lemma}

\begin{proof}
This follows from various versions of the principle that gluing simply-connected, non-positively-curved spaces along convex subspaces using isometries yields a simply-connected non-positively-curved result.  More precisely:

\textbf{Combinatorial version:}  A cube complex is CAT(0) if and only if its $1$--skeleton is a median graph \cite[Theorem 6.1]{Chepoi:median}.  Products of median graphs are median, and a graph decomposing as the union of two median graphs intersecting along a subgraph that is convex in each piece is median (see e.g. \cite{Isbell,Chepoi:median}).

\textbf{CAT(0) version:}  A complete geodesic metric space is CAT(0) if it is the union of two CAT(0) spaces whose intersection is convex in each piece, and products of CAT(0) spaces are CAT(0) \cite[Theorem 11.1, Exercise 1.16.(2)]{BridsonHaefliger}.

The part of $\cuco Z$ at left is obtained from two copies of $[0,2]\times[0,\infty)$, glued along a point, by gluing squares along a collection of edges.  A similar observation applies to the two pieces at right.  Apply the gluing principle once more, using that the red/blue line is convex in each piece where it appears.
\end{proof}

  % he two spaces at the right of Figure~\ref{fig:counterexample_image} are products of trees, and hence CAT(0).  The space at the left is obtained as follows: two copies of $[0,2]\times[0,\infty)$ are glued at a point, and then squares are glued to the resulting complex along edges.  This is CAT(0) because at each step, we are gluing two CAT(0) cube complexes by identifying a convex subcomplex of one with a convex subcomplex of the other using a combinatorial isometry; this preserves the CAT(00 property by \begin{com}ref\end{com}.  The red/blue lines in each space are convex, so another application of \begin{com}same ref\end{com} shows that $Z$ is CAT(0).

  Next, let us check quickly that $\cuco Z$ admits a factor system:

  \begin{lemma}\label{lem:counterexample-factor-system-sleazy}
  Let $\bar{\mathfrak h}$ be the hyperclosure in $\cuco Z$.  Then $\bar{\mathfrak h}_{/\sim}$ is a weak factor system.  More strongly, $\bar{\mathfrak h}$ is a factor system in the sense of \cite{HHS_I}.

  Finally, the HHS structure $(\cuco Z,\bar{\mathfrak h}_{/\sim})$ from Theorem~\ref{thm:hyperclosure_is_orthogonal} has the DPR property.
  \end{lemma}

  \begin{proof}
  Suppose that $[F_1],\ldots,[F_n]\in\bar{\mathfrak h}_{/\sim}$ satisfy $[F_i]\propnest [F_{i+1}]$ for all $i$, with each $F_i\subsetneq \cuco Z$. For each $i$ consider the parallelism class $[F_i^\orth]$ of its abstract orthogonal, which, by Lemma~\ref{lem:container-supporting}, coincides with the parallelism class $[F_{i}]^\orth$ of the maximal subcomplex orthogonal to $F_i$. Since $[F_i]\propnest [F_{i+1}]$ we have that $[F_i]^\orth\nest [F_{i-1}]^\orth$; moreover the nesting is proper, since by \cite[Corollary 3.4]{HagenSusse} we have that $[F_i]^{\orth\orth}=[F_i]$.

  Now, if $\bar{\mathfrak h}_{/\sim}$ had arbitrarily large $\nest$--chains, then for any $R\geq 0$, we could choose $n$ as above so that for some $m\leq n$, the subcomplexes $F_m$ and $F_m^\orth$ both have at least $R$ vertices (take $n$ much larger than $R$, and consider $m=\lfloor n/2\rfloor$).

  Next observe that there is a uniform bound on the degrees of vertices in $\cuco Z$.  Hence, for any $R_1\geq 0$, we can choose $R$ and thus the $F_i$ so that $F_m$ and $F_m^\orth$, chosen as above for the given $R$, have diameter more than $R_1$. 

  Thus $\cuco Z$ contains a convex subcomplex isometric to $F_m\times F_m^\orth$, with each factor having diameter at least $R_1$.  Now observe that the inclusion into $\cuco Z$ of the union of the red and blue rays is a quasi-isometry $\mathbb R\to \cuco Z$, so $\cuco Z$ is hyperbolic.  By taking $R_1$ sufficiently large in terms of the hyperbolicity constant, we contradict, say, \cite[Theorem 7.6]{Hagen:contact} or \cite{CDEHV}.

  Now, $\cuco Z$ is uniformly locally finite (i.e. the number of $0$--cubes in a ball is bounded in terms of the radius of the ball only), and $\bar{\mathfrak h}$ is closed under taking intersections (since the projection of $A$ to $B$ is $A\cap B$ when $A,B$ are convex subcomplexes with $A\cap B\neq\emptyset$), the bound on the length of $\nest$--chains implies a bound on the number of elements of the hyperclosure that can contain a given $0$--cube, as required by the definition of a factor system in \cite[Sec. 8]{HHS_I}.

  Finally, the DPR follows from Remark~\ref{rem:easy-DPR-cubes} in this example, since each element of the hyperclosure is uniformly coarsely covered by edges coloured yellow in Figure~\ref{fig:counterexample_image}.
  \end{proof}

  Crucially, the orthogonal for non-split property does not hold in $\bar{\mathfrak h}_{/\sim}$. We give here a concrete motivation, which we shall then revisit under a more conceptual light in Remark~\ref{rem:counterexample-discussion}.
  
  \begin{lemma}\label{lem:no_ofns}
      $\bar{\mathfrak h}_{/\sim}$ does not have the orthogonal for non-split domains property.
  \end{lemma}

  \begin{proof}
      Let $F$ be the hyperplane dual to the edge $\Delta$, and let $F'$ be the subcomplex obtained by projecting to $F$ the hyperplane dual to $\Gamma_1$. Notice that $[n]\nest[F]$ for all $n\in \mathbb{N}_{>0}$ and $[F']\nest [F]$. However $F$ is $1$--dimensional, and as a consequence the orthogonal complements \emph{in $F$} of these subcomplex are points and therefore do not belong to $\bar{\mathfrak h}_{/\sim}$. Hence $(\cuco Z,\bar{\mathfrak h}_{/\sim})$ fails to have orthogonals for non-split domains.
  \end{proof}

\subsubsection{Failure of the blow-up construction}
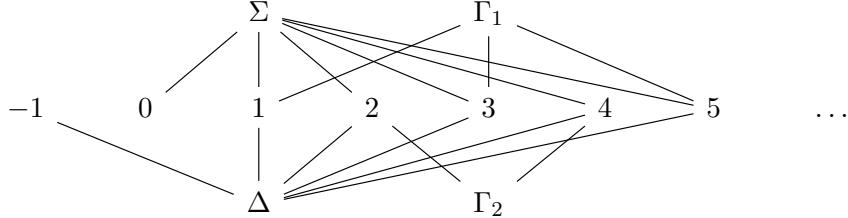
\begin{figure}[htp]
    \centering
    $$\begin{tikzcd}
        &&\Sigma\ar[dl, no head]\ar[d, no head]\ar[dr, no head]\ar[drr, no head]\ar[drrr, no head]\ar[drrrr, no head]&&\Gamma_1\ar[dll, no head]\ar[d, no head]\ar[drr, no head]&&&\\
        -1&0&1&2&3&4&5&\ldots\\
        &&\Delta\ar[ull, no head]\ar[u, no head]\ar[ur, no head]\ar[urr, no head]\ar[urrr, no head]\ar[urrrr, no head]&&\Gamma_2\ar[ul, no head]\ar[ur, no head]&&&
    \end{tikzcd}$$
    \caption{The minimal orthogonality graph $\bar X$.}
    \label{fig:counterexample_min_orth_graph}
\end{figure}

We now prove that the combinatorial HHS $(X,\mathcal W)$ associated to $(\cuco Z,\bar{\mathfrak h}_{/\sim})$ cannot be quasi-isometric to $\cuco Z$. We start by describing the minimal orthogonality graph:
  \begin{lemma}\label{lem:counterexample-min-orth-graph}
  The minimal orthogonality graph of $\bar{\mathfrak h}_{/\sim}$ is the graph shown in Figure~\ref{fig:counterexample_min_orth_graph}.
  \end{lemma}

  \begin{proof}
  Each edge labelled by an integer $n\geq -1$ in Figure~\ref{fig:counterexample_image} is parallel to a (yellow) combinatorial hyperplane, and hence its parallelism class $[h]$ belongs to $\bar{\mathfrak h}_{/\sim}$.  The same is true of the edges $\Gamma_1,\Gamma_2$.  

  Notice moreover that the edge $\Delta$ is (parallel to) the gate from the hyperplane dual to $-1$ to the hyperplane dual to $1$. Hence $[\Delta]\in\bar{\mathfrak h}_{/\sim}$, and the same argument with $0$ replacing $1$ shows that $[\Sigma]\in\bar{\mathfrak h}_{/\sim}$ as well.
  
  Now recall that $\bar{\mathfrak h}$ excludes singletons, so any element of $\bar{\mathfrak h}_{/\sim}$ represented by a subcomplex consisting of a single edge is $\nest$--minimal; this holds in particular for $[\Gamma_1],[\Gamma_2],[\Sigma],[\Delta]$ and $[n],\ n\ge -1$.  The other parallelism classes of edges do not belong to $\bar{\mathfrak h}$.  Moreover, any element of $\bar{\mathfrak h}$ contains a parallel copy of one of the above-named edges.

  Therefore, the minimal orthogonality graph $\bar X$ has vertex set 
  $$\{[\Sigma],[\Delta],[\Gamma_1],[\Gamma_2]\}\cup\{[n]\}_{n\geq -1}.$$
  From the definition of orthogonality, vertices of $\bar X$ are adjacent if and only if the corresponding edges of $\cuco Z$ have parallel copies spanning a square, and thus $\bar X$ is as in Figure~\ref{fig:counterexample_min_orth_graph}.
  \end{proof}

Now recall how to construct the simplicial complex $X$ from $\bar X$: for each $\nest$--minimal $[F]\in\bar{\mathfrak h}_{/\sim}$, the corresponding vertex of $\bar X$ is blown up to a cone over the vertex set of the coordinate space $\fontact [F]$ from the HHS structure of Theorem~\ref{thm:HHS_from_WFS}.  From \cite[Remark 13.2]{HHS_I}, $\fontact [F]$ is the \emph{factored contact graph} of the CAT(0) cube complex $F$.  In the present example, each such $F$ is an edge (one of $\Delta,\Sigma,\Gamma_i,\ i\in\{1,2\},$ or $n,\ n\geq -1$), so $\fontact [F]$ is a single vertex.  Hence $X$ is obtained from $\bar X$ by blowing up each vertex to an edge (the cone over a vertex); each edge of $\bar X$ therefore blows up to a $3$--simplex.

We now argue that this $\bar X$ cannot support a combinatorial HHS structure for $\cuco Z$.

\begin{prop}\label{prop:counterexample}
Suppose that $\mathcal W$ is any $X$--graph such that $(X,\mathcal W)$ is a combinatorial HHS.  Then $\mathcal W$ is bounded, and in particular cannot be quasi-isometric to $\cuco Z$.
\end{prop}

\begin{proof}
For each vertex $v$ of $\bar X$, let $\hat v$ be the $1$--simplex supported on $v$, and let $v, \mathcal Cv$ be the two $0$--simplices supported on $v$.  We use the same notation for a parallelism class of subcomplexes in $\bar{\mathfrak h}$ as for the corresponding vertex of $\bar X$.  So, for instance, in $X$, $[3]$ and $\fontact[3]$ are $0$--simplices, and $\widehat{[3]}$ and $\{[3],\fontact[\Gamma_1]\}$ are $1$--simplices.

If $(X,\mathcal W)$ is a CHHS with $\mathcal W$ unbounded, then by Theorem~\ref{thm:HHS_links}, the augmented links of simplices in $X$ must have arbitrarily large diameter.  However, their diameters are uniformly bounded, in terms of the constant $\delta$ from Definition~\ref{defn:combinatorial_HHS}, which exists by hypothesis.  For example:
\begin{itemize}
    \item $\link_X(\emptyset)^{+\mathcal W}=X^{+\mathcal W}$ has diameter bounded independently of $\delta$, by inspection.

    \item Let $\widehat{[n]}$ be the edge of $X$ projecting to the vertex $[n]$ of $\bar X$.  For $n\geq 1$, the saturation of $\widehat{[n]}$ consists of those $\widehat{[m]}$ with $m=n\ \mathrm{mod}\ 2$, and the link consists of $\widehat\Sigma \cup\widehat \Delta\cup\Gamma_i$ for one of the values of $i$.  For $n\leq 0$, the link is $\widehat\Sigma$ or $\widehat \Delta$. In any case, $\link_X(\widehat{[n]})^{+\mathcal W}$ is connected, by Definition~\ref{defn:combinatorial_HHS}.\eqref{item:cHHS_delta}, and there are only finitely many such graphs, each of which has finitely many vertices, so these links are uniformly bounded.

    \item Note that $\link_X(\widehat\Gamma_1)$ is the union of the edges $\widehat {[2k+1]},\ k\geq 0$, and $\Sat(\widehat{\Gamma_1})=\widehat{\Gamma_1}^{(0)}=\{[\Gamma_1],\fontact[\Gamma_1]\}$. Moreover $Y_{\widehat{\Gamma_1}}=(X-\Sat(\widehat{\Gamma_1}))^{+\mathcal{W}}$ is uniformly bounded, and so must be $\link_X(\widehat{\Gamma_1})^{+\mathcal W}$ which is quasi-isometrically embedded in $Y_{\widehat{\Gamma_1}}$ by Definition~\ref{defn:combinatorial_HHS}.\eqref{item:cHHS_delta}.  
   % In particular, $\widehat\Gamma$ and $\widehat\Sigma$ intersect $X-\Sat(\widehat{\Gamma_1})$ in vertices that are adjacent in $X$ to all vertices in $\link_X(\widehat{\Gamma_1})$.  Hence $\link_X(\widehat{\Gamma_1})^{+\mathcal W}$ is bounded in $Y_{\widehat{\Gamma_1}}$ and thus, by Definition~\ref{defn:combinatorial_HHS}.\eqref{item:cHHS_delta}, it must be bounded.  
   The same holds for $\widehat\Gamma_2$.
\end{itemize}
Similar arguments give boundedness of the remaining links (some are bounded by construction, recall Corollary~\ref{cor:bounded_links}, while the links of almost maximal simplices are points). Hence, $\mathcal W$ is bounded since it is a HHS whose coordinate spaces are all bounded (e.g. by the Distance Formula~\ref{thm:distance_formula}).
\end{proof}

\begin{rem}[What's wrong with this example?]\label{rem:counterexample-discussion}
     The HHS structure on $\cuco Z$ has wedges, clean containers and the dense product region property by Theorem~\ref{thm:hyperclosure_is_orthogonal} and Lemma~\ref{lem:counterexample-factor-system-sleazy}.  Moreover, orthogonality determines nesting, i.e. $[F]\propnest[F']$ if and only if $[F']^\orth\propnest[F]^\orth$.

     This illustrates that an HHS in which orthogonality determines nesting need not have the orthogonals for non-split domains property.  This is a typical phenomenon in CAT(0) cube complexes.  Indeed, if $[F]\in\bar{\mathfrak h}_{/\sim}$ (in an arbitrary CAT(0) cube complex), then as long as $[F]$ is not the unique $\nest$--maximal element, $[F]^\orth$ is defined and belongs to $\bar{\mathfrak h}_{/\sim}$.  However, the existence of orthogonals is not inherited by the sub-HHS structure on $F$, in general.  More precisely, we could consider the hyperclosure of the CAT(0) cube complex $F$, called $\bar{\mathfrak h}^F$. Since hyperplanes of $F$ are of the form $F\cap H$, where $H$ is a hyperplane of $\cuco Z$ intersecting $F$, by e.g. \cite[Lem. 3.1]{CapraceSageev},  and $F\cap H=\gate_F(H)$ by e.g. \cite[Lem. 2.6]{HHS_I}, we see that $\bar{\mathfrak h}^F$ naturally embeds in $\bar{\mathfrak h}$ (preserving parallelism), yielding a set of domains in $\bar{\mathfrak h}_{/\sim}$ nested in $[F]$.  However, the set of all $[F']\in \bar{\mathfrak h}_{/\sim}$ which are nested inside $[F]$ is in general larger than this, since it contains subcomplexes of the form $\gate_F(F'')$, where $F''$ is an element of the hyperclosure not parallel to any subcomplex in the “intrinsic” hyperclosure of $F$. Hence, given an element $[F']\nest[F]$, its orthogonal inside $[F]$ might not exist.
     
     Concretely, in our counterexample, we can translate Lemma~\ref{lem:no_ofns} in the above language as follows. Let $F$ be the hyperplane dual to the edge $\Delta$. Since $F$ is $1$--dimensional, its hyperplanes are points and so its intrinsic hyperclosure is empty.  But in $\bar{\mathfrak h}_{/\sim}$, we have, say, $[n]\nest[F]$ and $[F']\nest [F]$, where $[F']$ is the subcomplex obtained by projecting to $F$ the hyperplane dual to $\Gamma_1$.  Again since $F$ is $1$--dimensional, these subcomplexes do not have non-trivial orthogonal complements that belong to the hyperclosure but are nested in $F$.  Hence $(\cuco Z,\bar{\mathfrak h}_{/\sim})$ fails to have orthogonals for non-split domains.

     \textbf{Absent orthogonals and boundedness of $\mathcal W$:}  Let us see how the failure of orthogonals for non-split domains caused problems in the example.  Recall that the simplex $\widehat{\Gamma_1}$ had bounded link \emph{because} its saturation failed to contain, say, $\widehat\Sigma$, which then acted as a cone-point in $X-\Sat(\widehat\Gamma_1)$ over $\link_X(\widehat\Gamma_1)$.  Consider the hyperplanes $H_{\Gamma_1}$ and $H_\Sigma$ in $\cuco Z$ dual to the edges $\Gamma_1$ and $\Sigma$, and let $F$ be the projection of the former onto the latter, so $F\in\bar{\mathfrak h}$ and $F$ is the ray consisting of the red edges in Figure~\ref{fig:counterexample_image}.  Then $F$ is non-split, since any two $\nest$-minimal domains it contains (i.e. every two red edges) are transverse. Moreover $F\nest H_\Sigma$, but $F^\orth=\Gamma_1\cup\Sigma\cup\Delta$, so $F$ is not orthogonal to anything properly nested in $H_\Sigma$.  Back in $X$, the role of $\fontact [F]$ (which is a ray) should be played by the link of $\widehat\Gamma_1$, which is contained in the link of $\widehat\Sigma$.  So it seems reasonable that by adding the “missing” orthogonal domain, one could add a vertex $w$ to $X$ in such a way that $w\star\Sigma$ is defined and has the same link as $\widehat\Gamma_1$. This way, removing $\Sat(\widehat\Gamma_1)$ would now remove $\Sigma$, which as we mentioned is an obstruction to having an unbounded augmented link.  Of course, one would also need to deal with $\Delta$ similarly.  But then, if $\widehat\Gamma_1$ and $\widehat\Gamma_2$ were made to have unbounded links, a new difficulty would arise, because Definition~\ref{defn:combinatorial_HHS}.\eqref{item:cHHS_join} would then demand that $\widehat \Sigma$ and $\widehat \Delta$ have a common nested simplex whose link contains those of $\widehat\Gamma_1$ and $\widehat\Gamma_2$.  This would presumably be addressed by adding a single $w$ as above, joined to both $\Sigma$ and $\Delta$ and $\Gamma_1$, and corresponding to an orthogonal complement of $F$ inside the wedge of $H_\Sigma$ and $H_\Delta$.
\end{rem}

\begin{rem}[Connection to lattice theory]\label{rem:lattice}
Let $(\mathfrak S,\nest,\orth)$ be an orthogonal set and let $\emptyset$ be a symbol distinct from all elements of $\mathfrak S$. Then the operation $\wedge$ which, by Definition~\ref{defn:orth_set}, was partially defined on $\mathfrak S$, extends to a binary operation on $\mathfrak S\sqcup\{\emptyset\}$ if one sets $U\wedge \emptyset=\emptyset\wedge U=\emptyset$ for $U\in\mathfrak S\sqcup\{\emptyset\}$, and $U\wedge V=\emptyset$ if $U,V\in\mathfrak S$ do not have any common nested elements of $\mathfrak S$.  We extend $\nest$ so that $\emptyset$ is the unique $\nest$--minimal element.  Assuming that the complexity is finite, we also have a join operation: $U\vee V$ is the unique $\nest$--minimal $W$ such that $U,V\nest W$; that this is well-defined is an easy exercise using wedges and finite complexity.  In fact, the poset $(\mathfrak S\sqcup\{\emptyset\},\nest)$, equipped with the operations $\wedge$ and $\vee$, is a \emph{complete lattice}, as in e.g. \cite[Definition 4.1]{Sanpei}, where completeness also follows from finite complexity.  Moreover, the clean container assumption, together with Lemma~\ref{lem:orth_implies_nest_equiv}, gives an involution ${}^\orth:\mathfrak S-\{S\}\to\mathfrak S-\{S\}$, where $S$ is the unique $\nest$--maximal element, and we extend this to $\mathfrak S\cup\{\emptyset\}$ by declaring $S^\orth=\emptyset$ and $\emptyset^\orth= S$.\footnote{On $\mathfrak S$, the orthogonality \emph{relation} $U\orth V$ is still equivalent to $U\nest V^\orth$, i.e. $U\wedge V^\orth=U$.  However, $\nest$ and $\orth$ are not mutually exclusive on $\mathfrak S\cup\{\emptyset\}$, but the only failure is $\emptyset\nest S$ and $\emptyset\orth S$.}  This makes $\mathfrak S\sqcup \{\emptyset\}$, with the lattice relations $\wedge,\vee$ and the orthogonal complementation operator ${}^\orth$, an \emph{ortholattice}, defined in e.g. \cite[Sec. 1.2]{Stern}.

Strong orthogonality (Property~\ref{property:strong_orth}) then becomes: for all $U,V\in\mathfrak S$ such that $U\propnest V$, we have that $U^\orth \wedge V\neq \emptyset$.  Another formulation is: for any $V\in\mathfrak S$, the order ideal $\{U\in\mathfrak S:U\nest V\}$ is again an orthogonal set, with the involution $U^\orth_V=U^\orth \wedge V$.

Now, if $U\propnest V$, then $(U^\orth\wedge V)\vee U\nest V$.  If for some $W\propnest V$ we have $U\nest W$ and $U^\orth\wedge V\nest W$, then strong orthogonality provides $A:=W^\orth \wedge V\neq\emptyset$.  But then $(U^\orth \wedge V)\orth A$, since $A\orth W$, while on the other hand $W^\orth \wedge V=A\nest U^\orth\wedge V$ since $U\nest W$.  This is a contradiction.  We have showed that strong orthogonality implies the identity $(U^\orth \wedge V)\vee U=V$ whenever $U\nest V$.  Ortholattices satisfying this identity have a name: they are \emph{orthomodular} \cite{Stern}.

From the point of view of HHS structures, factor systems in CAT(0) cube complexes provide the main motivating examples of orthogonal sets, and we saw earlier that orthomodularity fails in general; in fact, one can already see this failure in HHS structures on, say, right-angled Artin groups (see \cite[Sec. 8]{HHS_I}).  This raises the following:

\begin{quest}\label{question:lattice}
Let $(\mathfrak L,\wedge,\vee,{}^\orth,\emptyset,S)$ be an ortholattice such that $\nest$--chains in $\mathfrak L$ have length at most $N<\infty$, where $U\nest V$ means $U\wedge V=U$.  Write $U\orth V$ to mean $U\nest V^\orth$. 
\begin{itemize}
    \item Does there exist an \textbf{orthomodular} ortholattice $\mathfrak L_1$ and an injective map $\mathfrak L\to\mathfrak L_1$ that preserves the relations $\nest$ and $\orth$, as well as the negations of those relations?  Under what conditions can $\mathfrak L_1$ be chosen so that chains in $\mathfrak L_1$ also have length at most $N$?

    \item If a group $G$ acts on $\mathfrak L$ cofinitely, preserving the relations $\nest,\orth,\not\nest,\not\perp$, when can $\mathfrak L_1$ be chosen as above so that the $G$--action extends to $\mathfrak L_1$ and $|{}_G\backslash{}^{\mathfrak L_1}|<\infty$?
\end{itemize}
\end{quest}
The goal would be to begin with an HHS/G whose index set is an orthogonal set (e.g. a compact special group) and produce a new HHS/G structure to which Theorem~\ref{thm:strong_iff} applies.  Answers to the above questions are not quite sufficient but appear necessary, and also of independent interest.  We suspect that a sufficient condition for constructing $\mathfrak L_1$ will involve the existence of an order-preserving, $\orth$--preserving map from $\mathfrak L$ to a finite boolean lattice. It is also possible that this sort of construction is known to lattice theorists, in which case we would be grateful for a reference.
\end{rem}

\bibliography{biblio}
\bibliographystyle{alpha}
\end{document}